\documentclass[moor, nonblindrev]{informs1}

\usepackage{thmtools,thm-restate}
\usepackage{natbib}

\usepackage{tikz}

\usepackage{caption}
\usepackage{subcaption}
 \NatBibNumeric
 \def\bibfont{\small}%
 \bibpunct[, ]{[}{]}{,}{n}{}{,}%
  \usepackage{bm}
\usepackage{amsmath}
    \usepackage{amssymb}
\usepackage{tikz}
\usetikzlibrary{external}
\usepackage{enumitem}
\usepackage{pgfplots}

\newcommand{\Real}{\mathbb{R}}
\definecolor{darkred}{rgb}{0.7,0,0}
\definecolor{lightred}{rgb}{0.5,0,0}
\definecolor{lightgreen}{rgb}{0,0.8,0}
\definecolor{mathset}{HTML}{CFD9E8}
\definecolor{mathboundary}{HTML}{5E81B5}

\newcommand{\blue}[1]{#1}
\newcommand{\green}[1]{{\green{blue}#1}}
\usepgfplotslibrary{fillbetween}
\pgfplotsset{compat=1.10}
\usetikzlibrary{arrows,fadings,pgfplots.patchplots,decorations.pathmorphing,backgrounds,fit,positioning,shapes.symbols,chains}
\tikzstyle{frac}=[rectangle,draw=black,fill=lightgreen,
inner sep=0pt,minimum size=0.1cm]
\tikzstyle{active}=[draw=darkred,very thick]
\tikzstyle{inactive}=[draw=lightred, dotted,very thick]
\tikzstyle{notactive}=[draw=mathboundary,thick, fill=mathset]
\tikzstyle{notactiveinner}=[draw=mathboundary,thick,fill=white]
\usepackage{url}
\usepackage{color}
\newcommand{\set}[1]{\left\{#1\right\}}                     
\newcommand{\norm}[1]{\left\|#1\right\|}
\newcommand{\abs}[1]{\left|#1\right|}                    
\newcommand{\bra}[1]{\left(#1\right)}
\usepackage{stmaryrd}
\newcommand{\sidx}[1]{\left\llbracket     #1 \right\rrbracket}
\usepackage{comment}
\DeclareMathOperator{\cl}{cl}
\DeclareMathOperator{\conv}{conv}

\DeclareMathOperator{\relint}{relint}

\DeclareMathOperator{\aff}{aff}
\DeclareMathOperator{\spann}{span}
\DeclareMathOperator{\Int}{int}

\DeclareMathOperator{\proj}{proj}
\definecolor{darkred}{rgb}{0.7,0,0}
\definecolor{lightred}{rgb}{0.5,0,0}
\definecolor{lightgreen}{rgb}{0,0.8,0}

\newcommand{\bw}[0]{\bm{w}}

\newcommand{\bx}[0]{\bm{x}}
\newcommand{\by}[0]{\bm{y}}
\newcommand{\bz}[0]{\bm{z}}

\newcommand{\pointtoproof}[1]{$\longrightarrow$ See page \pageref{#1} for the proof.}

\usepackage[colorlinks=true,breaklinks=true,bookmarks=true,urlcolor=blue,
     citecolor=blue,linkcolor=blue,bookmarksopen=false,draft=false]{hyperref}

\usepackage[normalem]{ulem}
\usepackage{cancel}

\def\EMAIL#1{\href{mailto:#1}{#1}}

\TheoremsNumberedBySection  

\EquationsNumberedBySection 

\begin{document}


 \RUNAUTHOR{Lubin, Vielma, and Zadik}

 \RUNTITLE{Mixed-integer convex representability}

\TITLE{Mixed-integer convex representability\footnote{An earlier version of this work appeared in the proceedings of IPCO 2017~\cite{Lubin2016b}.}}
\ARTICLEAUTHORS{%
\AUTHOR{Miles Lubin}
\AFF{Google Research, \EMAIL{mlubin@google.com}}
\AUTHOR{Juan Pablo Vielma}
\AFF{Sloan School of Management, M.I.T., \EMAIL{jvielma@mit.edu}}
\AFF{Google Research, \EMAIL{jvielma@google.com}}
\AUTHOR{Ilias Zadik}
\AFF{Operations Research Center, M.I.T., \EMAIL{izadik@mit.edu}}
\AFF{Center for Data Science, N.Y.U., \EMAIL{zadik@nyu.edu}}
} 

\ABSTRACT{%
Motivated by recent advances in solution methods for mixed-integer convex optimization (MICP), we study the fundamental and open question of which sets can be represented exactly as feasible regions of MICP problems. We establish several results in this direction, including the first complete characterization for the mixed-binary case and a simple necessary condition for the general case. We use the latter to derive the first non-representability results for various non-convex sets such as the set of rank-1 matrices and the set of prime numbers. Finally, in correspondence with the seminal work on mixed-integer linear representability by Jeroslow and Lowe, we study the representability question under rationality assumptions. Under these rationality assumptions,  we establish that representable sets obey strong regularity properties such as periodicity, and we provide a complete characterization of representable subsets of the natural numbers and of representable compact sets. Interestingly, in the case of subsets of natural numbers, our results provide a clear separation between the mathematical modeling power of mixed-integer linear and mixed-integer convex optimization. In the case of compact sets, our results imply that using unbounded integer variables is necessary only for modeling unbounded sets.

}%


\KEYWORDS{}
\MSCCLASS{90C11,90C25}
\ORMSCLASS{Primary: Programming:Integer:Theory  ; Secondary: Mathematics:Convexity }
\HISTORY{}

\maketitle

%


\section{Introduction}\label{intro:sec}

More than 60 years of development have made mixed-integer linear programming (MILP) an extremely successful tool \cite{Juenger2010}. MILP's modeling flexibility allows it to describe a wide range of business, engineering and scientific problems, and, while solving MILP is NP-complete\footnote{For MILP to be in \emph{bit complexity} class NP in addition to being NP-hard, we need some assumptions on the encoding of the MILP instance \cite[Chapters 2 and 17]{schrijver1998theory}.}, many of these problems are routinely solved in practice thanks to state-of-the-art solvers that nearly double their machine-independent speeds every year  \cite{Achterberg2013,bixby2012brief}. The last decade has seen a surge of activity on the solution and application of mixed-integer convex programming (MICP), which extends MILP's versatility by allowing the use of convex constraints in addition to linear inequalities. A small selection of applications of MICP includes sensor placement in dynamic networks \cite{nugroho2019sensor}, configuration of antennas in wireless communications \cite{Ni2018signalprocessing}, and grasp planning in robotics \cite{Liu2020grasp}; see also \cite{benson2013mixed} for a review of applications of the MICP sub-class known as mixed-integer second-order cone programming. From a theoretical complexity standpoint, solving MICP problems can be much harder than solving MILP problems---even simple convex quadratic mixed integer programming problems can fail to be in NP (e.g. \cite[Example 2, page 227]{del2017mixed})\footnote{This failure can hold under instance-encoding restrictions that are similar to those required for MILP to be in NP \cite{del2017mixed,koppe2012complexity}.}. Nonetheless, state-of-the-art solvers for MICP are rapidly closing on the effectiveness and robustness of their MILP counterparts \cite{Filmint,Bonami2008,coey2018outer,DuranGrossmann,Gally2016,DBLP:conf/ipco/LubinYBV16,Lubin2016a}.

Motivated by these developments, in this work we investigate the \textit{representability} power of MICP; that is, we study and classify which sets can be represented as the feasible regions of MICP problems. In the mid-1980s, Jerowslow and Lowe were similarly motivated by advancements in solution methods for MILP when they developed their iconic result on what sets can be modeled using \emph{rational} MILP formulations, i.e., MILP formulations defined by linear inequalities with rational coefficients (e.g. \cite[Section 11]{Vielma2015}). More recent results for mixed-integer programming (MIP) representability have included alternative algebraic characterizations of rational MILP representability \cite{BasuMILPRMOR} and extensions to  special cases of MICP representability~\cite{DBLP:conf/ipco/PiaP16,DelPia2016,DelPia2017,Dey2013,hemmecke2007representation,Lubin2016b,vielma2017small}. As surveyed below, many of the first natural questions that can be asked on general MICP representability remain open, and hence are a next step in this line of research.

These questions have both academic and practical implications, because knowing that a non-convex set is representable using MICP could open a new path for modeling it using MICP techniques. Negative answers, on the other hand, imply that some approximation or redefinition of the problem is needed to apply MICP techniques.
For example, the results of Jeroslow and Lowe \cite{Jeroslow1984} show that the set \blue{
$\{ \bx \in \mathbb{R}^2 : x_2 \ge 0, x_1 = 0 \} \cup \{ \bx \in \mathbb{R}^2 :
x_2 \ge 1, x_1 \ge 0 \}$}, a union of two polyhedra that models the epigraph of a
fixed cost in production, is not rational MILP representable. However, it is now common knowledge that we \textit{can} model such fixed costs by introducing an upper bound on $x$, i.e., by intersecting with $\mathbb{R} \times [0, M]$ for
some $M \in \mathbb{R}$. The study of MICP representability has the potential to yield similar
modeling insights.

In the following subsections, we outline existing results on MICP representability and identify open questions that frame the contributions of this work. To do this, we first formally define representability.

\begin{definition}\label{MICPFORMULATIONDEF}
   Let $n,p,d \in \mathbb{N}$, $S \subseteq \mathbb{R}^n$ and $M\subseteq \mathbb{R}^{n+p+d}$ be a closed convex set. We say \blue{that the pair $\bra{M,d}$\footnote{The value of $p$ is implied by $d$ and the number of variables used to describe  $S$ and $M$ ($n$ and $n+p+d$, respectively).}} induces an \textbf{MICP formulation} of $S$ if
   \begin{equation}\label{MICPdefexists}
   \bx \in S \quad \Leftrightarrow \quad \exists \by \in \mathbb{R}^{p}, \bz \in \mathbb{Z}^{d}  \quad \text{s.t.} \quad \bra{\bx,\by,\bz}\in  M.
   \end{equation}
   \blue{If $p=0$, we say that the formulation induced by $\bra{M,d}$ is \textbf{pure}. In addition, if $d$ is evident from the context, we say more concisely that $M$ induces the formulation. }
\end{definition}

 Note that the \emph{auxiliary variables} in \autoref{MICPFORMULATIONDEF} (i.e. those that are not the original $\bx$ variables) include both integer auxiliary variables $\bz$ and continuous auxiliary variables $\by$. Both classes of variables add to the modeling power of MICP formulations. \blue{In particular, the use of continuous auxiliary variables $\by$ allows for the modeling of certain non-closed sets despite the restriction for $M$ to be closed. For instance,  consider the closed convex set $M=\set{(x,y) \in \Real^2\,:\, x,y\geq 0,\, x y\geq 1}$ which induces a formulation of the open set $S=(0,\infty)=\set{x\in \Real\,:\, \exists y\in \Real \text{ s.t. } (x,y)\in M}$ without the use of integer variables. Furthermore, as we will see in Section~\ref{sec:bounded_intro}, the use of continuous auxiliary variables provides additional modeling power even if we consider only the case where the set $S$ is closed.}

We now introduce the followings notions of MICP representability, which comprise the core objects of study of the present work.
\begin{definition}
   We call a set $S \subseteq \mathbb{R}^n$ \textbf{MICP representable} (MICP-R) if there exists a closed convex set $M\subseteq \mathbb{R}^{n+p+d}$ that induces an MICP formulation of $S$.
\end{definition}
\begin{definition}

   We call a set $S$ \textbf{(rational) MILP representable} (MILP-R) if there exists  a (rational) polyhedron $M\subseteq \mathbb{R}^{n+p+d}$ that induces an MICP formulation of $S$.
\end{definition}
\begin{definition}\label{binarydef}
We call a set $S$  \textbf{binary MICP-R (MILP-R)}\footnote{\label{binaryboundedfootnote}Binary MILP-R sets were also denoted bounded MILP-R sets by Jeroslow and Lowe \cite{Jeroslow1984} because any formulation with integer variables that are bounded can be transformed to one with only binary variables through standard transformations. For example, the constraint $z\in [a,b]\cap \mathbb{Z}$ can be written as $z=\sum_{i=a}^b i \tilde{z}_i,\quad \sum_{i=a}^b \tilde{z}_i=1,\quad \tilde{z}_i\in \set{0,1}\quad \forall i\in [a,b]\cap \mathbb{Z}$. We avoid this notation as it can erroneously suggest that binary/bounded MILP- or MICP-R sets must be bounded (i.e., as sets).} if  there exists a closed convex set (polyhedron) $M\subseteq \mathbb{R}^{n+p+d}$ that induces an MICP formulation of $S$ and satisfies \[\proj_{\bz}\bra{M\cap \bra{\mathbb{R}^{n+p} \times \mathbb{Z}^d}}\subseteq\set{0,1}^d,\]where $\proj_{\bz}$ is the projection onto the last $d$ variables.
\end{definition}
\begin{definition}
   \blue{We call a set $S \subseteq \mathbb{R}^n$ \textbf{pure MICP representable (MICP-R)} if there exists a closed convex set $M\subseteq \mathbb{R}^{n+d}$ which induces a pure MICP formulation of $S$. We similarly define \textbf{pure (rational) MILP-R} and \textbf{pure binary MICP-R (MILP-R)}. }
\end{definition}

\subsection{Binary MICP representability}

The classical characterization by Jeroslow and Lowe  \cite{Jeroslow1984} established that  binary rational MILP-R  sets are exactly those sets that are unions of finitely many rational polyhedra that share the same recession cone\footnote{See \autoref{rationaldefcone} and \autoref{rationaldefcprop} for a definition of recession cone.}. Note that a potentially counter-intuitive implication of this result is that the simple set \blue{$S = (-\infty, -1] \cup \set{0}$} is not binary rational MILP-R. A non-polyhedral generalization of this result has been proven recently by Del Pia and Poskin  \cite{DBLP:conf/ipco/PiaP16} under rational ellipsoidal restrictions on $M$\footnote{The restriction requires $M$ to be the intersection of an ellipsoidal cylinder having a rational recession cone and a rational polyhedron.}. Under these rational ellipsoidal restrictions on $M$, binary MICP-R sets are finite unions of rational ellipsoidal sets that also must share the same recession cone.

It is known that a finite union of closed convex sets that share a common recession cone is binary MICP-R \cite{Ceria1999,jeroslow1987representability,stubbs1999branch}\footnote{For a self contained proposition summarizing these results, see \cite[Theorem~1]{vielma2017small}.}. However, whether the common recession cone condition is necessary for a set to be binary MICP-R has remained an open problem. For example, it was not known prior to the present paper whether the set \blue{$S = (-\infty, -1] \cup \set{0}$} \blue{ or the  epigraph of a
fixed cost in production given by
$S'=\{ \bx \in \mathbb{R}^2 : x_2 \ge 0, x_1 = 0 \} \cup \{ \bx \in \mathbb{R}^2 :
x_2 \ge 1, x_1 \ge 0 \}$ are} binary MICP-R or not. \blue{The closest result along these lines appears in \cite{jeroslow1987representability}, where Jeroslow showed that a set is \textit{pure} binary MICP-R if and only if it is a finite union of closed convex sets that share a common recession cone (see \autoref{lem:balasform_jeroslow} and \autoref{prop:jeroslowbinaryequiv}). Hence, neither of the example sets $S$ or $S'$ is pure binary MICP-R.}

\subsection{\blue{Obstructions to MICP representability}}

By enumerating over the possible values of $\bz$ in  \eqref{MICPdefexists}, we can check that any MICP-R set can be written as a countable union of projections of closed convex sets.
A very rich family of non-convex sets can be expressed as a countably infinite union
of projections of closed convex sets (e.g. \autoref{dcrepref} will show this for sets which are complements of convex bodies). Hence, it would be surprising if arbitrary countably infinite unions of closed convex sets are MICP-R; instead,
we expect the MICP-R sets to have fairly regular structures.

\blue{Despite these expectations, the simple-to-state question \begin{center} \textit{``Is a set $S$ MICP-R?"} \end{center} has not yet been studied, to the best of our knowledge. For example, the MICP-R status of the set of rank-1 matrices or even the set of integer points on the graph of $f(x)=x^2$ remains unknown. In this work, we make a first attempt to answer this question.}

\blue{Lacking an effective complete characterization of MICP-R sets, we take the route of studying \textit{obstructions} to representability; these are conditions that imply a set is not MICP-R. While obstructions to MICP representability are novel, they have been used in the context of polyhedral and positive semidefinite representability~\cite{Fawzi2020}.}

\subsection{MICP representability: a non-trivial ``rationality" question }\label{sec:rational_intro}

\blue{In the special cases of (non-binary) MICP-R studied so far, where $M$ is polyhedral~\cite{Jeroslow1984} or ellipsoidal~\cite{DelPia2016}, authors have obtained exact characterizations with nice structural properties only in the presence of rationality assumptions. Recall the characterization of Jeroslow and Lowe~\cite{Jeroslow1984} of the rational MILP-R case.}
\begin{definition}
For any finite set  $U\subseteq \mathbb{Z}^n$ we denote the \textbf{integer cone of $\bm{U}$} as the set $\operatorname{intcone}\bra{U}=\set{ \sum_{\bm{u}\in U} \lambda_u \bm{u} : \lambda \in \mathbb{Z}_+^U}$.
\end{definition}
\begin{theorem}[\cite{Jeroslow1984}]\label{JLTheo}
   A set $S\subseteq \Real^n$ is rational MILP-R if and only if there exist a finite set  $U\subseteq \mathbb{Z}^n$ and rational polytopes $\set{P_i}_{i=1}^k$ such that $P_i\subseteq \Real^n$ for all $i\in \sidx{k}$ and
   \begin{equation}\label{eq:jlregular}
   S= \bigcup_{i=1}^k P_i + \operatorname{intcone}\bra{U}.
   \end{equation}
   \end{theorem}
   \blue{It is not hard to show (e.g. see \autoref{boundedorperiodic:coro}) that any set that satisfies \eqref{eq:jlregular}  is either bounded (if $U\setminus \set{\bm{0}} =\emptyset$) or periodic according to the following standard definition.}
    \begin{definition}\label{periodicdef}
 A set $S$ is {\bf periodic} if there exists $ \bm{u} \in \Real^n\setminus \set{\bm{0}}$ such that   $\bm{x} + \lambda \bm{u} \in S$ for all $\bm{x} \in S$ and $\lambda \in \mathbb{N}$.
   \end{definition}
\blue{Similar results hold under a rational ellipsoidal restriction on the set $M$~\cite{DelPia2016} and for other restricted versions of MICP-R (e.g. see \autoref{boundedprop} and \autoref{boundedorperiodic:coro}).}

\blue{Why is rationality important for this characterization?}  Standard textbook exercises
(e.g. \cite[Exercise 4.30]{conforti2014integer}) show that \blue{\autoref{JLTheo} may not hold for sets with an MILP formulation induced by a polyhedron $M$ whose recession cone is not a \textit{rational} polyhedron.}  For instance, consider the  set
\blue{
\begin{equation}\label{badset2}
S=\set{x\in \mathbb{Z}: f\bra{x\sqrt{2} }\notin \bra{0.4, 1-0.4} }
\end{equation}}
where  $f(x)=x-\lfloor x \rfloor$ denotes the fractional part of $x$.
\blue{Using Kronecker's approximation theorem (\autoref{Kronecker}) we can show that the set $S$ is a countably infinite set of integers that is \emph{highly} \textit{non-periodic} in that: (1) it is not periodic, and (2) none of its subsets is periodic (\autoref{exampleslemmaperiodicref:lemma}).
}

\blue{Because $S$ is unbounded and not periodic, \autoref{JLTheo} (e.g. through \autoref{boundedorperiodic:coro}) implies that it cannot be rational MILP-R. Nevertheless, $S$ \textit{does} have} the simple MILP formulation induced by the (non-rational) polyhedron
\blue{
\begin{equation}\label{irrationalex2}
    M=\set{\bra{x,\bz}\in \mathbb{R}\times \Real^2\,:\, z_1=x,\; -0.4\leq \sqrt{2} z_1 - z_2 \leq     0.4}\footnote{For validity of this formulation, let $\bra{x,\bz}\in M\cap\bra{\mathbb{R}\cap \mathbb{Z}^2}$. If $0\leq \sqrt{2} z_1 - z_2 \leq     0.4$, then $z_2=\lfloor \sqrt{2} z_1\rfloor$ and $f\bra{z_1\sqrt{2} }\leq 0.4$. If $0\leq  z_2 -\sqrt{2} z_1 \leq     0.4$, then $z_2=\lceil \sqrt{2} z_1\rceil=\lfloor \sqrt{2} z_1\rfloor+1$ and $f\bra{z_1\sqrt{2} }\geq 1-0.4$. The result follows from $x=z_1$. }.
    \end{equation}
}
Using the irrational coefficient $\sqrt{2}$ is crucial to construct such a formulation, so such pathological cases are \blue{excluded} in MILP formulations by requiring rational coefficients defining the underlying polyhedron. We infer from this MILP example that any hope of obtaining a \blue{result like \autoref{JLTheo}} for general MICP-R is dependent on understanding how to place rationality restrictions on convex sets.

Developing such rationality restrictions for general MICP formulations is a challenge that is \textit{previously largely unaddressed}. \blue{Natural first attempts might consider} restricting to MICP formulations described by polynomial inequalities with rational coefficients, or \blue{to} \textit{conic formulations} \cite{Ben-Tal2001,Boyd2004,MOSEKModeling} of the form $\bm{A}\bm{x}+\bm{B}\bm{y}+\bm{C}\bm{z}-\bm{f}\in K$, where $\bm{A}$, $\bm{B}$, and $\bm{C}$ are appropriately sized rational matrices, $\bm{f}$ is a rational vector, and $K$ is a specially structured closed convex cone. This unfortunately does not solve the problem; we can show (e.g. see \autoref{nonperiodicirrationalex:lem}) that the set $M$  described by  \eqref{irrationalex2} is the projection onto the $x$ and $\bz$ variables of the set described by
\blue{
\begin{subequations}\label{irrationalexalt2}
\begin{alignat}{5}
z_1&=x,\quad& -0.4&\leq y_5 - z_2 \leq     0.4,\quad& \norm{\bra{y_1,y_1}}_2&\leq y_2-y_5,\quad\\
 \norm{\bra{y_2,y_2}}_2&\leq  2 (z_1-y_1),\quad&
  \norm{\bra{y_3,y_3}}_2&\leq y_5-y_4 ,\quad& \norm{\bra{y_4,y_4}}_2&\leq 2 (y_3-z_1)
\end{alignat}
\end{subequations}}
\blue{Restricting to polynomial or conic constraints with rational coefficients is therefore not restrictive enough to exclude the non-periodic set \eqref{badset2}}\footnote{Let $\mathcal{L}^{1+n} = \{ (t_0,\bm{t}) \in \mathbb{R}^{1+n} : ||\bm{t}||_2 \le t_0\}$ be the second-order or Lorentz cone \cite{alizadeh2003second}. Then formulation
\eqref{irrationalexalt2} can be written as a conic formulation using a cone $K$ which is the Cartesian product of Lorentz cones (with the identification $\mathcal{L}^{1+0} =\Real_+$).
A version with polynomial inequalities follows by noting that $\mathcal{L}^{1+n} = \{ (t_0,\bm{t}) \in \mathbb{R}^{1+n} : ||\bm{t}||_2^2 \le t_0^2,\quad t_0\geq 0\}$.}. The most technical developments of this work are made to obtain a sufficiently restrictive, yet also sufficiently expressive, rationality restriction \blue{under which one can obtain structural results on MICP-R sets}.

\paragraph{Organization} The remainder of the paper is structured as follows. In Section~\ref{sec:contributions},
we outline the contributions of this work. In Section~\ref{sec:definitions}, we
introduce our notation and some background material on convex analysis. In Section~\ref{sec:main_results}, we formally state and prove our results concerning characterizations of MICP-R sets and \blue{obstructions to} MICP representability.
In Section~\ref{MICPrat}, we formally state our results concerning rationality and periodicity. Section~\ref{MICPrat}  includes proofs of all these results except for  \autoref{periodictheoref} whose proof we postpone to Sections~\ref{secondproofsection}  because of its highly technical nature. Finally, in Section~\ref{conclusion:sec}, we present some final remarks.  Appendices~\ref{app:conicratpolyirrat}--\ref{rationalnotaffine} include additional complementary materials. Any deferred proofs are marked with a $\longrightarrow$ symbol pointing to the page where the proof can be found.

\section{Contributions}\label{sec:contributions}
\subsection{Binary MICP representability}\label{sec:BinaryMICPcontribution}

\begin{figure}[htpb]\centering
 \begin{subfigure}{.49\linewidth}
 \centering
  \includegraphics{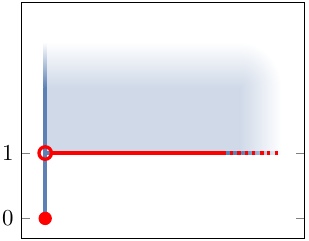}

\captionsetup{font=small}
\caption[size=small]{In red, the graph of $f(x) = \{ 0 \text{ if } x = 0, 1 \text{ if } x > 0 \}$. In blue, its epigraph~\eqref{eq:fixedcostprod}.  \label{fig:fixedcost}}
\end{subfigure}
         \begin{subfigure}{.49\linewidth}
         \centering
         \includegraphics{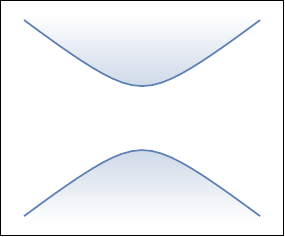}
%
%
%
%
\captionsetup{font=small}
\caption[size=small]{In blue, the two-sheet hyperboloid~\eqref{eq:hyperboloid}. \vspace{\baselineskip}\label{fig:hyperboloid}}
\end{subfigure}
\caption{\blue{Sets not previously known to be binary MICP-R.}}
\end{figure}

As noted in the introduction, \blue{a set is pure binary MICP-R if and only if it is a finite union of closed convex sets that share the same recession cone \cite{jeroslow1987representability}. In addition, all known non-pure} binary MICP formulations for finite unions of closed convex sets also require the convex sets to share the same recession cone. We show that the identical recession cone condition is \textit{not} necessary for binary MICP representability, and provide the first complete characterization of binary MICP-R. More specifically, we show that \textit{a set is binary MICP-R if and only if it is the finite union of convex sets, where the convex sets are
a projection of \blue{arbitrary} closed convex sets}. This result is
constructive, i.e., we provide an explicit procedure to construct an MICP formulation for these sets. \blue{In particular, this result positively answers the question on the binary MICP representability of $S = (-\infty, -1] \cup \set{0}$ and the  epigraph of a
fixed cost in production given by
\begin{equation}\label{eq:fixedcostprod}
S'=\{ \bx \in \mathbb{R}^2 : x_2 \ge 0, x_1 = 0 \} \cup \{ \bx \in \mathbb{R}^2 :
x_2 \ge 1, x_1 \ge 0 \}
\end{equation}and depicted in \autoref{fig:fixedcost}.
The result also gives the first MICP formulation for finite unions of non-polyhedral closed convex sets with different recession cones such as the} two sheet hyperboloid (see \autoref{fig:hyperboloid} for the case $n=2$),
\begin{equation}\label{eq:hyperboloid}
    S''=\set{\bx\in\mathbb{R}^n\,:\, \sum_{i=1}^{n-1} x^2_i +1 \leq x_n^2 }.
\end{equation}

\subsection{\blue{Obstructions} to MICP representability}

We develop a powerful geometric \blue{obstruction} to general MICP representability, which we use to prove the first non-representability results. The property is a natural extension of the elementary characterization of closed convex sets as those sets that are closed under taking the midpoint between any pair of points in the set. In its contrapositive form, this standard characterization states that a closed set $S\subseteq \mathbb{R}^n$ is non-convex if and only if \[\exists\; \bx,\bx'\in S \quad\text{ such that }\quad \frac{\bx+\bx'}{2} \notin S.\]
Although this condition can be satisfied by binary MICP-R sets (consider the two-point set $\{0, 1\}$, which is clearly binary MICP-R), its pairwise infinite extension given by
\begin{equation}\label{infinite}\exists\; R\subseteq S,\; \abs{R}=\infty  \quad\text{ such that }\quad\frac{\bx+\bx'}{2} \notin S \quad\forall \bx,\bx'\in R\blue{,\; \bx\neq \bx'}\quad \end{equation}
clearly shows that $S$ is not binary MICP-R, because, based on the results described in \autoref{sec:BinaryMICPcontribution}, any binary MICP-R set is a finite union of convex sets and hence cannot satisfy \eqref{infinite}.

\blue{However, it is not immediately clear if it is possible for a set to both satisfy~\eqref{infinite} and be MICP-R, i.e., representable using unbounded integer variables. The primary contribution of this line of investigation is the result that the ``midpoint'' condition~\eqref{infinite} is in fact incompatible with MICP-R; i.e., it is not satisfied by any MICP-R set.}

  We use this fact to conclude \blue{that various non-convex sets are not MICP-R,} including \begin{itemize}
    \item[(1)] the set of rank-1 matrices, \item[(2)]the set of integer points on the parabola, and
    \item[(3)] the set of prime numbers.\end{itemize} These sets respectively have non-convex quadratic, (mixed) integer non-convex quadratic and (mixed) integer non-convex polynomial formulations\footnote{The formulation for the set of rank-1 $m\times n$ matrices is given by $\bm{X}=\bm{u}\bm{v}^T,\bm{u}\in \Real^m,\bm{v}\in \Real^n$, the formulation for the integer points in the parabola is given by $\bx=(z_1,z_1^2), z_1 \in \mathbb{Z}$, and a integer non-convex polynomial formulation for the prime numbers can be found in \cite{diophantine}.}.  Hence our negative result provides a separation between representability under these classes of formulations and \blue{MICP-R}.

\subsection{MICP under rationality assumptions}
\blue{We develop a slightly technical ``rationality'' condition that requires that rational affine transformations of a set $M$, of a particular type, produce bounded sets or sets with at least one rational recession direction. We refer to sets with MICP formulations induced by sets that satisfy this rationality condition as \textbf{rational MICP representable} (rational MICP-R) and show they have various interesting properties. These are the first such results for MICP-R without strong shape assumptions on the set $M$ such as the ellipsoidal restrictions of \cite{DelPia2016}.}

\textit{Our key structural result for rational MICP-R sets} is that a \blue{set $S\subseteq \Real^n$ that is closed and rational MICP-R} must be a finite union of sets that are \blue{either \{compact and convex\} or \{closed and periodic\}} if there exists an upper bound on the diameter\footnote{See \autoref{diameterdef}.} of any convex subset of $S$.
In particular, this result excludes the pathological example \eqref{badset2} and other similar examples from the class of rational MICP-R sets, without excluding interesting,  \blue{more well-behaved sets such as $S=\set{\bx\in \mathbb{Z}^3\,:\, \norm{\bra{x_1,x_2}}_2\leq x_3}$ (this set is periodic, but does not satisfy the ellipsoidal restrictions of \cite{DelPia2016}).}

\blue{This structural result allows us} to fully characterize rational MICP-R compact sets and rational MICP-R subsets of the natural numbers. Of potential interest to modelers, a corollary of our characterization of rational MICP-R compact sets is that for rational MICP-R \textit{closed} sets, binary variables are sufficient when the set itself is bounded (i.e., compact). In the case of rational MICP-R subsets of the natural numbers, our results allow us to clearly separate the modeling power of rational MICP formulations and rational MILP formulations.

Finally, we prove that the class of rational MICP-R sets is furthermore closed under unions, but, interestingly, not under intersections. In contrast, the class of rational MILP-R sets is closed under intersections, but not under unions. This disparity suggests that the class of rational MICP-R sets is structurally quite different from the class of rational MILP-R sets, even though the class of rational MILP-R sets is a sub-class of the class rational MICP-R sets and both classes are associated with periodic sets.

\section{Notation and Preliminaries}\label{sec:definitions}

Before formally stating our main results, we review some notation and results from convex analysis.

We let $\mathbb{Z}_+=\set{0,1,2,\ldots}$ be the set of non-negative integers, $\mathbb{N}=\set{1,2,\ldots}$ be the set of natural numbers, and $\sidx{k}=\set{1,2,\ldots,k}$. We use bold for vectors and matrices (e.g. $\bx \in \Real^n$ and $\bm{U}\in \Real^{n\times n}$) and non-bold for scalars (e.g. $x\in \Real$). In particular, we let $\bm{e}(d,i)\in \Real^d$ be the $i$-th unit vector (i.e. $\bm{e}(d,i)_j=0$ if $i\neq j$ and $\bm{e}(d,i)_i=1$ if $i=j$) and simplify this notation to  $\bm{e}(i)$ when the dimension $d$ is evident from the context. Otherwise, we follow the standard notation from \cite{Hiriart-Urruty1996}. We will often work with projections of a set $M\subseteq \mathbb{R}^{n+p+d}$ for some $n,p,d \in \mathbb{Z}_+$,  so we identify the variables in $\mathbb{R}^n$, $\mathbb{R}^p$ and $\mathbb{R}^d$ of this set as $\bm{x}$, $\bm{y}$ and $\bm{z}$ and we let
\[\proj_{\bx}\bra{M}=\set{\bm{x} \in \mathbb{R}^n\,:\,\exists \bra{\bm{y},\bm{z}}\in \mathbb{R}^{p+d} \text{ s.t. } \bra{\bm{x},\bm{y},\bm{z}}\in M}. \]
We similarly define $\proj_{\by}\bra{M}$ and $\proj_{\bz}\bra{M}$, dropping the bold for scalar variables (e.g. $\proj_z\bra{M}$ when $d=1$).

The following notation on set-valued maps (e.g \cite[Section 5.4]{borwein2010convex} and \cite[Chapter 5]{rockafellar2009variational}) will be useful for defining alternative characterizations of MICP-R sets.
\begin{definition}
   \blue{ Let $I\subseteq  \mathbb{R}^d$  and $F:I\rightrightarrows \mathbb{R}^{n+p}$ be a set-valued map, i.e., $F(\bz)\subseteq \mathbb{R}^{n+p}$ for all $\bz\in I$. The graph of $F$ is defined as $\operatorname{gr}(F)=\bigcup_{\bz\in I} F\bra{\bz}\times \set{\bz}$. We say $F$ is \textbf{convex} if $\operatorname{gr}(F)$ is convex and \textbf{closed} if $\operatorname{gr}(F)$ is closed. Note that $F$ is convex if and only if $I$ is convex and for all $\bm{z},\bm{z}'  \in I$ and $\lambda \in [0,1]$ we have $\lambda F\bra{\bm{z}}+ (1-\lambda) F\bra{\bm{z}'} \subseteq F\bra{\lambda \bm{z}+ (1-\lambda) \bm{z}'}$.    However,  $F$ may be closed even if $I$ is not closed\footnote{e.g. consider $I=(0,1]$ and $F\bra{\bz}=[1/z_1,\infty)$.}.}
\end{definition}

To formally define rational MICP representability we will need the following notion of recession directions for a convex set that is not necessarily closed  (cf. \cite[Section 8]{Rockafellar1997}).

\begin{definition}\label{rationaldefcone}
 Let $C \subseteq \mathbb{R}^d$ be a convex set, and for each $\bx\in C$ let $C_{\infty}\bra{\bx}=\set{\bm{r}\in \Real^d\,:\, \bm{x}+\lambda \bm{r} \in C \quad \forall \lambda\geq 0}$. We define the \textbf{recession cone} of $C$ as \[C_\infty=\bigcap_{\bx\in C} C_{\infty}\bra{\bx} = \set{\bm{r}\in \Real^d\,:\, \bm{x}+\lambda \bm{r} \in C \quad \forall \bm{x} \in C,\quad \lambda\geq 0}\] and the \textbf{closed recession cone} of $C$ as the recession cone of the topological closure of $C$ given by $C_{\overline{\infty}}=\bra{\cl\bra{C}}_\infty$. An element $\bm{r} \in C_\infty \setminus \set{\bm{0}}$ is a \textbf{recession direction}.
\end{definition}

\blue{As illustrated by  \autoref{recessionconeexex}, some common properties for recession cones of closed convex sets may fail to hold for  non-closed convex sets, even if they are projections of closed convex sets  (e.g. $\bra{\cdot}_\infty$ may not preserve set containment, and for $\bx\neq \bx'$ we may have $C_{\infty}\bra{\bx}\neq C_{\infty}\bra{\bx'}$). The following proposition, whose proof we include in Appendix~\ref{sec:recessionconeex} for completeness, provides a useful connection to the more familiar properties for closed convex sets.}
\begin{restatable}{proposition}{recessionconeprop}
\label{rationaldefcprop} Let $C \subseteq \mathbb{R}^d$ be a convex set. Then $C_\infty$ is a convex cone containing the origin and $C_{\overline{\infty}}=\bra{\cl\bra{C}}_\infty=\bra{\relint\bra{C}}_\infty$. If $C$ is additionally closed, then $C_\infty=C_{\overline{\infty}}$ is a non-empty closed convex cone, and $C_\infty=C_\infty\bra{\bx}$ for all $\bx \in C$. Furthermore, if  $C$ and $C'$ are closed convex sets such that $C\subseteq C'$, then $C_\infty\subseteq C'_\infty$.
\end{restatable}
\pointtoproof{recessionconeprop:proof}

Finally, we use the following definition for the diameter of a set.
\begin{definition}[Diameter]\label{diameterdef}
The \textbf{diameter of a set} $C\subseteq \Real^n$ is \blue{the maximum distance between any two points in the set; that is $D(C)=\sup\set{\norm{\bx - \bx'}_2 : \bx, \bx' \in C}.$}
\end{definition}

\section{MICP-R characterizations without rationality assumptions}\label{sec:main_results}

From \autoref{MICPFORMULATIONDEF} we see that if $M\subseteq \mathbb{R}^{n+p+d}$ induces an MICP formulation of $S\subseteq \mathbb{R}^{n}$, then
   \begin{equation}\label{infiniteunion}
   S=\bigcup_{\bm{z} \in I\cap \mathbb{Z}^d} \proj_{\bx}\bra{B_{\bm{z}}},\end{equation}
   where $I=\proj_{\bz}\bra{M}$ and  $B_{\bm{z}}=M \cap (\mathbb{R}^{n+p} \times \{\bm{z}\})$ for any $\bm{z} \in I$. Hence, MICP representable sets can be seen as a specially structured countable unions of projections of convex sets. We can precisely describe this special structure using properties of set-valued maps.

\begin{restatable}{theorem}{reptheotheo}\label{reptheo}
 A set $S \subseteq \mathbb{R}^n$ is MICP-R if and only if there exist $d,p\in \mathbb{Z}_+$, a  set $I\subseteq  \mathbb{R}^d$ and a closed convex set-valued map $F:I\rightrightarrows \mathbb{R}^{n+p}$ such that
\begin{equation}\label{infiniteunionfunction}S=\bigcup_{\bm{z} \in I\cap \mathbb{Z}^d} \proj_{\bx}\bra{F\bra{\bm{z}}}.\end{equation}
\end{restatable}
\proof{\textbf{Proof}}
\label{reptheotheo:proof}
\blue{
Let $S\subseteq \Real^n$ be an MICP-R set and  $M \subseteq \mathbb{R}^{n+p+d}$ be a closed convex set inducing a formulation for $S$ and $I=\proj_{\bz}\bra{M}$.
We define $F:I\rightrightarrows \mathbb{R}^{n+p}$  to be such that $F\bra{\bz}=\proj_{\bx,\by}\bra{M \cap (\mathbb{R}^{n+p} \times \{\bm{z}\})}$. Then,
\[\operatorname{gr}(F)=\bigcup_{\bz\in I}F\bra{\bz}\times \set{\bz}=\bigcup_{\bz\in I} \proj_{\bx,\by}\bra{M \cap (\mathbb{R}^{n+p} \times \{\bm{z}\})}\times \set{\bz}=M\]
and hence $F$ is a closed and convex set-valued map because $M$ is closed and convex. Finally, for any $\bz\in I\cap\mathbb{Z}^d $ we have $\proj_{\bx}\bra{F(\bz)}=  \proj_{\bx}\bra{\proj_{\bx,\by}\bra{M \cap (\mathbb{R}^{n+p} \times \{\bm{z}\})}}=  \proj_{\bx}\bra{B_{\bm{z}}}$,
where $B_{\bm{z}}=M \cap (\mathbb{R}^{n+p} \times \{\bm{z}\})$ so \eqref{infiniteunionfunction} follows from \eqref{infiniteunion}.

For the converse, let
$M=\operatorname{gr}(F)$. Then $M$ is closed and convex because $F$ is closed and convex. In addition, because $M=\bigcup_{\bz\in I}F\bra{\bz}\times \set{\bz}$, for all $\bz\in I$ we have $M \cap (\mathbb{R}^{n+p} \times \{\bm{z}\})=F\bra{\bz}\times \set{\bz}$ and, in particular,  $\proj_{\bx}\bra{M \cap (\mathbb{R}^{n+p} \times \{\bm{z}\})}=\proj_{\bx}\bra{F\bra{\bz}\times \set{\bz}}=\proj_{\bx}\bra{F\bra{\bz}}$. Then,    \eqref{infiniteunionfunction} implies \eqref{infiniteunion} with $B_{\bm{z}}=M \cap (\mathbb{R}^{n+p} \times \{\bm{z}\})$ and $M$ induces an MICP formulation of $S$. }
 \Halmos\endproof

Given \autoref{reptheo}, the main goal of this paper is to understand properties of this specially structured union. The sets $I$ and $\proj_{\bx}\bra{B_{\bz}}=\proj_{\bx}\bra{F\bra{\bz}}$ from  \eqref{infiniteunion}--\eqref{infiniteunionfunction} will play a central role in the analysis.

  \begin{definition}\label{def:index}
      Let $M\subseteq \mathbb{R}^{n+p+d}$ be a closed, convex set that induces an MICP formulation of $S \subseteq \mathbb{R}^n$. We refer to $I=\proj_{\bz}\bra{M}$ as the \textbf{index set} of the MICP formulation and to \blue{the collection of sets $\set{A_{\bm{z}}}_{\bz\in I}$ with $A_{\bm{z}}=\proj_{\bx}\bra{M \cap (\mathbb{R}^{n+p} \times \{\bm{z}\})}$ for each $\bz \in I$, as its \textbf{$\bm{z}$-projected sets}. }
  \end{definition}

\subsection{Binary MICP representability}\label{sec:bounded_intro}

Jeroslow's characterization of pure binary MICP-R from \cite{jeroslow1987representability} can be shown through the following result that highlights the non-necessity of introducing any continuous auxiliary variables in the formulations of such sets\footnote{The characterization from \cite{jeroslow1987representability} focuses on a function-based restriction of binary MICP-R sets that guarantees the represented sets are closed. \autoref{prop:jeroslowbinaryequiv} gives a formal proof of the equivalence between this function-based notation and our definition of pure binary MICP-R.}.

\begin{restatable}{proposition}{balasformlem}\label{lem:balasform_jeroslow}
 $S\subseteq  \mathbb{R}^{n}$ is pure binary MICP-R if and only if there exist nonempty, closed, convex sets $\set{S_i}_{i=1}^d$ such that $S_i\subseteq \mathbb{R}^{n}$ for each $i\in\sidx{d}$, $\bra{S_i}_\infty = \bra{S_j}_\infty$ for all $i,j\in \sidx{d}$, and $S=\bigcup_{i =1}^d S_i$. For such an $S$ we have that $\bx\in S$ if and only if
 \begin{equation}\label{form:extendedbounded_jeroslow}
    \exists \bz\in\mathbb{Z}^d \text{ s.t. } \bra{\bx,\bz}\in M= \conv\bra{\bigcup_{i=1}^d \bra{S_i \times \set{\bm{e}(d,i)}}}\subseteq \Real^{n+d}.
 \end{equation}
  Furthermore, $M$ defined in \eqref{form:extendedbounded_jeroslow} is a closed convex set such that $\proj_{\bz}\bra{M\cap \bra{\mathbb{R}^{n} \times \mathbb{Z}^d}}\subseteq\set{0,1}^d$ and hence induces a pure binary MICP formulation of $S$.
\end{restatable}
\proof{\textbf{Proof}}
Let $\set{S_i}_{i=1}^d$ be a family of non-empty, closed convex sets such that $S_i\subseteq \mathbb{R}^{n}$ for each $i\in\sidx{d}$ and $\bra{S_i}_\infty = \bra{S_j}_\infty$ for all $i,j\in \sidx{d}$. Then Corollary 9.8.1 \cite{Rockafellar1997} shows that $M$ is closed. Furthermore, $M\subseteq \Real^n\times [0,1]^d$ so $\proj_{\bz}\bra{M\cap \bra{\mathbb{R}^{n} \times \mathbb{Z}^d}}\subseteq\set{0,1}^d$. In addition, if
 $\bra{\bx,\bz}\in M\cap\bra{\Real^{n}\times \mathbb{Z}^d}$, then there exist $i\in\sidx{d}$ such that $\bz=\bm{e}(d,i)$ and $\bx\in S_i$. Hence $M$ induces a pure binary MICP formulation of $S$.

For the converse let  $M\subseteq\Real^{n+d}$ be a closed convex set such that $S=\proj_{\bx}
    \bra{M\cap \bra{\mathbb{R}^{n} \times \mathbb{Z}^d}}$ and $\proj_{\bz}\bra{M\cap \bra{\mathbb{R}^{n} \times \mathbb{Z}^d}}\subseteq\set{0,1}^d$. For each $\bz\in \set{0,1}^d$ let $M_{\bz}=M \cap (\mathbb{R}^{n} \times \{\bm{z}\})$ and $A_{\bz}=\proj_{\bx}\bra{M_{\bz}}$ so that  $M_{\bz}=A_{\bz}\times \{\bm{z}\}$. Then, $M_{\bz}$ is closed, and if $M_{\bz}\neq\emptyset$, then by Corollary 8.3.3 in \cite{Rockafellar1997} we have $\bra{M_{\bz}}_{\infty}=M_{\infty} \cap (\mathbb{R}^{n} \times \{\bm{0}\})$. In addition,  if $A_{\bz}\neq \emptyset$ by Theorem 9.1 in \cite{Rockafellar1997} we have that $A_{\bz}$ is a nonempty closed convex set and $\bra{A_{\bz}}_{\infty}=\proj_{\bx}\bra{M_{\infty} \cap (\mathbb{R}^{n} \times \{\bm{0}\})}$. Hence,
    for any  $\bz,\bz'\in \set{0,1}^d$ such that  $A_{\bz}$ and  $A_{\bz'}$ are non-empty we have that $\bra{A_{\bz}}_{\infty}=\bra{A_{\bz'}}_{\infty}$. Finally, $S=\bigcup_{\bz\in \set{0,1}^d\,:\, A_{\bz}\neq \emptyset}A_{\bz}$.
\Halmos\endproof

\begin{figure}[htpb]\centering
  \includegraphics{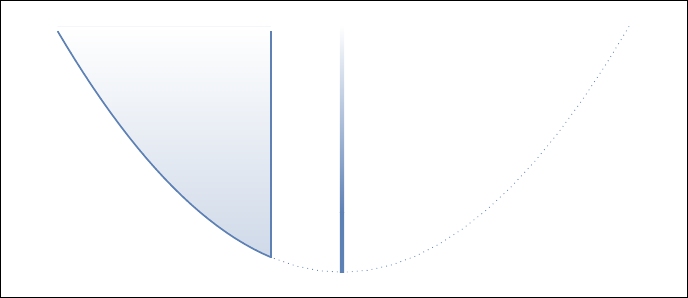}
%
%
\caption{The union of two closed convex sets with equal recession cones from \autoref{Jeroslowexample}.}\label{fig:jeroslowrecession}
\end{figure}

The following example illustrates how the class of pure binary MICP-R sets is strictly contained in the class of binary MICP-R sets.
\begin{example}\label{Jeroslowexample}
Consider the closed convex set
\[M= \set{\bra{x,y,z}\in \Real \times \Real^2_+\,:\, x^2 \le yz,\quad x \leq -z,\quad 0 \leq z \leq 1}.\]
Closure and convexity of $M$ follows because the set $\{ \bra{x, y, z} \in \Real \times \Real^2_+  \,:\, x^2 \le yz \}$, called the rotated second order cone, is known to be an invertible linear transformation of the second order cone or quadratic cone, which itself is closed and convex~\cite{alizadeh2003second,MOSEKModeling}. Then the following observations are in order.
\begin{enumerate}
    \item We have that for $z=0,$ $S_1=\proj_{x}\bra{M\cap\bra{\Real^2\times \set{0}}}=\set{0}$ and for $z=1,$ $S_2=\proj_{x}\bra{M\cap\bra{\Real^2\times \set{1}}}=(-\infty,-1]$.
    \item The closed set $S=S_1\cup S_2$ is binary MICP-R but not pure binary MICP-R. Indeed, based on the previous observation, $M$ induces a binary MICP formulation of exactly $S_1\cup S_2=S$. However, $S_1$ and $S_2$ have different recession cones and have empty intersection, which based on Proposition \ref{lem:balasform_jeroslow} allows us to conclude that the set $S$ is not pure binary MICP-R. In particular,  $\proj_{x,z}\bra{M\cap\bra{\Real^2\times [0,1]}}=\set{\bra{0,0}}\cup \bra{\Real_+\times (0,1]}$, which is not closed.
    \item \label{Jeroslowexamplep3} The set $T=\proj_{x,y}\bra{M\cap\bra{\Real^2\times \set{0,1}}}$ (depicted in \autoref{fig:jeroslowrecession}) is the union of $T_1=\set{\bra{x,y}\in \Real^2\,:\, x=0,\quad y\geq 0}$ and $T_2=\set{\bra{x,y}\in \Real^2\,:\, x^2\leq y,\quad x\leq -1}$,  which are closed convex sets with the same recession cone. Hence $T$ is pure binary MICP-R by \autoref{lem:balasform_jeroslow}. Finally, $S_1=\proj_{x}\bra{T_1}$, $S_2=\proj_{x}\bra{T_2}$, and $S=\proj_{x}\bra{T}$.
\end{enumerate}
\Halmos
\end{example}

\autoref{Jeroslowexample} shows that the class of pure binary MICP-R sets is a strict subset of the class of binary MICP-R sets, even if we restrict to MICP-R sets that are also closed. Furthermore,  \autoref{Jeroslowexamplep3} in \autoref{Jeroslowexample} shows that the class of pure binary MICP-R sets is not closed under orthogonal projections, as the projection of finite family of closed convex sets with the same recession cones (e.g. $T_1$ and $T_2$) may lose this latter property even if they remain closed (e.g.  $S_1$ and $S_2$). Finally, while \autoref{Jeroslowexample} shows that binary MICP-R sets, such as $S$, may be unions of finitely many closed convex sets with different recession cones,  \autoref{Jeroslowexamplep3} suggests that a formulation for such sets may be obtained by lifting the closed convex sets into a higher dimensional space (i.e. by adding continuous auxiliary variables) where they do share the same recession cone. After the lifting step, a formulation could be used based on Proposition \ref{lem:balasform_jeroslow}. The following proposition shows that this lifting can indeed always be achieved and provides a complete characterization of binary MICP-R sets.

\begin{restatable}{proposition}{balasformlem_new}\label{lem:balasform}
 $S\subseteq  \mathbb{R}^{n}$ is binary MICP-R if and only if there exists a $q \in \mathbb{Z}_+$ and nonempty, closed, convex sets $\set{\tilde{S}_i}_{i=1}^d$ (without any restriction on their recession cones) such that $\tilde{S}_i\subseteq \mathbb{R}^{n+q}$ for each $i\in\sidx{d}$ and $S=\bigcup_{i =1}^d S_i$, where $S_i= \operatorname{proj}_{\bx} \bra{ \tilde{S}_i}$ for all $i\in \sidx{d}$. For such an $S$ we have that $\bx\in S$ if and only if
 \begin{equation}\label{form:extendedbounded_new}
  \exists \by\in \Real^p,\,\bz\in\mathbb{Z}^d \text{ s.t. } \bra{\bx,\by,\bz}\in   M= \conv\bra{\bigcup_{i=1}^d {T}_i \times \set{\bm{e}(d,i)}}\subseteq \Real^{n+p+d},
 \end{equation}
 where $p=q+1$, and for each $i\in \sidx{d}$
 \[{T}_i=\set{\bra{\bm{w},w_0}\in \Real^{n+q+1}\,:\,\bm{w}\in \tilde{S}_i,\quad \norm{\bm{w}}_2^2\leq w_0}. \]
 Furthermore, $M$ defined in \eqref{form:extendedbounded_new} is a closed convex set such that $\proj_{\bz}\bra{M\cap \bra{\mathbb{R}^{n+p} \times \mathbb{Z}^d}}\subseteq\set{0,1}^d$ and hence induces a binary MICP formulation of $S$.
\end{restatable}
\proof{\textbf{Proof}}
Closure of $M$ follows from Corollary 9.8.1 \cite{Rockafellar1997} because for each $i\in \sidx{d}$,  we have that ${T}_i \times \set{\bm{e}(d,i)}$ is a nonempty, closed convex set and  \[\bra{{T}_i \times \set{\bm{e}(d,i)}}_\infty=\set{\bra{\bm{w},w_0,\bz}\in \Real^{n+q+1+d}\,:\, \bm{w}=\bm{0},\quad \bz=\bm{0},\quad w_0\geq 0 }
.\]
Furthermore, $M\subseteq \Real^{n+p}\times [0,1]^d$ so $\proj_{\bz}\bra{M\cap \bra{\mathbb{R}^{n+p} \times \mathbb{Z}^d}}\subseteq\set{0,1}^d$. In addition, through the identification $\bra{\bx,\by}=\bra{\bm{w},w_0}\in \Real^{n+p}=\Real^{n+q+1}$
we have that if $\bra{\bx,\by,\bz}\in M\cap\bra{\Real^{n+p}\times \mathbb{Z}^d}$, then there exists $i\in\sidx{d}$ such that $\bz=\bm{e}(d,i)$, $\bra{\bx,\by}\in T_i$ and hence  $\bx\in \operatorname{proj}_{\bx} \bra{ \tilde{S}_i}$. Then $\proj_{\bx}
    \bra{M\cap \bra{\mathbb{R}^{n+p} \times \mathbb{Z}^d}}\subseteq S$. For the reverse containment, it suffices to note that $(\bm{w}, ||\bm{w}||_2^2) \in T_i$ for any $i$ and $\bm{w} \in \tilde{S}_i$. Hence,  \[S_i=\operatorname{proj}_{\bx} \bra{ \tilde{S}_i}\subseteq \operatorname{proj}_{\bx} \bra{ T_i}=\proj_{\bx}
    \bra{M\cap \bra{\mathbb{R}^{n+p} \times \set{\bm{e}(d,i)}}} \subseteq \proj_{\bx}
    \bra{M\cap \bra{\mathbb{R}^{n+p} \times \mathbb{Z}^d}}.\]
    Then, $M$ induces a binary MICP formulation of $S$.

For the converse let $M\subseteq\Real^{n+p+d}$ be a nonempty, closed convex set such that $S=\proj_{\bx}
    \bra{M\cap \bra{\mathbb{R}^{n+p} \times \mathbb{Z}^d}}$ and $\proj_{\bz}\bra{M\cap \bra{\mathbb{R}^{n+p} \times \mathbb{Z}^d}}\subseteq\set{0,1}^d$. Then for each $\bz\in \set{0,1}^d$ we have that $B_{\bz}=\proj_{\bx,\by}\bra{M \cap (\mathbb{R}^{n+p} \times \{\bm{z}\})}$ is either empty, or a nonempty closed convex set.  Finally, $S=\bigcup_{\bz\in \set{0,1}^d\,:\, B_{\bz}\neq \emptyset}\proj_{\bx}\bra{B_{\bz}}$.
\Halmos\endproof

Note that for a finite union of closed convex sets, \autoref{lem:balasform} applies with $q=0$. In such case, $p=1$, and hence to remove the equal-recession-cone restriction from \autoref{lem:balasform_jeroslow} it is sufficient to add a single continuous auxiliary variable.

\subsection{\blue{Obstructions} to MICP representability}

\autoref{reptheo} states that MICP-R sets are countable unions of projections of closed convex sets. However, as illustrated by the following proposition, some highly non-convex sets can be written in this way.

\begin{restatable}{proposition}{dcrep}\label{dcrepref}
The complement of any convex body (i.e. a full-dimensional compact convex set) is a countable union of projections of closed convex sets.
\end{restatable}
\proof{\textbf{Proof}}
Let $K\subseteq \Real^n$ be a convex body and $\sigma_K:\Real^n\to\Real\cup \set{+\infty}$ be its support function defined by $\sigma_K\bra{\bm{c}}=\operatorname{sup}\{ {\bm{c}}^T\bm{x} : \bm{x} \in K \}$.
 By Proposition 2.1 in \cite{The-Chvatal-Gomory-Closure-of-a-Strictly} we have $K=\set{\bx\in \Real^n\,:\ \bm{a}^T\bx \leq \sigma_K\bra{\bm{a}}, \quad \forall \bm{a}\in \mathbb{Z}^n}$.
Then \[\Real^n\setminus K=\bigcup_{a\in \mathbb{Z}^n} \set{\bx\in \Real^n\,:\ \bm{a}^T\bx > \sigma_K\bra{\bm{a}}}=\bigcup_{a\in \mathbb{Z}^n} \proj_{\bx}\bra{\set{\bra{\bx,y}\in \Real^{n}\times \Real_+\,:\ \bm{a}^T\bx \geq 1/y+ \sigma_K\bra{\bm{a}}}}.\]
\Halmos\endproof

Not surprisingly, many such countable unions of closed convex sets are not MICP-R. We develop a geometric obstruction to MICP representability to prove this. The result is based on the following notion of nonconvexity of a set and a definition of the MICP rank.
\begin{definition}\label{def:wstrongnonconvex}
    Let $w \in \bra{[2,\infty)\cap\mathbb{Z}} \cup \{ + \infty \}$. We say that a set $S \subseteq \mathbb{R}^n$ is \textbf{$w$-strongly nonconvex}, if there exists a subset $R \subseteq S$ with $|R| = w$ such that for all pairs $\bm{x},\bm{y} \in R$,
    $\bm{x}\neq \bm{y}$,
\begin{equation}
\frac{\bm{x}+\bm{y}}{2} \not \in S,
\end{equation}
that is, a subset of points in $S$ of cardinality $w$ such that the midpoint between any pair is not in $S$.
\end{definition}
\begin{definition}
   \blue{The \textbf{MICP rank} of an MICP-R set $S$ is the smallest integer $d$ such that there exists a closed convex set $M$ such that the pair $(M,d)$ induces an MICP formulation of $S$. If $S$ is not MICP-R, its rank is defined to be $\infty$.}
\end{definition}

The following lemma gives necessary conditions for a set to have a given finite MICP rank, and hence provides necessary conditions for a set to be MICP-R. We call this \textit{the Midpoint Lemma}.
\begin{restatable}[The Midpoint Lemma]{lemma}{midpointlemma}\label{midpointlemmaref}
\label{midpoint}
Let $S \subseteq \mathbb{R}^n$. If $S$ is $w$-strongly nonconvex, then the MICP \blue{rank} of $S$ must be at least $\lceil \log_2(w) \rceil$. In particular, $S$ is not MICP-R if either
\begin{itemize}
    \item $S$ is $\infty$-strongly nonconvex, or
    \item $S$ is $w$-strongly nonconvex for all $w \in [2,\infty)\cap\mathbb{Z}$.
\end{itemize}
\end{restatable}
\proof{\textbf{Proof}}
Suppose we have $R$ as in \autoref{def:wstrongnonconvex}.

First, note that if the set $S$ is not MICP-R, the statement of the result follows immediately. Otherwise, we proceed with the proof by contradiction. Suppose the existence of \blue{a pair $(M,d)$ that induces} an MICP formulation of $S$ with \blue{$d < \lceil \log_2(w) \rceil$}. \blue{By definition, $M \subseteq \mathbb{R}^{n+p+d}$ is a closed convex set} such that $\bm{x} \in S$ iff $\exists \bm{z} \in \mathbb{Z}^d,\by \in \mathbb{R}^p$ such that $(\bm{x},\bm{y},\bm{z}) \in M$.
Then for each point $\bm{x} \in R$ we associate at least one integer point $\bm{z}^{\bm{x}} \in \mathbb{Z}^d$ and a $\bm{y}^{\bm{x}} \in \mathbb{R}^p$ such that $(\bm{x},\bm{y}^{\bm{x}},\bm{z}^{\bm{x}}) \in M$. If there are multiple such pairs of points $\bm{z}^{\bm{x}},\bm{y}^{\bm{x}}$ then for the purposes of the argument we may choose one arbitrarily.

Now recall we assume $d < \lceil \log_2(w) \rceil$. We will derive a contradiction by proving that there exist two points $\bm{x},\bm{x}' \in R$ such that the associated integer points $\bm{z}^{\bm{x}},\bm{z}^{\bm{x}'}$ satisfy
\begin{equation}\label{eq:vecparity}
\frac{\bm{z}^{\bm{x}}+\bm{z}^{\bm{x}'}}{2} \in \mathbb{Z}^d.
\end{equation}
Indeed, this property combined with convexity of $M$, i.e., $\left(\frac{\bm{x}+\bm{x}'}{2},\frac{\bm{y}^{\bx}+\bm{y}^{\bm{x}'}}{2},\frac{\bm{z}^{\bm{x}}+\bm{z}^{\bm{x}'}}{2}\right) \in M$
would imply that $\frac{\bm{x}+\bm{x}'}{2} \in S$, which contradicts the definition of $R$.

Recall a basic property of integers that if $i,j \in \mathbb{Z}$ and $i \equiv j \text{ (mod 2)}$, i.e., $i$ and $j$ are both even or odd, then $\frac{i+j}{2} \in \mathbb{Z}$. We say that two integer vectors $\bm{\alpha},\bm{\beta} \in \mathbb{Z}^d$ have the same \textit{parity} if $\alpha_i$ and $\beta_i$ are both even or odd for each component $i \in \sidx{d}$. Trivially, if $\bm{\alpha}$ and $\bm{\beta}$ have the same parity, then $\frac{\bm{\alpha}+\bm{\beta}}{2} \in \mathbb{Z}^d$. Given that we can categorize any integer vector according to the $2^d$ possible choices for whether its components are even or odd, and we notice that from any collection of integer vectors of size greater than $2^d+1$ we must have at least one pair that has the same parity. Therefore, using our assumption $d < \lceil \log_2(w) \rceil$, we have $|R| = w \ge 2^d+1$, which implies the existence of a pair $\bm{x},\bm{x}' \in R$ such that their associated integer points $\bm{z}^{\bm{x}},\bm{z}^{\bm{x}'}$ have the same parity and thus satisfy~\eqref{eq:vecparity}, leading to the desired contradiction. It follows that the MICP \blue{rank} of $S$ must be at least $\lceil \log_2(w) \rceil$.

The remaining statements in the lemma follow directly from the first result as the MICP \blue{rank} of any such set $S$ cannot be finite.
\Halmos\endproof

We can use the last statement of \autoref{midpoint} to show that many classes of sets fail to be MICP-R. The following corollary is a direct consequence of \autoref{midpointlemmaref}.

\begin{restatable}{corollary}{negativeresult}\label{negativeresultref}
 The following sets are not MICP-R.
 \begin{itemize}
  \item The set $\{ X \in \mathbb{R}^{m\times n} : \operatorname{rank}(X) \le 1 \}$ of $m$ by $n$ matrices with rank at most $1$ for $m,n \geq 2$.
\item The complement of a strictly convex body.
\item The spherical shell $\set{\bx\in \Real^n\,:\, 1\leq \norm{\bx}_2\leq 2}$.
\item The set $\bigcup_{i\in \mathbb{Z}}\set{\bra{i,i^2},\bra{i+1,(i+1)^2}}$ of integer points in the parabola and its piecewise linear interpolation given by  $\bigcup_{i\in \mathbb{Z}}\conv\bra{\set{\bra{i,i^2},\bra{i+1,(i+1)^2}}}$.

 \end{itemize}
\end{restatable}
\proof{\textbf{Proof}}
In what follows, for any  $n\in \mathbb{N}$ and $R\subseteq \Real^n$, let $\mathcal{A}\bra{R}=\set{\frac{\bx+\by}{2}\,:\, \bx,\by\in R,\quad \bx\neq \by}$ be the set of all midpoints between elements of $\blue{R}$.

Let $C_1=\{ X \in \mathbb{R}^{m\times n} : \operatorname{rank}(X) \le 1 \}$. \blue{It suffices to prove the statement for} $m=2$. We set for all $k \in \mathbb{N}$ the matrix
$A_k = \begin{bmatrix}
    1      & k   & O_{1 \times n-2}\\
    k       & k^2 & O_{1 \times n-2}
\end{bmatrix} \in C_1$. We then set  $R=\{A_k | k \in \mathbb{N}  \}$. Clearly $|R|=\infty$. It is easy to verify that $\operatorname{rank}(\frac{1}{2} (A_k+A_{k'})) = 2$ for $k \ne k'$. Therefore $\mathcal{A}\bra{R}\subseteq \Real^{ m \times n}\setminus C_1$ and hence $C_1$ is $\infty$-strongly nonconvex and in particular not MICP-R.

Let $C_2=\mathbb{R}^n\setminus K$ where $K\subseteq \Real^n$ is a strictly convex body, $w \in \bra{[2,\infty)\cap\mathbb{Z}}$ and $R'\subseteq K\setminus \Int\bra{K}$ with $\abs{R'}=w$ be a set of $w$ distinct points in the boundary of $K$. Because $K$ is strictly convex, then $\mathcal{A}\bra{R'}\subseteq\Int\bra{K}$. Without loss of generality (by possibly translating $K$) we may assume that $0\in \Int\bra{K}\setminus \mathcal{A}\bra{R'}$. Then, because $\mathcal{A}\bra{R'}\subseteq\Int\bra{K}$, $\abs{\mathcal{A}\bra{R'}}<\infty$ and $0\notin \mathcal{A}\bra{R'}$, we have that  $\varepsilon=\sup\set{\lambda\,:\, \mathcal{A}\bra{R'}\subseteq (1-\lambda) K}\in (0,1)$. This implies that $\frac{1}{1-\lambda/2}\mathcal{A}\bra{R'}\subseteq K$. Furthermore, because $\frac{1}{1-\lambda/2}>1$ we have that   $R=\frac{1}{1-\lambda/2}R'\subseteq C_2=\mathbb{R}^n\setminus K$. Finally, because $\mathcal{A}\bra{R}=\frac{1}{1-\lambda/2}\mathcal{A}\bra{R'}$ we have that  $\mathcal{A}\bra{R}\subseteq \mathbb{R}^n \setminus C_2$. Therefore $C_2$ is $w$-strongly nonconvex for all $w \in \bra{[2,\infty)\cap\mathbb{Z}}$ and hence is not MICP-R.

Let $C_3=\set{\bx\in \Real^n\,:\, 1\leq \norm{\bx}_2\leq 2}$ and $R=\set{\bx\in \Real^n\,:\, \norm{\bx}_2=1}$. Then $R\subseteq C_3$, $\abs{R}=\infty$ and $\mathcal{A}\bra{R}\subseteq \mathbb{R}^n \setminus C_3$. Therefore $C_3$ is $\infty$-strongly nonconvex and in particular not MICP-R.

Let $C_4=\bigcup_{i\in \mathbb{Z}}\set{\bra{i,i^2},\bra{i+1,(i+1)^2}}$ or $C_4=\bigcup_{i\in \mathbb{Z}}\conv\bra{\set{\bra{i,i^2},\bra{i+1,(i+1)^2}}}$. In both cases,  $R=\bigcup_{i\in \mathbb{Z}}\set{\bra{2i,(2i)^2}}\subseteq C_4$, $\abs{R}=\infty$ and $\mathcal{A}\bra{R}\subseteq \mathbb{R}^n \setminus C_4$. Therefore both versions of $C_4$ are $\infty$-strongly nonconvex and in particular not MICP-R.
\Halmos\endproof

Notice that according to \blue{the second bullet of \autoref{negativeresultref}}, we can exclude the sub-class of the sets from Proposition~\ref{dcrepref} that is obtained by requiring that the convex body is additionally \emph{strictly convex} (i.e. any strict convex combination of two points in the set lies in the interior of the set). This strict convexity condition cannot be significantly relaxed as \autoref{lem:balasform} shows that the complement of a full-dimensional polytope is MICP-R.

One can also ask which subsets of the natural numbers are MICP-R. A subset of distinct interest is the set of prime numbers, whose regularity properties have been an object of interest at least since 300BC with the work of Euclid. The following result is a nontrivial application of \autoref{midpointlemmaref} and shows that the set of prime numbers lacks enough regularity to be  MICP-R.

\begin{restatable}{proposition}{primes}\label{primesref}
The set of prime numbers is not MICP-R.
\end{restatable}
\proof{\textbf{Proof}}
To construct $R$ as in \autoref{def:wstrongnonconvex}, we will inductively construct a subset of primes such that no midpoint of any two elements in the set is prime.

Let $\set{p_i}_{i=1}^n$ be a set of such primes. We will find a prime $p$ such that $\set{p_i}_{i=1}^n\cup\set{p}$ has no prime midpoints. We may start the induction with $p_1 = 3$, $p_2=5$.

Set $M = \prod_{i=1}^n p_i$. Choose any prime $p$ (not already in our set and not equal to 2) such that $p \equiv 1 \pmod{M!}$. By Dirichlet's theorem on arithmetic progressions \cite[Theorem 15]{hardy1979introduction}, there exists an infinite number of primes of the form $1 + kM!$ because 1 and $M!$ are coprime, so we can always find such $p$.

Suppose for some $i$ we have $q = \frac{p+p_i}{2}$ is prime. By construction, we have $p + p_i \equiv 1 + p_i\pmod{M!}$, so $\exists\, k$ such that $p + p_i = k\cdot M! + 1 + p_i$. Note that $M$ is larger than $p_i$, so $M!$ will contain $(1+p_i)$ as a factor; in other words, $(1+p_i)$ divides $M!$, so it divides also $k \cdot M!+1+p_i=p+p_i$. In fact, we can write $p+p_i = k'(1+p_i)$ for some $k' \in \mathbb{Z}_{+}$. We claim that $k'=1$. Indeed $q =\frac{p+p_i}{2}= k'\frac{1+p_i}{2}$. Note $1+p_i$ is even, so $\frac{1+p_i}{2}$ is an integer bigger than 1 as $p_i>1$. But $q$ is prime, and therefore, since it is written as the product of $k'$ and $\frac{1+p_i}{2}>1$, it must be the case that $k'=1$ as claimed. But $k'=1$ implies that $p+p_i =1+p_i$, i.e., $p = 1$ which is a contradiction.
\Halmos\endproof

Finally, with finite $w$ we  obtain the following interesting result  on modeling subsets $S$  of the binary hypercube $\{0,1\}^n$. It is clear that a formulation with binary integer variables requires at least $\lceil\log_2 \abs{S}\rceil$ binary integer variables (e.g. \cite[Proposition 1]{huchette2016small}). In addition, the following simple corollary of \autoref{midpoint} shows  that the same lower bound holds if we use unbounded integer variables. That is, using general integer variables instead of binary variables \textit{does not provide any advantage} in this modeling task.
\begin{corollary} Let $S\subseteq \{0,1\}^n$. Then any MICP formulation of $S$ requires at least $\lceil\log_2 \abs{S}\rceil$ integer variables.
\end{corollary}
\proof{\textbf{Proof}}
By picking $R=S$ in \autoref{def:wstrongnonconvex} we have that $S$ is $\lceil\log_2 \abs{S}\rceil$-strongly non-convex. The result then follows from \autoref{midpoint}.
\Halmos\endproof

\section{Rationality and periodicity}\label{MICPrat}

\autoref{JLTheo} and its  rational ellipsoidal extension in \cite{DelPia2016} can be further generalized as follows.
\begin{restatable}{proposition}{boundedpropres}\label{boundedprop}
   If $M$ induces an MICP formulation of $S\subseteq \Real^n$ and $M=B+K$ where $B$ is a compact convex set and $K$ is a rational polyhedral cone, then  there exist a finite set  $U\subseteq \mathbb{Z}^n$ and compact convex sets $\set{S_i}_{i=1}^k$ such that $S_i\subseteq \Real^n$ for all $i\in \sidx{k}$ and
   \begin{equation}\label{eq:boundedrep}
   S= \bigcup_{i=1}^k S_i + \operatorname{intcone}\bra{U}.
   \end{equation}
   \end{restatable}
\proof{\textbf{Proof}}
\blue{The result follows from a straightforward extension of Theorem 11.6 of~\cite{Vielma2015}}.\Halmos\endproof

We recover the rational MILP-R result from \cite{Jeroslow1984} when both $B$ and $\set{S_i}_{i=1}^k$ from \autoref{boundedprop} are rational polytopes (i.e. rational bounded polyhedra). In such case, we have that structure \eqref{eq:boundedrep} is both necessary and sufficient for rational MILP-R (e.g. see \autoref{JLTheo}). However, it is not hard to find even periodic MICP-R sets that do not satisfy the characterization from \autoref{boundedprop}\footnote{e.g. Corollary 1.4 in \cite{hemmecke2007representation} shows that $S=\set{\bx\in \mathbb{Z}^3\,:\, \norm{\bra{x_1,x_2}}_2\leq x_3}$ does not satisfy \eqref{eq:boundedrep}, but $S$ is periodic according to \autoref{periodicdef} with $\bm{u}=\bra{0,0,1}$.}. For this reason we focus on conditions on MICP-R sets that guarantee some level of periodicity without necessarily having the structure from \eqref{eq:boundedrep}.
\blue{ Specifically, we aim for a condition similar to the following consequence of \eqref{eq:boundedrep}.

\begin{corollary}\label{boundedorperiodic:coro}
If $S$ satisfies condition \eqref{eq:jlregular} from \autoref{JLTheo} or condition \eqref{eq:boundedrep} from  \autoref{boundedprop}, then either $S$ is a finite union of compact convex sets or $S$ is a closed periodic set. \end{corollary}
\proof{\textbf{Proof}}
First note that both \eqref{eq:jlregular} and  \eqref{eq:boundedrep} imply that $S$ is closed. Furthermore, for both \eqref{eq:jlregular} and  \eqref{eq:boundedrep}, $S$ is a finite union of compact convex sets if $U\setminus \set{\bm{0}} =\emptyset$. Otherwise, for any $\bm{u} \in U\setminus \set{\bm{0}}$ we have  $\bm{x} + \lambda \bm{u}\in S$ for all $\bx \in S$ and $\lambda \in \mathbb{N}$, and hence $S$ is periodic.
\Halmos\endproof}

   \subsection{Rational MICP}\label{ratsection}

   For $x\in \mathbb{R}$ let $f(x)=x-\lfloor x \rfloor$, and consider the  non-periodic set from~\eqref{badset2} given by $S=\set{x\in \mathbb{Z}: f\bra{x\sqrt{2} }\notin \bra{0.4, 1-0.4} }$, which, as noted in \eqref{irrationalex2}, has an MICP formulation  induced by  $M=\set{\bra{x,\bz}\in \mathbb{R}\times \Real^2\,:\, z_1=x,\; -0.4\leq z_2 - \sqrt{2} z_1 \leq     0.4}$. A notable property of the formulation induced by $M$ is that its index $I=\set{\bz\in \Real^2\,:\,   -0.4\leq z_2 - \sqrt{2} z_1 \leq     0.4}$ has a recession cone $I_\infty=\set{\bz\in \Real^3\,:\, z_2=\sqrt{2} z_1}$ that is a non-rational subspace such that $\mathbb{Z}^3 \cap {I}_\infty =\set{\bm{0}}$. Avoiding this property for $I$ and any rational affine transformation of $I$ is precisely the technical rationality condition that yields a generalization of \autoref{boundedorperiodic:coro}.

\begin{definition}\label{rationaldefc}
    We say a set $I \subseteq \mathbb{R}^d$ is \textbf{rationally unbounded} if for any rational affine image  $I'\subseteq \Real^{d'}$ of $I$, either $I'$ is a bounded set or it holds that $ \mathbb{Z}^{d'} \cap (I'_\infty \setminus \set{\bm{0}}) \not = \emptyset$.
    \end{definition}
\begin{definition}\label{rationaldefcformulation}
  We say that a set $S$ is \textbf{rational MICP representable} (rational MICP-R) if it has an MICP representation induced by the set $M$ whose index set  $I=\proj_{\bz}(M)$ is rationally unbounded.
\end{definition}

As illustrated by the following proposition and example, \autoref{rationaldefc} restricts irrational recession directions in the recession cone of $I$. However, it additionally restricts irrational directions that appear after certain projections of $I$ (a redundant requirement if $I$ is a polyhedron, such as when $I$ is the index set of a rational MILP formulation; cf. \autoref{rationalmicpmilp}).

\begin{restatable}{proposition}{rationalsubspace}\label{rationalsubspaceref}
The affine hull and the affine hull of the recession cone of a rationally unbounded set are rational subspaces.
\end{restatable}
\proof{\textbf{Proof}}
\label{rationalsubspace:proof}
Let $I\subseteq \Real^d$ be a rationally unbounded set and $A\subseteq \Real^d$ be equal to either the set $I$ or its recession cone $I_\infty$. Let $F=\aff\bra{\aff\bra{A}\cap \mathbb{Z}^d}$ be the maximal rational affine subspace contained in $\aff\bra{A}$ and let $L\subseteq  \Real^d$ be the rational linear subspace parallel to $F$. Let $P:\Real^d\to \Real^d$ be the projection onto the orthogonal complement $L^\perp$ of $L$. Then $P$ is a rational linear transformation and $P\bra{I}$ is unbounded if $L^\perp\neq \set{\bm{0}}$. Furthermore, because $P\bra{I}\subseteq L^\perp$, we have that  $P\bra{I}\cap \mathbb{Z}^d=P\bra{I}_\infty\cap \mathbb{Z}^d=\set{\bm{0}}$. Hence, because $I$ is rationally unbounded we have $L^\perp=\set{\bm{0}}$ and hence $F$ is a rational subspace.
\Halmos\endproof

\begin{figure}
    \centering
    \includegraphics{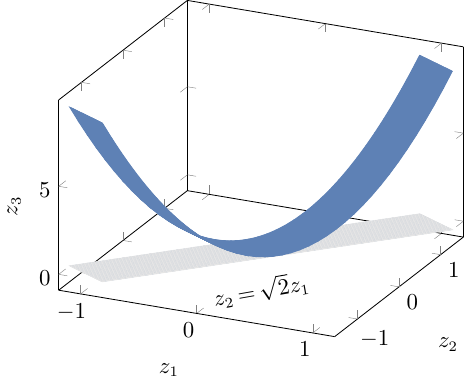}
%
%
%
%
%
%
%
%
%
%

    \caption{In blue, a portion of the lower boundary of the set $I$ in  Example~\ref{almostfinalexampleagain}. The affine hull of the set is $\mathbb{R}^3$ and the affine hull of its recession cone is $\set{0}^2\times \Real$, yet the projection of the set onto $z_1$ and $z_2$ (in grey) has a recession cone that is a irrational subspace. The irrationality is ``revealed" only after projection.}
    \label{fig:hidden_irrational}
\end{figure}
 The following example shows that the property described in \autoref{rationalsubspaceref} does not fully characterize rationally unbounded sets. It also demonstrates why it is necessary to consider all \textit{rational affine images} of index sets in \autoref{rationaldefc}, rather than imposing a simpler condition  in  like $ \mathbb{Z}^{d'} \cap (I_\infty \setminus \set{\bm{0}}) \not = \emptyset$ (note $I_\infty$ in place of $I'_\infty$).
 \begin{example}\label{almostfinalexampleagain}
For $x\in \mathbb{R}$ let $f(x)=x-\lfloor x \rfloor$, and let  $S=\set{x\in \mathbb{Z}: f\bra{x\sqrt{2} }\notin \bra{0.4, 1-0.4} }$ be the non-periodic set from~\eqref{badset2}.  An alternative MICP formulation of $S$ is induced by
\[
M=\set{\bra{x,\bz}\in \mathbb{R}\times  \Real^3\,:\,
 z_1=x,\quad -0.4\leq z_2 - \sqrt{2} z_1 \leq     0.4,\quad  \bra{z_1+\sqrt{2}z_2}^2\leq z_3}
 \]
  and  the index set of $M$ is the set  $I=\set{\bz\in \Real^3\,:\, -0.4\leq z_2 - \sqrt{2} z_1 \leq     0.4,\quad  \bra{z_1+\sqrt{2}z_2}^2\leq z_3 }$ depicted in \autoref{fig:hidden_irrational}. We can check that $I_\infty=\set{0}^2\times \Real_+$, $\mathbb{Z}^3 \cap ({I}_\infty \setminus \set{\bm{0}})\neq\emptyset$,  $\spann\bra{I}=\Real^3$ and  $\spann\bra{I_\infty}=\spann\bra{\set{\bm{e}(3)}}$. However, \blue{$\bra{\proj_{z_1,z_2}\bra{I}}_\infty=\set{\bz\in \Real^2_+\,:\, z_2=\sqrt{2} z_1}$} hence $I$ is not rationally unbounded because $\proj_{z_1,z_2}\bra{I}$ is an image of $I$ under a rational affine mapping, namely $g(z_1, z_2, z_3) = (z_1, z_2)$. \Halmos
\end{example}

Our main technical result in   \autoref{periodictheoref} below uses Definitions~\ref{rationaldefc} and \ref{rationaldefcformulation}  to prove \blue{that, under certain geometric conditions, rational MICP-R sets are finite unions of sets that are either \{compact and convex\} or \{closed and periodic\}.}
 \autoref{periodictheoref} differs from \autoref{boundedorperiodic:coro} in that it proves the rational MICP-R sets can be \emph{finite unions} of closed periodic sets. The need for this union stems from the fact that (as shown in \autoref{basicproppropref} below), the class of rational MICP-R sets is closed under finite unions. Hence, we cannot expect a rational MICP-R set to have a unique direction of periodicity (i.e. consider the union of rational MICP-R sets that have unique and different directions of periodicity).

   \begin{restatable}{theorem}{periodictheo}\label{periodictheoref}
  \blue{ Let $S$ be a closed and rational MICP-R set. If there exists a uniform upper bound on the diameter of any convex subset of $S$, then there exist compact convex sets $\set{C_j}_{j=1}^{k_1}$ and closed periodic sets $\set{S_i}_{i=1}^{k_2}$ such that
  \begin{equation}\label{eq:periodictheoref}
      S = \left( \bigcup_{j=1}^{k_1} C_j\right) \cup \left(\bigcup_{i=1}^{k_2} S_i\right).
  \end{equation}
  }
\end{restatable}
\pointtoproof{periodictheo:proof}

\blue{An interesting contrast between \ref{contrast1} condition \eqref{eq:boundedrep} from  \autoref{boundedprop}, and \ref{contrast2} condition \eqref{eq:periodictheoref} from \autoref{periodictheoref} is that they respectively characterize $S$ as
\begin{enumerate}[label=\emph{\alph*})]
    \item\label{contrast1} The \emph{Minkowski sum} of a finite union of compact convex sets and a \emph{single} closed periodic set.
    \item\label{contrast2} The \emph{union} of a finite union of compact convex sets and a \emph{finite union} of closed periodic sets.
\end{enumerate}

In particular, we cannot truly divide rational MICP-R sets into bounded sets and periodic sets as in \autoref{boundedorperiodic:coro}. Then, to show that the pathological set from \eqref{badset2} is not rational MICP-R, it is not sufficient to note that it is both  unbounded and non-periodic. }
\blue{Fortunately, because the non-periodic set from \eqref{badset2} is a subset of the integer numbers, it is a countable union of sets with diameter equal to zero and hence satisfies the diameter assumption for \autoref{periodictheoref}. We can then conclude that the set from  \eqref{badset2} is not rational MICP-R by showing that it does not satisfy \eqref{eq:periodictheoref}. To achieve this we first use the following version of  Kronecker's Approximation Theorem to show that the set from  \eqref{badset2}  fails to contain any periodic subsets.}
\begin{theorem}[Theorems 438 and 439 in \cite{hardy1979introduction}]\label{Kronecker}
  Let $\alpha\in \Real$, $\varepsilon>0$ and $N\in \mathbb{N}$ and $\theta\in \mathbb{R}\setminus \mathbb{Q}$. Then, there exist $n,p\in\mathbb{Z}$ such that $n>N$ and
  \[\abs{n\theta -p -\alpha}<\varepsilon.\]
  In particular, if $f(x)=x-\lfloor x\rfloor $ and $\theta$ is irrational, then $S\bra{\theta}=\set{ f\bra{\theta z}\,:\, z\in \mathbb{N}}$ is dense in $(0,1)$.
\end{theorem}

\begin{lemma}\label{exampleslemmaperiodicref:lemma}
Let $S=\set{x\in \mathbb{Z}: f\bra{x\sqrt{2} }\notin \bra{0.4, 1-0.4} }$ be the non-periodic set from~\eqref{badset2} \blue{and $T\subseteq S$}. Then the set \blue{$T$} is not periodic.
\end{lemma}
\proof{\textbf{Proof}}
\blue{Assume for a contradiction that $S$ has a periodic subset. Then there exist}
 some $b\in S$ and $a\in \mathbb{Z}$ \blue{such that}  $a \lambda +b \in S$ for all $\lambda\in \mathbb{N}$. Let $\varepsilon=0.01$ so that $(0.42-\varepsilon,0.42+\varepsilon)\subset(0.4, 1-0.4)$.
  By \autoref{Kronecker} for $\theta=\sqrt{2}a$, $\alpha=0.42-\sqrt{2}b$, $\varepsilon=0.01$ and $N=1$ we have that there exist $n>1$ and $p\in \mathbb{Z}$ such that
 \[\abs{\sqrt{2}\bra{an+b}-\bra{p+0.42}}=\abs{n\sqrt{2}a -p -(0.42-\sqrt{2}b)}<\varepsilon=0.01.\]
 This implies that $f\bra{\sqrt{2}\bra{an+b}}\in(0.42-\varepsilon,0.42+\varepsilon)\subset(0.4, 1-0.4)$, which contradicts $a \lambda +b \in S$ for all $\lambda\in \mathbb{N}$.
  \Halmos\endproof

\blue{Finally, we use  the discreteness and cardinality of the  set from \eqref{badset2} to reach our desired conclusion.
\begin{corollary}\label{nonperiodicisnotmicpr}
Let $S=\set{x\in \mathbb{Z}: f\bra{x\sqrt{2} }\notin \bra{0.4, 1-0.4} }$ be the non-periodic set from~\eqref{badset2}. $S$ is not rational MICP-R.
\end{corollary}
\proof{\textbf{Proof}}
Assume for a contradiction that  $S$ is rational MICP-R. All convex subsets of $S$ are elements of $\mathbb{Z}$ with diameter equal to zero, so \autoref{periodictheoref} is applicable and hence $S$ is a finite union of sets that are either periodic or convex.   \autoref{exampleslemmaperiodicref:lemma} further implies that $S$ can only be a finite union of convex sets. However, this contradicts the fact that $S$ is a countably infinite subset of the integer numbers.
 \Halmos\endproof

Finally, using Corollary 1.4 in \cite{hemmecke2007representation} we can check that  $S=\set{\bx\in \mathbb{Z}^3\,:\, \norm{\bra{x_1,x_2}}_2\leq x_3}$ is an example of a rational MICP-R set that is periodic, but does not satisfy  \eqref{eq:boundedrep}.
 }

\subsection{Basic properties of rational MICP-R sets}
\blue{We now present a proposition that summarizes which operations preserve the different classes of representability. To prove this proposition we need the following auxiliary lemma. }
\begin{lemma}\label{lem:unbdproduct}
If $C_1 \subseteq \mathbb{R}^{n_1}$ and $C_2 \subseteq \mathbb{R}^{n_2}$ are  rationally unbounded sets, then $C_1 \times C_2$ is  rationally unbounded.
\end{lemma}
\proof{\textbf{Proof}}
Let $C=C_1\times C_2=\set{\bra{\bx,\by}\,:\, \bx\in C_1, \quad \by \in C_2}$, and let $\mathcal{R}:\Real^{n_1+n_2}\to \Real^m$ be a rational affine transformation. Then $\mathcal{R}\bra{C}=\mathcal{R}_1\bra{C_1}+\mathcal{R}_2\bra{C_2}$ where $\mathcal{R}_1\bra{C_1}=\mathcal{R}\bra{C_1\times \set{0}^{n_2}}-\mathcal{R}\bra{\set{0}^{n_1+n_2}}/2$ and $\mathcal{R}_2\bra{C_2}=\mathcal{R}\bra{\set{0}^{n_1} \times C_2 }-\mathcal{R}\bra{\set{0}^{n_1+n_2}}/2$. The result follows by noting that $\mathcal{R}_1$ and $\mathcal{R}_2$ are rational affine transformations and $\mathcal{R}\bra{C}_\infty=\mathcal{R}_1\bra{C_1}_\infty+\mathcal{R}_2\bra{C_2}_\infty$.
\Halmos\endproof

\begin{restatable}{proposition}{basicpropprop}\label{basicproppropref}
 The classes of rational MILP-R, binary MICP-R, MICP-R, and rational MICP-R sets are closed under rational affine transformations, finite Cartesian products, and Minkowski sums.

 The classes of rational MILP-R, binary MICP-R, and MICP-R sets are closed under finite intersection, but the class of rational MICP-R sets is not closed under finite intersection.

The classes of binary MICP-R, MICP-R, and rational MICP-R sets are closed under finite unions, but the class of rational MILP-R sets is not closed under finite unions.

\end{restatable}
\proof{\textbf{Proof}}
\label{basicpropprop:proof}
\blue{Without loss of generality, we may restrict to two sets $S_1,S_2\subseteq \Real^n$. We may also assume that for each $i\in \sidx{2}$ there exists a closed convex set $M_i\subseteq \Real^{n+p+d}$ such that $(M_i,d)$ induces an MICP formulation of $S_i$ and whose index set is $I_i\subseteq \Real^d$. Finally, we may restrict to a rational affine transformation $\mathcal{R}:\Real^{n}\to \Real^m$.

For the rational affine transformation, let $M\subseteq \Real^{m+(n+p)+d}$ in variables $\bx\in \Real^m$, $\by=\bra{\by^{y},\by^{x}}\in \Real^{p}\times \Real^n$ and $\bz\in \Real^{d}$ be the closed convex set given by
\[
M=\set{\bra{\bx,\by,\bz}\in \Real^{m+(n+p)+d}\,:\, \begin{alignedat}{4} \bra{\by^{x},\by^{y},\bz}&\in M_1,&\quad \mathcal{R}\bra{\by^x}&= \bx\end{alignedat}}.
\]
Then $(M,d)$ induces an MICP formulation of $\mathcal{R}\bra{S_1}$ and the index set of the formulation induced by $(M,d)$ is $I=I_1$. For the Cartesian product, let $M\subseteq \Real^{2n+2(n+p)+2d}$ in variables $\bx\in \Real^{2n}$, $\by=\bra{\by^{y,1},\by^{x,1},\by^{y,2},\by^{x,2}}\in \Real^{p}\times \Real^n\times \Real^{p}\times \Real^n$ and $\bz=\bra{\bz^{z,1},\bz^{z,2}}\in \Real^{d}\times \Real^{d}$ be the closed convex set given by
\[
M=\set{\bra{\bx,\by,\bz}\in \Real^{2n+2(p+n)+2d}\,:\, \begin{alignedat}{4} \bra{\by^{x,1},\by^{y,1},\bz^{z,1}}&\in M_1,&\quad\bra{\by^{x,2},\by^{y,2},\bz^{z,2}}&\in M_2,\quad
\bx=\bra{\by^{x,1},\by^{x,2}}\end{alignedat}}.
\]
Then $(M,2d)$ induces an MICP formulation of $S_1\times S_2$, and the index set of the formulation induced by $(M,2d)$ is $I=I_1\times I_2$. For the Minkowski sum, let $M\subseteq \Real^{n+2(n+p)+2d}$ in variables $\bx\in \Real^{n}$, $\by=\bra{\by^{y,1},\by^{x,1},\by^{y,2},\by^{x,2}}\in \Real^{p}\times \Real^n\times \Real^{p}\times \Real^n$ and $\bz=\bra{\bz^{z,1},\bz^{z,2}}\in \Real^{d}\times \Real^{d}$ be the closed convex set given by
\[
M=\set{\bra{\bx,\by,\bz}\in \Real^{n+2(p+n)+2d}\,:\, \begin{alignedat}{4} \bra{\by^{x,1},\by^{y,1},\bz^{z,1}}&\in M_1,&\quad\bra{\by^{x,2},\by^{y,2},\bz^{z,2}}&\in M_2,\quad
\bx=\by^{x,1}+\by^{x,2}\end{alignedat}}.
\]
Then $(M,2d)$ induces an MICP formulation of $S_1\times S_2$, and the index set of the formulation induced by  $(M,2d)$ is $I=I_1\times I_2$.
}
In addition (using \autoref{lem:unbdproduct} for rational MICP-R), we have that the formulations \blue{associated with the rational affine transformation, Cartesian product, and Minkowski sum} preserve the additional properties that make the formulations rational MILP-R, binary MICP-R, and rational MICP-R.

\blue{For the intersection, let $M\subseteq \Real^{n+2p+2d}$ in variables $\bx\in \Real^{n}$, $\by=\bra{\by^{y,1},\by^{y,2}}\in \Real^{p}\times \Real^{p}$ and $\bz=\bra{\bz^{z,1},\bz^{z,2}}\in \Real^{d}\times \Real^{d}$ be the closed convex set given by
\[
M=\set{\bra{\bx,\by,\bz}\in \Real^{n+2p+2d}\,:\, \begin{alignedat}{4} \bra{\bx,\by^{y,1},\bz^{z,1}}&\in M_1,&\quad\bra{\bx,\by^{y,2},\bz^{z,2}}&\in M_2\end{alignedat}}.
\]
Then  $(M,2d)$ induces an MICP formulation of $S_1\cap S_2$. This}
formulation preserves the additional properties that make the formulations rational MILP-R or binary MICP-R. However, it is not clear if the properties for rational MICP-R are preserved. To check that they may indeed not be preserved, consider the sets
\[M_1=\set{\bra{\bm{x},\bz}\in \mathbb{R}^2\times \mathbb{R}^2\,:\, \bz=\bx,\quad z_2 - \sqrt{2} z_1 \leq     0.4}
\]
and
\[M_2=\set{\bra{\bm{x},\bz}\in \mathbb{R}^2\times \mathbb{R}^2\,:\,  \bz=\bx,\quad \sqrt{2} z_1 - z_2  \leq     0.4}.
\]
We can check that $\proj_{\bz}\bra{M_1}$ and $\proj_{\bz}\bra{M_1}$  are rationally unbounded so ${S}_i=\proj_{\bx}\bra{M_i}$ is rational MICP-R for each $i\in \set{1,2}$. However,
\[\proj_{x_1}\bra{{S}_1\cap{S}_2}=\set{x\in \mathbb{Z}: f\bra{x\sqrt{2} }\notin \bra{0.4, 1-0.4} },\]
where  $f(x)=x-\lfloor x \rfloor$, is exactly the \blue{non-periodic set from \eqref{badset2}, which is not rational MICP-R by \autoref{nonperiodicisnotmicpr}}.
Then ${S}_1\cap{S}_2$ also fails to be rational MICP-R, because the orthogonal projection is a rational affine transformation, and we have already proven that such transformations preserve rational MICP-R representability.

\blue{Finally, for the union operation let $M\subseteq \Real^{n+(2p+2n+1)+(2d+1)}$ in variables $\bx\in \Real^n$, $\by=\bra{\by^{y,1},\by^{x,1},\by^{y,2},\by^{x,2},y_0}\in \Real^{p}\times \Real^n\times \Real^{p}\times \Real^n\times \Real$ and $\bz=\bra{\bz^{z,1},\bz^{z,2},z_0}\in \Real^{d}\times \Real^{d}\times \Real $ be the closed convex set given by
\[
M=\set{\bra{\bx,\by,\bz}\in \Real^{n+(2p+2n+1)+(2d+1)}\,:\, \begin{alignedat}{4} \bra{\by^{x,1},\by^{y,1},\bz^{z,1}}&\in M_1,&\quad\bra{\by^{x,2},\by^{y,2},\bz^{z,2}}&\in M_2\\
\norm{\bx-\by^{x,1}}_2^2 &\leq y_0 z_0 ,&\quad \norm{\bx-\by^{x,2}}_2^2 &\leq y_0 (1-z_0)\\
y_0&\geq 0,&\quad 0\leq z_0&\leq 1\end{alignedat}}.
\]
If $\bra{\bx,\by,\bz}\in M\cap \bra{\Real^{3n+p_1+p_2}\times \mathbb{Z}^{d_1+d_2+1}}$ is such that $z_0=1$, then $\bx=\by^{x,2}$ and $\bx\in S_2$. Similarly, if instead $z_0=0$, then $\bx=\by^{x,1}$ and $\bx\in S_1$. Hence,  $(M,2d+1)$ induces an MICP formulation of $S_1\cup S_2$. In addition, if $I_1$ and $I_2$ are the index sets of $M_1$ and $M_2$ respectively, then the index set of the formulation induced by $(M,2d+1)$ is equal to $I_1\times I_2\times [0,1]$. By \autoref{lem:unbdproduct}, we have }
 that the additional properties are preserved for binary MICP-R and rational MICP-R. \blue{However, $M$ includes two additional non-linear inequalities so it} does not preserve the properties for rational MILP-R. The fact that this cannot be resolved with another version of the formulation follows by noting that a rational polyhedron is rational MILP-R, but the union of two rational polyhedra with different recession cones is not rational MILP-R  by \autoref{JLTheo}.
\Halmos\endproof

We end this section with the following simple connection between MILP-R and MICP-R sets.

\begin{lemma}\label{rationalmicpmilp}
Any rational MILP-R set is also rational MICP-R.
\end{lemma}
\proof{\textbf{Proof}}
Let $M\subseteq \Real^{n+p+d}$ be a rational polyhedron that induces a rational MILP formulation and whose index set is  $I\subseteq \Real^d$. Then, $I$ is the projection of a rational polyhedron so $I$ and any of its rational affine images is also a rational polyhedron. Hence, $I$ is rationally unbounded and $M$ induces a rational MICP formulaiton.
\Halmos\endproof

\subsection{Applications of \autoref{periodictheoref}}

\subsubsection{Compact sets}
We can also use \autoref{periodictheoref} to fully characterize rational MICP-R sets that are compact sets.

\begin{restatable}{proposition}{compactfinitetheo}\label{thm:compactfinite}
Let $S \subseteq \mathbb{R}^n$ be a compact set. Then the following are equivalent:
\begin{enumerate}
    \item[(a)]\label{thmcompacta} $S$ is rational MICP-R.
    \item[(b)]\label{thmcompactb} $S$ is a finite union of compact convex sets.
    \item[(c)]\label{thmcompactc} \blue{$S$ is pure binary MICP-R.}
     \item[(d)]\label{thmcompactd} $S$ is binary MICP-R.
\end{enumerate}
\end{restatable}
\proof{\textbf{Proof}}
\label{compactfinitetheo:proof}
(\emph{(a)}$\Rightarrow$\emph{(b)}):
Compactness of $S$ implies that there exists a uniform upper bound on the diameter of any convex subset of $S$. Hence, we may apply  \autoref{periodictheoref} to conclude that $S$ is a finite union of sets that are either compact and convex, or closed and periodic. However, periodic sets must be unbounded; therefore $S$ is a finite union of compact convex sets.

(\emph{(b)}$\Rightarrow$\emph{(c)}): \blue{Direct from \autoref{lem:balasform_jeroslow}.}

(\emph{(c)}$\Rightarrow$\emph{(d)}): Direct from the definitions of binary and pure binary MICP-R.

(\emph{(d)}$\Rightarrow$\emph{(a)}): \blue{Let $M\subseteq \Real^{n+p+d}$ be a closed convex set such that $(M, d)$ induces a binary MICP formulation of $S$ and let $M'=M\cap \bra{\Real^{n+p}\cap [0,1]^d}$. Because $\proj_{\bz}\bra{M\cap \bra{\mathbb{R}^{n+p} \times \mathbb{Z}^d}}\subseteq\set{0,1}^d$ we have that $(M', d)$ also induces a binary MICP formulation of $S$. The index set $I'=\proj_{\bz}\bra{M'}$ of the formulation induced by $(M', d)$ is such that $I'\subset [0,1]^d$ and hence it is bounded. } Therefore every rational affine transformation of $I'$ is bounded, so $I'$ is rationally unbounded and hence $(M', d)$ induces a rational MICP formulation of $S$.
\Halmos\endproof
\autoref{thm:compactfinite} \blue{provides} a simpler alternative to the midpoint lemma to prove that sets are not \emph{rational} MICP-R. For instance,
 \autoref{thm:compactfinite} trivially implies that the  set $\{1/n : n \in \mathbb{N}  \} \cup \{0\}$, the spherical shell $\set{\bm{x}\in \Real^n\,:1\leq \norm{\bm{x}}_2\leq 2\,}$, and the set of rank 1 contained in some compact domain are not rational MICP-R, \blue{because neither satisfies condition (b)}\footnote{The midpoint lemma also can be used for the first set by noting that for any $n,m\in \mathbb{N}$ with $n\neq m$ there is no $k\in \mathbb{N}$ such that $(2^{-n}+2^{-m})/2=1/k$.}. Another interesting interpretation of \autoref{thm:compactfinite} is that for closed rational MICP-R sets,  unbounded integer variables are needed \textit{only} for modeling unbounded sets.

\subsubsection{Subsets of the natural numbers}

The following lemma gives a characterization of rational MILP-R subsets of the natural numbers, providing context for the corresponding MICP-R characterization given by \autoref{thm:naturalrep}.

\begin{restatable}{lemma}{natmilpper}\label{thm:naturalreplem}
Let $S \subseteq \mathbb{N}$ with $|S| = \infty$. Then the following are equivalent:
\begin{enumerate}
\item[(a)]  $S$ is rational MILP-R.
    \item[(b)] $S$ is periodic.
    \item[(c)] There exist  a non-empty finite set $S_0\subset \mathbb{N}$ and $u_0\in \mathbb{N}\setminus \set{0}$ such that $S = S_0 + \operatorname{intcone}(\set{u_0})$.

\end{enumerate}
\end{restatable}
\proof{\textbf{Proof}}
\label{natmilpper:proof}
(\emph{(b)}$\Rightarrow$\emph{(c)}): Suppose $S$ is periodic and let $ u_0 \in \Real\setminus \set{0}$ be such that
\begin{equation}\label{naturalperiod}
x + \lambda u_0 \in S \quad  \forall x \in S,\quad \lambda \in \mathbb{Z}_+.
\end{equation}
We clearly also have that $u_0\in \mathbb{N}$. For each $i \in \sidx{u_0-1}$ let  $\lambda_i =\min\set{\lambda \in \mathbb{Z}_+\,:\, i+\lambda u_0 \in S}$ with the convention that $\lambda_i=\infty$ if $i+\lambda u_0 \not \in S$ for all  $\lambda\in \mathbb{Z}_+$. If   $\lambda_i<\infty$, then by \eqref{naturalperiod} we also have that  $(i+\lambda_i u_0)+\lambda u_0= i+(\lambda_i+\lambda)u_0 \in S$ for all $\lambda \in \mathbb{Z}_+$. Furthermore, every integer in $S$ with remainder $i$ modulo $u_0$ is of this form. Hence, $S = S_0 + \operatorname{intcone}(\set{u_0})$ for  $S_0=\set{i+\lambda_i u_0\,:\, i\in \sidx{u_0-1}, \lambda_i<\infty}$.

(\emph{(c)}$\Rightarrow$\emph{(a)}): Follow directly from \autoref{JLTheo} as \eqref{eq:jlregular} holds with $U=\set{u_0}$.

(\emph{(a)}$\Rightarrow$\emph{(b)}):  If $S$ is rational MILP-R, then by \autoref{JLTheo} we have that $S = S_0 + \operatorname{intcone}(R)$ for some finite sets $S_0 \subset \mathbb{N}$ and $R \subset \mathbb{N}$. Because $\abs{S}=\infty$ there exists $u_0\in R\setminus\set{0}$. That $u_0$ satisfies \eqref{naturalperiod} and hence $S$ is periodic.
\Halmos\endproof

\autoref{thm:naturalreplem} shows that for infinite subsets of the natural numbers, Jeroslow's rational MILP-R characterization \eqref{eq:jlregular} from \autoref{JLTheo} holds for a set $U$ that contains a single non-zero element. The following proposition shows that the corresponding rational MICP-R characterization is nearly identical: the equivalence between (a) and (d) in \autoref{thm:naturalrep} states that for any  infinite subset of the natural numbers that is rational MICP-R,  there is a rational MILP-R set that differs by at most finitely many points.

\begin{restatable}{proposition}{naturalreptheo}\label{thm:naturalrep}
Let $S \subseteq \mathbb{N}$ with $|S| = \infty$. Then the following are equivalent:
\begin{enumerate}
    \item[(a)] $S$ is rational MICP-R.
    \item[(b)] There exists a finite set $S_1\subset \mathbb{N}$ and a non-empty periodic set $S_2\subset \mathbb{N}$ such that $S = S_1 \cup S_2$.
        \item[(c)]  There exists a finite set $S_1\subset \mathbb{N}$, a non-empty finite set $S_0\subset \mathbb{N}$ and $u_0\in \mathbb{N}\setminus \set{0}$ such that \mbox{$S = S_1 \cup\bra{S_0 + \operatorname{intcone}(\set{u_0})}$}.
    \item[(d)]\label{thm:naturalrep:d}  There exists a finite set $S_1\subset \mathbb{N}$ and an infinite rational MILP-R set $S_2\subset \mathbb{N}$ such that $S = S_1 \cup S_2$.
\end{enumerate}
\end{restatable}
\proof{\textbf{Proof}}
\label{naturalreptheo:proof}
(\emph{(a)}$\Rightarrow$\emph{(b)}): All convex subsets of $S$ are elements of $\mathbb{N}$ with diameter equal to zero, so  \autoref{periodictheoref} is applicable and hence $S$ is a finite union of sets $\set{\tilde{S}_i}_{i=1}^k$ that are either convex or periodic. Let $S_1=\bigcup_{i\in \sidx{k}\,:\, \abs{\tilde{S}_i}=1} \tilde{S}_i$ and $S_2=\bigcup_{i\in \sidx{k}\,:\, \abs{\tilde{S}_i}=\infty} \tilde{S}_i$. The result follows by noting that $\abs{S}=\infty$ implies $S_2\neq \emptyset$ and that a finite union of periodic sets in $\mathbb{N}$ is itself periodic (e.g. if $u_i\in \mathbb{N}$ is such that $\tilde{S}_i+\operatorname{intcone}(\set{u_i})\subseteq \tilde{S}_i$ for each $i\in I$,
then $\bigcup_{i\in I}\tilde{S}_i+\operatorname{intcone}\bra{ \set{\prod_{i\in I}u_i} }\subseteq \bigcup_{i\in I}\tilde{S}_i$).

(\emph{(b)}$\Leftrightarrow$\emph{(c)}$\Leftrightarrow$\emph{(d)}): Follows by \autoref{thm:naturalreplem} applied to $S_2$.

(\emph{(d)}$\Rightarrow$\emph{(a)}): Both $S_1$ and $S_2$ are rational MICP-R so the result follows from \autoref{basicproppropref}.
\Halmos\endproof

\section{Proof of \autoref{periodictheoref}}\label{secondproofsection}

\subsection{Roadmap of the proof}
\begin{figure}[htpb]
  \includegraphics{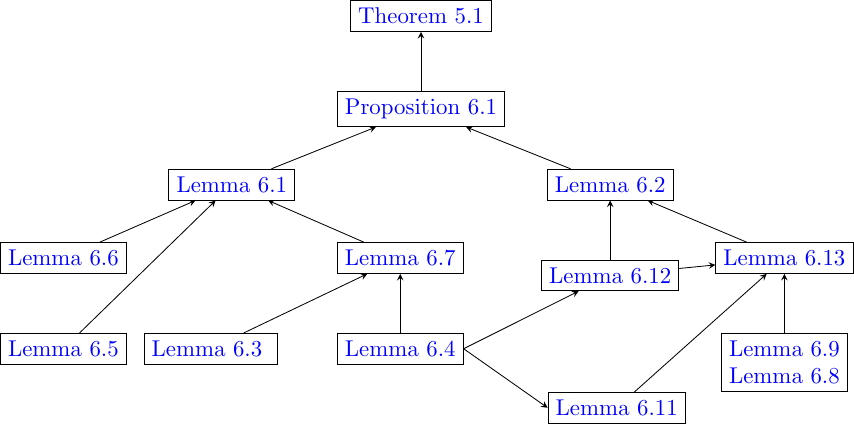}
%
%
%
%
%
%
%
%
%
%
%
%
%
%
	\caption{Causal connection between the Lemmas and Propositions used to prove \autoref{periodictheoref}.}\label{lemmapropdiagramfig}
\end{figure}

In this section, we provide the proof for \autoref{periodictheoref}. The proof is organized as follows. In the next subsection, we provide three important technical definitions used throughout the proof. In the following subsection, we state two key lemmas: \autoref{interiorpointind} and \autoref{lem:generaldecompsimple}, deferring their proof for later sections. Following that, in the same subsection, we use the lemmas to prove \autoref{periodictheoref}. Then, in the subsequent two subsections we prove \autoref{interiorpointind} and \autoref{lem:generaldecompsimple}. To alleviate some of the complexity of the proof of \autoref{periodictheoref}, we offer the following diagram for the reader's convenience. The diagram in Figure~\ref{lemmapropdiagramfig} shows the causal connection between the different Lemmas and Propositions presented in this Section.

\subsection{Key Technical Definitions}
In this subsection, we provide three technical definitions of instrumental importance for the proof of \autoref{periodictheoref}.

We first define the notion of  MICP diameter for an MICP formulation.
\begin{definition}\label{uniformbounddef}
Let $M\subseteq \Real^{n+p+d}$ be a closed convex set inducing an MICP formulation with index set $I\subseteq\Real^d$ and $\bz$-projected sets $\set{A_{\bz}}_{\bz\in I}$. We define the \textbf{MICP diameter} of the MICP formulation induced by \blue{the pair $(M,d)$} as
\[D_{\text{MICP}}\bra{M, d}=\sup_{\bz\in I\cap\mathbb{Z}^d} D\bra{A_{\bz}},\]
where $D(\cdot)$ is the diameter (see \autoref{diameterdef}). \blue{The argument $d$ is omitted when implied by the context.}
\end{definition}
Notice that if $M$ induces an MICP formulation of $S$, then for each $\bz\in I\cap \mathbb{Z}^d$ the $\bz$-projected set $A_{\bz}$ is a convex subsets of $S$. Then the geometric assumption of \autoref{periodictheoref} implies that $D_{\text{MICP}}\bra{M}< \infty$.

We further define the notion of a nearly periodic MICP formulation.
\begin{definition}\label{nearlyperiodicdef}We say an MICP formulation induced by $M \subseteq \mathbb{R}^{n+p+d}$ with index set $I$ and $\bz$-projected sets   $\set{A_{\bz}}_{\bz\in I}$, is  \textbf{nearly periodic}  if there exist $\bm{u}\in \Real^n\setminus\set{\bm{0}}$ and $\bm{r}\in \mathbb{Z}^d \cap (I_\infty \setminus \set{\bm{0}})$ such that
    \begin{equation}\label{nearlyperiodiceq}\bx+\lambda \bm{u} \in \operatorname{cl}\bra{A_{\bz+\lambda \bm{r}}} \quad \forall \; \lambda\in \mathbb{Z}_+, \;\bm{z} \in\mathbb{Z}^d\cap I, \;\bx\in A_{{\bz}}.
    \end{equation}
\end{definition}
Notice that if $S$ has a nearly periodic MICP formulation
with $\bm{u}$ satisfying \eqref{nearlyperiodiceq}, then  $\bx + \lambda \bm{u} \in \operatorname{cl}\bra{S}$ for all $ \lambda\in \mathbb{Z}_+$ and $\bx\in S$. Hence, a closed set with a nearly periodic MICP formulation is periodic according to \autoref{periodicdef}. The additional technicalities in \autoref{nearlyperiodicdef} are primarily artifacts of our proofs, which are unfortunately hard to avoid. For instance, \autoref{nearlyperiodicclosureex} shows that replacing $\operatorname{cl}\bra{A_{\bz+\lambda \bm{r}}}$ by $A_{\bz+\lambda \bm{r}}$ in \eqref{nearlyperiodiceq} may not be possible even if $S$ is a closed convex set.

Finally, to deal with sets whose affine hull may fail to be a rational subspace we define an ad-hoc notion for a rational affine hull as is common in MICP theory  \cite{The-Chvatal-Gomory-Closure-of-a-Strictly,On-the-Chvatal-Gomory-Closure-of-a-Compact-FULL,Dey2013,diego2011maximal,moran2016closedness}.
\begin{definition}
For $I\subseteq \Real^d$  with $I\cap \mathbb{Z}^d\neq \emptyset$ we let $\aff_{\mathbb{Z}}\bra{I}=\aff\bra{I\cap \mathbb{Z}^d} $ be the \textbf{rational affine hull} of $I$. We also let the intersection of $I$ and its rational affine hull be $I_{\mathbb{Z}}=I\cap \aff_{\mathbb{Z}}\bra{I}$.
\end{definition}

\subsection{Statement of two key lemmas and derivation of \autoref{periodictheoref}}

First, \autoref{interiorpointind} decomposes rational MICP-R sets into a binary MICP-R set and sets that have rational MICP-R formulations with an additional technical condition on their index sets.
\begin{restatable}{lemma}{interiorpointindprop}\label{interiorpointind}
Let $M\subseteq \Real^{n+p+d}$ be a closed convex set \blue{such that the pair $(M,d)$ induces} a rational MICP formulation of a non-empty set $S \subseteq \mathbb{R}^n$. Then, there exist a binary MICP-R set $S_0\subseteq \Real^n$ and closed convex sets $\set{{M}_i}_{i=1}^k$ \blue{such that for each $i \in \sidx{k}$, the pair $(M_i, d)$ induces a rational MICP formulation of a non-empty set $S_i$, and together these sets decompose $S$ as the union} $$S=S_0\cup\bigcup_{i=1}^k {S}_i.$$ For each $i \in \sidx{k}$ the following properties additionally hold:
\begin{subequations}\label{interiorpointindprop:c13}
\begin{align}
        \label{interiorpointindprop:c2} D_{\text{MICP}}\bra{{M}_i}&\leq D_{\text{MICP}}\bra{M},\\
    \label{interiorpointindprop:c3} \bra{{I}_i}_{\mathbb{Z}} &\text{ has an integer point in its relative interior},
\end{align}
\end{subequations}
where ${I}_i$ is the index set of the formulation induced by $\bra{{M}_i,d}$.
\end{restatable}
\pointtoproof{interiorpointindprop:proof}

Next, \autoref{lem:generaldecompsimple} shows that sets with rational MICP formulations that comply with the  technical index-set condition from \autoref{interiorpointind} satisfy additional properties that we will use for induction on the number of integer variables of the MICP formulation.
\begin{restatable}{lemma}{generaldecompsimpleprop}\label{lem:generaldecompsimple}
Let $M\subseteq \Real^{n+p+d}$ be a closed convex set \blue{such that the pair $(M, d)$ induces} a \blue{rational} MICP formulation of a non-empty set $S \subseteq \mathbb{R}^n$, and let $I\subseteq \Real^d$ be the index set of the formulation induced by $M$.
If $D_{\text{MICP}}\bra{M}<\infty$ and $I_{\mathbb{Z}}$ has an integer point in its relative interior, then one of the following holds
\begin{enumerate}[label=\textnormal{(\alph*)}]
    \item\label{lem:generaldecompc1simple} $S$ is binary MICP-R, or
    \item\label{lem:generaldecompc2simple} $S$ has a nearly periodic MICP formulation, or
    \item\label{lem:generaldecompc3simple} there exists a closed convex set $\tilde{M}$ \blue{such that the pair $(\tilde{M}, d-1)$ induces a} rational MICP formulation of a non-empty set $\tilde{S}$ with the properties
    \begin{subequations}\label{lem:generaldecompc3simplea}
\begin{align}
        \label{lem:generaldecompc3simpleab} D_{\text{MICP}}\bra{\tilde{M}}&\leq D_{\text{MICP}}\bra{{M}},\\
    \label{lem:generaldecompc3simpleac} S\subseteq  \tilde{S}\subseteq \cl\bra{\tilde{S}}&\subseteq \cl\bra{S}.
\end{align}
\end{subequations}
\end{enumerate}
\end{restatable}
\pointtoproof{generaldecompsimpleprop:proof}

We now use the two lemmas to establish the following technical proposition as the final intermediate step towards proving \autoref{periodictheoref}. We note that the need for  \autoref{lem:generaldecompsimple} and \autoref{peroidictheomoregeneral} to consider possibly non-closed MICP-R sets is primarily an artifact of the use of \autoref{interiorpointind} in our proofs. Unfortunately, \autoref{nearlyperiodicclosureex} again shows this technicality is also hard to avoid.

\begin{proposition}\label{peroidictheomoregeneral}
Let $M\subseteq \Real^{n+p+d}$ be a closed convex set \blue{such that the pair $(M,d)$ induces} a rational MICP formulation of $S \subseteq \mathbb{R}^n$. If $D_{\text{MICP}}\bra{M}<\infty$, then there exist a binary MICP-R set $S_0$ and sets $\set{T_i}_{i=1}^k$ with nearly periodic MICP formulations such that
\[
S \subseteq S_0 \cup \bigcup_{i=1}^k T_i\subseteq \cl\bra{S}.
\]
\end{proposition}

\proof{\textbf{Proof}}
\label{peroidictheomoregeneral:proof}
We prove the statement by induction on $d$. The base case $d=0$ follows directly because in that case $S$ is a convex set.

Now assume the result for $d-1$. \blue{By \autoref{interiorpointind},  there exists a binary MICP-R set $S_0\subseteq \Real^n$ and closed convex sets $\set{{M}_i}_{i=1}^k$ \blue{such that for each $i \in \sidx{k}$, the pair $(M_i,d)$ induces a rational MICP formulation of a non-empty set $S_i$, and together these sets decompose $S$ as the union} $S=S_0\cup\bigcup_{i=1}^k {S}_i$. \blue{Additionally,} for each $i\in \sidx{k}$ we have \eqref{interiorpointindprop:c13}. In particular, because of \eqref{interiorpointindprop:c2}, \eqref{interiorpointindprop:c3}, and $D_{\text{MICP}}\bra{M}<\infty$, we have that $S_i$ satisfies the conditions of \autoref{lem:generaldecompsimple} for all  $i\in \sidx{k}$. Then, for each $j\in \sidx{k}$ we have that $S_i$ satisfies one of the  options \ref{lem:generaldecompc1simple}--\ref{lem:generaldecompc3simple} of \autoref{lem:generaldecompsimple}. We can then partition $\sidx{k}$ into sets $J_a$, $J_b$ and $J_c$ according to which one of this options each $S_j$ satisfies. Because the union of binary MICP-R sets is binary MICP-R by \autoref{basicproppropref} we have that $S_0\cup \bigcup_{j\in J_a} S_j$ is binary MICP-R. Redefine  $S_0$ to be equal to  $S_0\cup \bigcup_{j\in J_a} S_j$, and let $k_1=\abs{J_b}$, $k_2=\abs{J_c}$, $\set{{T}_j}_{j=1}^{k_1}=\set{S_j}_{j\in J_b}$  and  $\set{{R}_i}_{i=1}^{k_2}=\set{S_j}_{j\in J_c}$. Then, $S_0$ is a non-empty binary MICP-R set,  $\set{{T}_j}_{j=1}^{k_1}$ are non-empty sets with nearly periodic MICP formulations and  $\set{{R}_i}_{i=1}^{k_2}$ are non-empty sets such that}
\begin{equation}\label{peroidictheomoregeneral:eq1}
    S=S_0\cup \bigcup_{j=1}^{k_1} {T}_j\cup \bigcup_{i=1}^{k_2} {R}_i
\end{equation}
and for each $i\in \sidx{k_2}$ there exists a closed convex set $\tilde{M}_i$ \blue{such that the pair $(\tilde{M}_i, d-1)$ induces a} rational MICP formulation of a non-empty set $\tilde{R}_i$ with the properties
        \begin{subequations}\label{peroidictheomoregenerallem:generaldecompc3simplea}
\begin{align}
        \label{peroidictheomoregenerallem:generaldecompc3simpleab} D_{\text{MICP}}\bra{\tilde{M}_i}&\leq D_{\text{MICP}}\bra{{M}},\\
    \label{peroidictheomoregenerallem:generaldecompc3simpleac} R_i\subseteq  \tilde{R}_i\subseteq \cl\bra{\tilde{R}_i}&\subseteq \cl\bra{R_i}.
\end{align}
\end{subequations}

Note that  \eqref{peroidictheomoregeneral:eq1} implies that $R_i\subseteq S$ and hence $\cl\bra{R_i}\subseteq \cl\bra{S}$. Hence,  \eqref{peroidictheomoregeneral:eq1} and    \eqref{peroidictheomoregenerallem:generaldecompc3simpleac}  imply
\begin{equation}\label{peroidictheomoregeneral:eq2a}
    S \subseteq S_0\cup \bigcup_{j=1}^{k_1} {T}_j\cup \bigcup_{i=1}^{k_2}   \tilde{R}_i \subseteq S_0\cup \bigcup_{j=1}^{k_1} {T}_j\cup \bigcup_{i=1}^{k_2} \cl\bra{\tilde{R}_i}\subseteq S_0\cup \bigcup_{j=1}^{k_1} {T}_j\cup \bigcup_{i=1}^{k_2} \cl\bra{R_i}\subseteq \cl\bra{S}.
\end{equation}

Furthermore, by the induction hypothesis for each $i\in\sidx{k_2}$ there exists a binary MICP-R set $R_{i,0}$ and sets $\set{T_{i,l}}_{l=1}^{q_i}$ with nearly periodic MICP formulations such that
\begin{equation}\label{peroidictheomoregeneral:eq3a}
   \tilde{R}_i \subseteq R_{i,0} \cup \bigcup_{l=1}^{q_i} T_{i,l}\subseteq \cl\bra{\tilde{R}_i }.
\end{equation}
 Combining \eqref{peroidictheomoregeneral:eq2a} with \eqref{peroidictheomoregeneral:eq3a}  we get
\begin{equation*}
    S\subseteq \bra{S_0\cup  \bigcup_{i=1}^{k_2} R_{i,0}} \cup \bigcup_{j=1}^{k_1}{T}_j \cup \bigcup_{i=1}^{k_2}  \bigcup_{l=1}^{q_i} T_{i,l}\subseteq \cl\bra{S}.
\end{equation*}
Then the result follows for $M$ by noting that from \autoref{basicproppropref}, the set $\bra{S_0\cup  \bigcup_{i=1}^{k_2} R_{j,0}} $ is binary MICP-R.
\Halmos\endproof

We now establish \autoref{periodictheoref} from \autoref{peroidictheomoregeneral}.

\periodictheo*
\proof{\textbf{Proof}}
\label{periodictheo:proof}
Let  $M\subseteq \Real^{n+p+d}$ be a closed convex set inducing a rational MICP formulation of $S$, $I$ be its index set and $\set{A_{\bz}}_{\bz\in I}$ be its $\bz$-projected sets. Then, for each $\bz\in I\cap \mathbb{Z}^d$ the $\bz$-projected set $A_{\bz}$ is a convex subset of $S$, so under the theorem's assumptions we have $D_{\text{MICP}}\bra{M}<\infty$.

Then by \autoref{peroidictheomoregeneral} there exist $k_2 \in \mathbb{Z}_+$, a binary MICP-R set $S_0$, and sets $\set{T_i}_{i=1}^{k_2}$ with nearly periodic MICP formulations such that
\begin{equation}\label{periodictheo:eq1}
    S \subseteq S_0 \cup \bigcup_{i=1}^{k_2} T_i\subseteq \cl\bra{S}.
\end{equation}
However, because $S$ is closed we actually have \begin{equation}\label{periodictheo:eq2}S = S_0 \cup \bigcup_{i=1}^{k_2} T_i. \end{equation}

Fix some $i\in\sidx{k_2}$. Let $M_i \subseteq \mathbb{R}^{n+p_i+d_i}$ be a closed convex set such that $(M_i,d_i)$ induces a nearly periodic MICP formulation of $T_i$, with index set $I_i$ and $\bz$-projected sets $\set{A_{i,\bz}}_{\bz\in I_i}$. That is, $T_i=\bigcup_{\bm{z} \in I_i \cap \mathbb{Z}^{d_i} } A_{i,\bz} \subseteq \mathbb{R}^n$. Since $T_i\subseteq S$ we have that $\operatorname{cl}\bra{A_{i,\bz}}\subseteq S$ because $S$ is closed. Hence, we have that for the set $\tilde{T}_i=\bigcup_{\bm{z} \in I_i \cap \mathbb{Z}^{d_i}} \operatorname{cl}\bra{A_{i\bz}}$ it holds $T_i\subseteq \tilde{T}_i\subseteq S$.

We claim that because $M_i$ induces a nearly periodic formulation of $T_i$, then $\tilde{T}_i$ is periodic. To establish this, note that from the definition of nearly-periodic MICP formulation there exist $\bm{u}(i)\in \Real^n\setminus\set{\bm{0}}$ and $\bm{r}(i)\in \mathbb{Z}^{d_i}\setminus\set{\bm{0}}$ such that
\begin{equation}\label{nearperiodictoperiodic}
\bx+\lambda \bm{u}(i) \in \operatorname{cl}\bra{A_{i,\bz+\lambda \bm{r}(i)}},  \bz+\lambda \bm{r}(i) \in I\cap \mathbb{Z}^d \quad \forall \; \lambda\in \mathbb{Z}_+, \;\bm{z} \in\mathbb{Z}^{d_i}\cap I_i, \;\bx\in A_{i,{\bz}}.
\end{equation}
Let $\bx \in \tilde{T}_i$. That means that for some $\bm{z} \in\mathbb{Z}^{d_i}\cap I_i$ it holds $\bx\in \operatorname{cl}\bra{A_{i,{\bz}}}$. Let $\set{\bx^m}_{m\in\mathbb{N}}$ be such that $\bx^m\in A_{i,{\bz}}$ for all $m\in \mathbb{N}$ and $\lim_{m\to\infty} \bx^m=\bx$. Then, by \eqref{nearperiodictoperiodic},  for all $\lambda\in \mathbb{Z}_+$  we have $\bx^m+\lambda \bm{u}(i) \in \operatorname{cl}\bra{A_{i,\bz+\lambda \bm{r}(i)}}$ from which we can derive $\bx+\lambda \bm{u}(i) \in \operatorname{cl}\bra{A_{\bz+\lambda \bm{r}(i)}} \subseteq  \tilde{T}_i$, establishing the periodicity of $\tilde{T}_i$.

\blue{We additionally claim that $S_i=\operatorname{cl}\bra{\tilde{T}_i}$ is also periodic. To establish this, let
$\set{\bm{x}^m}_{m\in \mathbb{N}}$ be a convergent sequence contained in $\tilde{T}_i$ and $\overline{\bm{x}}=\lim_{m \to \infty} \bm{x}^m$. Because $\tilde{T}_i$ is periodic, there exists $ \bm{u} \in \Real^n\setminus \set{\bm{0}}$ such that   $\bm{x} + \lambda \bm{u} \in \tilde{T}_i$ for all $\bm{x} \in \tilde{T}_i$ and $\lambda \in \mathbb{N}$. In particular, for any $\lambda \in \mathbb{N}$ we have $\bm{x}^m +  \lambda \bm{u} \in \tilde{T}_i$ for all $i\in \mathbb{N}$. Then,  $\overline{\bm{x}}+\lambda \bm{u}= \lim_{m \to \infty} \bra{\bm{x}^m +  \lambda \bm{u}} \in \operatorname{cl}\bra{\tilde{T}_i}=S_i$ and hence $S_i$ is closed and periodic. Finally, because $S$ is closed we also have $T_i\subseteq \tilde{T}_i\subseteq S_i\subseteq S$.
}

Now because $T_i\subseteq  S_i \subseteq S$, we have that \eqref{periodictheo:eq2} implies
\begin{equation}\label{periodictheo:eq33}
   S = S_0 \cup \bigcup_{i=1}^{k_2} S_i.
\end{equation}
 By \autoref{lem:balasform} there exist convex sets $\set{\tilde{C}_i}_{j=1}^{k_1}$ such that $S_0=\bigcup_{j=1}^{k_1} \tilde{C}_j$. \blue{In particular, because $S$ is closed,  for all $j\in \sidx{{k_1}}$ we have  $\tilde{C}_j\subseteq C_j \subseteq S$ where  $C_j=\operatorname{cl}\bra{C_j}$ is a closed convex set. Combining this with \eqref{periodictheo:eq33} and $S_0=\bigcup_{j=1}^{k_1} \tilde{C}_j$ we finally have $S = \bigcup_{j=1}^{k_1} C_j \cup \bigcup_{i=1}^k S_i$.}  \Halmos\endproof

\subsection{Proof of \autoref{interiorpointind}}
\subsubsection{Tools from the geometry of numbers}

We first list several geometric tools we use to prove \autoref{interiorpointind}.

\begin{lemma}\label{intafflemma}
For any convex set $I\subseteq \Real^d$ such that $I\cap \mathbb{Z}^d\neq \emptyset$ we have that $\aff_{\mathbb{Z}}\bra{I}$ is a rational affine subspace,  $I_{\mathbb{Z}}\cap\mathbb{Z}^d=I\cap\mathbb{Z}^d$ and $\aff_{\mathbb{Z}}\bra{I_{\mathbb{Z}}}=\aff_{\mathbb{Z}}\bra{I}=\aff\bra{I_{\mathbb{Z}}}$. Furthermore, if $\mathcal{R}:\Real^d\to \Real^d$ is an invertible rational affine transformation such that $\mathcal{R}\bra{\mathbb{Z}^d}=\mathcal{R}^{-1}\bra{\mathbb{Z}^d}=\mathbb{Z}^d$, then $\mathcal{R}\bra{I_{\mathbb{Z}}}=\mathcal{R}\bra{I}_{\mathbb{Z}}$ and $\mathcal{R}\bra{\aff_{\mathbb{Z}}\bra{I}}=\aff_{\mathbb{Z}}\bra{\mathcal{R}\bra{I}}$.
\end{lemma}
\proof{\textbf{Proof}}
First recall the definitions $\aff_{\mathbb{Z}}\bra{I}=\aff\bra{I\cap \mathbb{Z}^d} $ and $I_{\mathbb{Z}}=I\cap \aff_{\mathbb{Z}}\bra{I}$.
Rationality of the space $\aff_{\mathbb{Z}}\bra{I}$ follows because $I\cap \mathbb{Z}^d$ is a set of rational points. Now, $\mathbb{Z}^d\cap I\subseteq \aff_{\mathbb{Z}}\bra{I}$, so $\mathbb{Z}^d\cap I\cap \aff_{\mathbb{Z}}\bra{I}=\mathbb{Z}^d\cap I$, which proves $I_{\mathbb{Z}}\cap\mathbb{Z}^d=I\cap\mathbb{Z}^d$. This last equation and the definition of $\aff_{\mathbb{Z}}\bra{I}$ prove $\aff_{\mathbb{Z}}\bra{I_{\mathbb{Z}}}=\aff_{\mathbb{Z}}\bra{I}$.

We now establish $\aff_{\mathbb{Z}}\bra{I}=\aff\bra{I_{\mathbb{Z}}}$ .
For the left to right containment note that  $I\cap \mathbb{Z}^d=I_{\mathbb{Z}}\cap\mathbb{Z}^d\subseteq I_{\mathbb{Z}}$ implies $\aff\bra{I\cap \mathbb{Z}^d}\subseteq \aff\bra{I_{\mathbb{Z}}}$, which shows the containment by the definition of $\aff_{\mathbb{Z}}\bra{I}$. For the reverse containment note that $I\cap \aff_{\mathbb{Z}}\bra{I}\subseteq \aff_{\mathbb{Z}}\bra{I}$ implies $\aff\bra{I\cap \aff_{\mathbb{Z}}\bra{I}}\subseteq \aff\bra{\aff_{\mathbb{Z}}\bra{I}}$, which shows the containment by the definition of $I_{\mathbb{Z}}$ and the fact that $\aff\bra{\aff_{\mathbb{Z}}\bra{I}}=\aff_{\mathbb{Z}}\bra{I}$.

For the final statement note that  $\mathcal{R}\bra{\aff\bra{J}}=\aff\bra{\mathcal{R}\bra{J}}$ for any $J\subseteq \Real^d$ because $\mathcal{R}$ is an affine function and $\mathcal{R}\bra{J\cap J'}=\mathcal{R}\bra{J}\cap\mathcal{R}\bra{J'}$ for any $J,J'\subseteq \Real^d$ because $\mathcal{R}$ is invertible.  Then
\[\mathcal{R}\bra{I_{\mathbb{Z}}}=\mathcal{R}\bra{I\cap \aff\bra{I\cap \mathbb{Z}^d}}= \mathcal{R}\bra{I}\cap \aff\bra{\mathcal{R}\bra{I}\cap \mathcal{R}\bra{\mathbb{Z}^d}} =\mathcal{R}\bra{I}_{\mathbb{Z}}\]
where the last equation follows from $\mathcal{R}\bra{\mathbb{Z}^d}=\mathbb{Z}^d$. Similarly, we also have
\[\mathcal{R}\bra{\aff_{\mathbb{Z}}\bra{I}}=\mathcal{R}\bra{\aff\bra{I\cap \mathbb{Z}^d}}=\aff\bra{\mathcal{R}\bra{I\cap \mathbb{Z}^d}}=\aff\bra{\mathcal{R}\bra{I}\cap \mathcal{R}\bra{\mathbb{Z}^d}}=\aff_{\mathbb{Z}}\bra{\mathcal{R}\bra{I}}\]
where the last equation again follows from $\mathcal{R}\bra{\mathbb{Z}^d}=\mathbb{Z}^d$.
\Halmos\endproof
We also use the following standard technique to align sub-spaces of dimension $k$ in $\Real^d$ with $ \mathbb{R}^{k}\times \set{\bm{0}}^{d-k}$ while preserving the integer lattice $\mathbb{Z}^d$ (\blue{proved in Appendix~\ref{Khintchine:sec} for completeness}).
\begin{definition}
A matrix $\bm{U}\in \mathbb{Z}^{d\times d}$ is \textbf{unimodular} if and only if it is invertible and $\bm{U}\mathbb{Z}^d=\bm{U}^{-1}\mathbb{Z}^d=\mathbb{Z}^d$.
\end{definition}
\begin{restatable}{lemma}{fulldimredlemma}\label{fulldimred} Let $F\subseteq \Real^d$ be a rational affine subspace such that $\dim\bra{F}=k$ and let $L\subseteq \Real^d$ be the $k$-dimensional rational linear subspace parallel to $F$ (i.e. $L=F-\bm{z}$ for any $\bm{z}\in F$). Then there exist $\bm{c}\in \mathbb{Z}^d\cap F$ and a unimodular matrix $\bm{U}\in \mathbb{Z}^{d\times d}$ such that
\begin{itemize}
    \item $\bm{U}L=\bm{U}\bra{F-\bm{c}}= \mathbb{R}^{k}\times \set{\bm{0}}^{d-k}$, and
    \item $(\bm{U}^T)^{-1} L^{\perp}=\set{\bm{0}}^{k}\times  \mathbb{R}^{d-k}$.
\end{itemize}
Furthermore, if $F=\bm{r}\mathbb{R}$ for $\bm{r}\in \mathbb{Z}^d\setminus \set{\bm{0}}$ with $\operatorname{gcd}(r_1,\ldots,r_d) = 1$, then we can choose $\bm{U}$ that also satisfies  $\bm{U}^{-1}\bm{e}(1)= \bm{r}$.
\end{restatable}
\pointtoproof{fulldimredlemma:proof}

Finally, the induction step in the proof of \autoref{interiorpointind} will use the following variant of Khintchine's Flatness Theorem. This variant follows easily from standard results from the geometry of numbers, but for completeness, we provide a proof in Appendix~\ref{Khintchine:sec}.

\begin{restatable}{lemma}{latticefreecorostatement}\label{latticefreecoro}
Let $I\subseteq \Real^d$ be a not necessarily full dimensional convex set such that $I\cap \mathbb{Z}^d\neq \emptyset$.
Then either $I_{\mathbb{Z}}$ has an integer point in its relative interior or there exist an invertible rational affine transformation $\mathcal{R}:\Real^d\to \Real^d$ such that $\mathcal{R}\bra{\mathbb{Z}^d}=\mathcal{R}^{-1}\bra{\mathbb{Z}^d}=\mathbb{Z}^d$ and $l,u\in \mathbb{Z}$ such that $\mathcal{R}\bra{I}_{\mathbb{Z}}\subseteq \set{\bz\in \mathbb{R}^{d}\,:\, l\leq z_1 \leq u}$ and $\aff_{\mathbb{Z}}\bra{\mathcal{R}\bra{I}}=\mathbb{R}^{\dim\bra{I_{\mathbb{Z}}}}\times \set{\bm{0}}^{d-\dim\bra{I_{\mathbb{Z}}}}$.
\end{restatable}
\pointtoproof{latticefreecorostatement:proof}
\subsubsection{Auxiliary Lemmas}

We also use the following lemma that allows us to rearrange the index set without loss of generality.

\begin{lemma}\label{unimodularequivformulation}
Let $M\subseteq \Real^{n+p+d}$ be a closed convex set \blue{such that the pair $(M,d)$ induces} a rational MICP formulation of $S \subseteq \mathbb{R}^n$ with index set  $I\subseteq \Real^d$  and $\bz$-projected sets $\set{A_{\bz}}_{\bz\in I}$. \blue{In addition, let} $\mathcal{R}:\Real^d\to\Real^d$ be an invertible rational affine transformation \blue{such that $\mathcal{R}^{-1}\bra{\mathbb{Z}^d}\subseteq \mathbb{Z}^d$} and
\[M'=\blue{{\mathcal{L}}\bra{M}=}\set{\bra{\bx,\by,\bz}\in \Real^{n+p+d}\,:\, \bra{\bx,\by,\mathcal{R}^{-1}\bra{\bz}}\in M},\]
\blue{where $\mathcal{L}:\Real^{n+p+d}\to \Real^{n+p+d}$ is such that  ${\mathcal{L}}\bra{\bx,\by,\bz}=\bra{\bx,\by,\mathcal{R}\bra{\bz}}$.}
\blue{Then $M'$ is a closed convex set and the pair $(M', d)$ induces a rational MICP formulation of
\begin{equation}\label{unimodularequivformulation:eqS}
    S'=\bigcup_{\bz'\in \mathcal{R}\bra{I}\cap \mathbb{Z}^d} A_{\mathcal{R}^{-1}\bra{\bz'}}
\end{equation}
with index set $I'=\mathcal{R}\bra{I}$ and $D_{\text{MICP}}\bra{M'}\leq D_{\text{MICP}}\bra{M}$.

If in addition, $\mathcal{R}\bra{\mathbb{Z}^d}=\mathcal{R}^{-1}\bra{\mathbb{Z}^d}=\mathbb{Z}^d$, then $S'=S$ and $D_{\text{MICP}}\bra{M}=D_{\text{MICP}}\bra{M'}$.
}
\end{lemma}
\proof{\textbf{Proof}}
\blue{Because $M$ is closed and convex and $\mathcal{L}$ is affine and invertible we have that $M'=\mathcal{L}\bra{M}$ is closed and convex. Hence, $(M',d)$ induces an MICP formulation. Let $I'$ be the index set of this formulation and  $\set{A'_{\bz'}}_{\bz'\in I'}$ be its $\bz$-projected sets. Because $I'=\mathcal{R}\bra{I}$ and  $\mathcal{R}$ is an invertible rational affine transformation we have that $I'$ is rationally unbounded and hence the formulation induced by $(M',d)$ is a rational MICP formulation. In addition, for any $\bz'\in I'=\mathcal{R}\bra{I}$ we have $\mathcal{R}^{-1}\bra{\bz'}\in I$ and $A'_{\bz'}=A_{\mathcal{R}^{-1}\bra{\bz'}}$. Then, $(M',d)$ induces an MICP formulation of the set $S'$ defined in \eqref{unimodularequivformulation:eqS} and  $D_{\text{MICP}}\bra{M'}\leq D_{\text{MICP}}\bra{M}$.

 For the last sentence of the lemma's statement, first note that invertibility of $\mathcal{R}$ and $\mathcal{R}\bra{\mathbb{Z}^d}=\mathbb{Z}^d$ imply \[\mathcal{R}\bra{I}\cap \mathbb{Z}^d=\mathcal{R}\bra{I}\cap \mathcal{R}\bra{\mathbb{Z}^d}=\mathcal{R}\bra{I\cap \mathbb{Z}^d}.\]
 In such case, the left hand side of
 \eqref{unimodularequivformulation:eqS} is equal to $\bigcup_{\bz'\in \mathcal{R}\bra{I\cap \mathbb{Z}^d}} A_{\mathcal{R}^{-1}\bra{\bz'}}=\bigcup_{\bz\in I\cap \mathbb{Z}^d} A_{\bz}=S$.}
  \Halmos
  \endproof

We also need the following auxiliary lemma. \blue{The goal of the lemma is to partition an index set $I\subseteq \Real^d$ into lower dimensional slices. Ideally, we would like these slices to be of the form $\set{\set{\bz\in I\,:\, z_1=i}}_{i\in [l,u]\cap \mathbb{Z}}$ for appropriately chosen $l,u$. Unfortunately, we require the slices remain rationally unbounded and this property may be lost when adding a restriction of the form $z_1=i$ (e.g. see \autoref{rationalaffinesections}).  }

\begin{lemma}\label{interiorpointindprop:auxlemma}
Let $l,u\in \mathbb{Z}$ and $I\subseteq \Real^d$ be a convex set such that
\begin{subequations}\label{interiorpointindprop:auxlemma:cond}
\begin{align}
\aff_{\mathbb{Z}}\bra{I}&=\mathbb{R}^{\dim\bra{\aff_{\mathbb{Z}}\bra{I}}}\times \set{\bm{0}}^{d-\dim\bra{\aff_{\mathbb{Z}}\bra{I}}}\label{interiorpointindonezeroother}\\
    I_{\mathbb{Z}}&\subseteq \set{\bz\in \Real^d\,:\,l \leq  z_1 \leq u}.\label{interiorpointindonezero}
\end{align}
\end{subequations}
Then,  there exists invertible rational affine transformations $\mathcal{R}_i:\Real^d\to \Real^d$ for each $i\in [l,u]\cap \mathbb{Z}$ such that
\begin{subequations}
\begin{alignat}{3}
\label{interiorpointindprop:claim1}    \dim\bra{\aff_{\mathbb{Z}}\bra{\mathcal{R}_i\bra{I}}}&\leq \dim\bra{\aff_{\mathbb{Z}}\bra{I}}-1,&\quad &\forall i\in [l,u]\cap \mathbb{Z}, \\
\label{interiorpointindprop:claim2} \bigcup_{i=l}^{u} \mathcal{R}_i\bra{I}\cap\mathbb{Z}^d&=   I\cap\mathbb{Z}^d, \text{ and}\\
\label{interiorpointindprop:claim3}\blue{\mathcal{R}^{-1}_i\bra{{\bz}}}&\blue{={\bz}}, &\quad &\blue{\forall  i\in [l,u]\cap \mathbb{Z},\quad  \bz\in \mathcal{R}_i\bra{I}\cap\mathbb{Z}^d}
\end{alignat}
\end{subequations}
\end{lemma}
\proof{\textbf{Proof}}
For each $i\in [l,u]\cap \mathbb{Z}$ let $\mathcal{R}_i:\Real^d\to \Real^d$ defined by
\[\mathcal{R}_i\bra{\bz}=\bra{\bm{I}-\bm{e}(1)\bm{e}(1)^T+\frac{\bm{e}(1)\bm{e}(1)^T}{2\max\set{u-i,i-l}}}\bra{\bz-i\bm{e}(1)}+i\bm{e}(1)\]
be the  invertible rational affine transformation that shrinks along the first coordinate so that
\begin{subequations}\label{interiorpointindonetwo}
\begin{alignat}{3}
\label{interiorpointindonetwo1}\mathcal{R}_i\bra{{\bz}}={\bz} \quad\quad&\Leftrightarrow\quad\quad z_1=i, \\
\mathcal{R}_i\bra{\set{\bz\in \Real^d\,:\, l\leq  z_1\leq u}}&\subseteq \set{\bz\in \Real^d\,:\, i -1 <  z_1< i +1 }.\label{interiorpointindonetwo2}
\end{alignat}
\end{subequations}
 We start by noting that \autoref{intafflemma} implies that
\begin{equation}\label{interiorpointindprop:auxlemmaIIZ}
    I_{\mathbb{Z}}\cap\mathbb{Z}^d=I\cap\mathbb{Z}^d,
\end{equation}
 and  \eqref{interiorpointindonezero} implies
\begin{equation}\label{interiorpointindone}
    I_{\mathbb{Z}}\cap \mathbb{Z}^{d} =\bigcup_{i=l}^{u} \set{\bz\in I_{\mathbb{Z}}\cap \mathbb{Z}^{d}\,:\, z_1=i}.
\end{equation}

We claim that for all $i\in [l,u]\cap \mathbb{Z}$, we have
\begin{equation}\label{finalslice} \set{\bz\in I_{\mathbb{Z}}\cap\mathbb{Z}^d\,:\, z_1=i}=\mathcal{R}_i\bra{I_{\mathbb{Z}}}\cap\mathbb{Z}^d=\mathcal{R}_i\bra{I}\cap\mathbb{Z}^d.\end{equation}
The first equation follows from  \eqref{interiorpointindonezero} and \eqref{interiorpointindonetwo}. For the second equation, note that $ \mathbb{Z}^d \subseteq \mathcal{R}_i\bra{\mathbb{Z}^d}$. Then
\begin{align*}
    \mathcal{R}_i\bra{I}\cap \mathbb{Z}^d &= \mathcal{R}_i\bra{I}\cap  \mathcal{R}_i\bra{\mathbb{Z}^d}\cap \mathbb{Z}^d\\
    &= \mathcal{R}_i\bra{I\cap\mathbb{Z}^d}\cap \mathbb{Z}^d\\
    &= \mathcal{R}_i\bra{I_{\mathbb{Z}}\cap\mathbb{Z}^d}\cap \mathbb{Z}^d\\
    &= \mathcal{R}_i\bra{I_{\mathbb{Z}}}\cap\mathcal{R}_i\bra{\mathbb{Z}^d}\cap \mathbb{Z}^d\\
    &= \mathcal{R}_i\bra{I_{\mathbb{Z}}}\cap \mathbb{Z}^d
\end{align*}
where the first and last equations follow from  $ \mathbb{Z}^d \subseteq \mathcal{R}_i\bra{\mathbb{Z}^d}$,  the second and fourth equations follow from the invertibility of $ \mathcal{R}_i$, and the third equation follows from \eqref{interiorpointindprop:auxlemmaIIZ}.

Then \eqref{interiorpointindprop:claim2} follows from \eqref{interiorpointindprop:auxlemmaIIZ}, \eqref{interiorpointindone} and \eqref{finalslice}, and  \eqref{interiorpointindprop:claim3} follows from \eqref{interiorpointindonetwo1} and \eqref{finalslice}.

To show \eqref{interiorpointindprop:claim1}, first note that \eqref{finalslice} implies
\[
\aff_{\mathbb{Z}}\bra{\mathcal{R}_i\bra{I}} =\aff\bra{\mathcal{R}_i\bra{I}\cap \mathbb{Z}^d}= \aff\bra{\set{\bz\in I_{\mathbb{Z}}\cap\mathbb{Z}^d\,:\, z_1=i}}\subseteq \set{\bz\in \aff\bra{I_{\mathbb{Z}}\cap\mathbb{Z}^d}\,:\, z_1=i}.
\]
Then \eqref{interiorpointindprop:claim1} holds because of \eqref{interiorpointindonezeroother}.\Halmos\endproof

\subsubsection{Proof of \autoref{interiorpointind}}

We now proceed with the proof of \autoref{interiorpointind}.
\interiorpointindprop*
\proof{\textbf{Proof}}
\label{interiorpointindprop:proof}
Let
$I\subseteq \Real^d$ be the index set of the MICP formulation induced by $M$ and $\set{A_{\bz}}_{\bz \in I}$
 \blue{its $\bz$-projected sets.} We prove the result by induction on $q=\dim\bra{I_{\mathbb{Z}}}\blue{=\dim\bra{\aff_{\mathbb{Z}}\bra{I}}}$.

For the base case $q=0$ we have that $\abs{I\cap \mathbb{Z}^d}=1$, \blue{and then $S$ is the projection of a closed convex set. Hence, $S$ is trivially binary MICP-R. Then the result follows with $S_0=S$ and $k=0$ (i.e. there are no sets $S_i$ with $i\neq 0$)}.

Now assume the result for $q-1$, and consider $M$ with index set $I$ such that $q=\dim\bra{I_{\mathbb{Z}}}$. If $\abs{I\cap \mathbb{Z}^d}<\infty$\blue{, then $S$ is a finite union of projections of closed convex sets, so  \autoref{lem:balasform} implies that $S$ is binary MICP-R, and} the result \blue{again} follows with $S_0=S$ and $k=0$. If $I_{\mathbb{Z}}$ has an integer point in its relative interior, the result follows with $S_0=\emptyset$, $k=1$, $M_1=M$ and $S_1=S$. Otherwise, \blue{by applying the rational affine transformation from \autoref{latticefreecoro} to $I$ (and $M$), we may assume that $I$ satisfies conditions \eqref{interiorpointindprop:auxlemma:cond} of \autoref{interiorpointindprop:auxlemma}. This assumption is without loss of generality because \autoref{unimodularequivformulation} ensures that the rational affine transformation from \autoref{latticefreecoro} preserves $D_{\text{MICP}}\bra{M}$ and that the pair $(M, d)$ continues to induce a rational MICP formulation of $S$. }

Let $\mathcal{R}_i:\Real^d\to \Real^d$ for each $i\in [l,u]\cap \mathbb{Z}$ be the invertible rational affine transformations from \autoref{interiorpointindprop:auxlemma}. For each $i\in [l,u]\cap \mathbb{Z}$, define the extended invertible rational affine transformation $\mathcal{L}_i:\Real^{n+p+d}\to \Real^{n+p+d}$ by ${\mathcal{L}}_i\bra{\bx,\by,\bz}=\bra{\bx,\by,\mathcal{R}_i\bra{\bz}}$, and let $\tilde{M}_i={\mathcal{L}}_i\bra{M}$. Then, \blue{by \autoref{unimodularequivformulation}}, the index set $\tilde{I}_i$ of $\tilde{M}_i$ is  $\tilde{I}_i=\mathcal{R}_i\bra{I}$, and \blue{the pair $(\tilde{M}_i,d)$ induces a rational MICP formulation of}
\begin{equation}\label{interiorpointindprop:zproj:tmp}
  \blue{  \tilde{S}_i=\bigcup_{\tilde{\bz}\in \tilde{I}_i\cap \mathbb{Z}^d} A_{\mathcal{R}^{-1}_i\bra{\tilde{\bz}}}}
\end{equation}
and
\begin{equation}\label{interiorpointindprop:c2:eq}
    D_{\text{MICP}}\bra{\tilde{M}_i}\leq D_{\text{MICP}}\bra{M}.
\end{equation}
\blue{In addition, \eqref{interiorpointindprop:claim3} in \autoref{interiorpointindprop:auxlemma} and \eqref{interiorpointindprop:zproj:tmp} imply }
\begin{equation}\label{interiorpointindprop:zproj}
 \blue{ \tilde{S}_i=\bigcup_{\tilde{\bz}\in \tilde{I}_i\cap \mathbb{Z}^d} A_{\tilde{\bz}}.}
\end{equation}
\blue{In addition, by \eqref{interiorpointindprop:claim1}--\eqref{interiorpointindprop:claim2} in} \autoref{interiorpointindprop:auxlemma}, we have
\begin{subequations}
\begin{align}
\label{interiorpointindpropcopy:claim1} \dim\bra{\aff_{\mathbb{Z}}\bra{\mathcal{R}_i\bra{I}}}=\dim\bra{\aff_{\mathbb{Z}}\bra{\tilde{I}^i}}&\leq \dim\bra{\aff_{\mathbb{Z}}\bra{I}}-1,\quad \forall i\in [l,u]\cap \mathbb{Z}, \quad\text{and}\\
\label{interiorpointindpropcopy:claim2} \bigcup_{i=l}^{u} \mathcal{R}_i\bra{I}\cap\mathbb{Z}^d=\bigcup_{i=l}^{u} \tilde{I}_i\cap\mathbb{Z}^d&=   I\cap\mathbb{Z}^d.
\end{align}
\end{subequations}
Then  \eqref{interiorpointindprop:zproj} and \eqref{interiorpointindpropcopy:claim2} imply
\begin{equation}\label{interiorpointindprop:almostlastunion}
    S=\bigcup_{i=l}^u \tilde{S}_i.
\end{equation}

Finally, \eqref{interiorpointindpropcopy:claim1} and the induction hypothesis imply that for each $i\in [l,u]\cap \mathbb{Z}$, there exists a binary MICP-R set $S_{i,0}\subseteq \Real^n$ and closed convex sets $\set{{M}_{i,j}}_{j=1}^{k_i}$ \blue{such that for each $j \in \sidx{k_i}$, the pair $(M_{i,j}, d)$ induces a rational MICP formulation of a non-empty set ${S}_{i,j}$, and together these sets decompose $\tilde{S}_i$ as the union} $\tilde{S}_i = S_{i,0} \cup \bigcup_{j=1}^{k_i} S_{i,j}$. For each $j\in \sidx{k_i}$,
we additionally have
\begin{subequations}\label{interiorpointindprop:c13:induction}
\begin{align}
         D_{\text{MICP}}\bra{{M}_{i,j}}&\leq D_{\text{MICP}}\bra{\tilde{M}_i},\\
     \bra{{I}_{i,j}}_{\mathbb{Z}} &\text{ has an integer point in its relative interior}
\end{align}
\end{subequations}
where ${I}_{i,j}$ is the index set of the formulation induced by ${M}_{i,j}$.

By \eqref{interiorpointindprop:almostlastunion} we have
\[
S=S_0 \cup \bigcup_{i=l}^u \bigcup_{j=1}^{k_i} {S}_{i,j}.
\]
where $S_0=\bigcup_{i=l}^u S_{i,0}$ is a binary MICP-R set by \autoref{basicproppropref}. The result then follows by noting that \eqref{interiorpointindprop:c13:induction} and \eqref{interiorpointindprop:c2:eq} imply \eqref{interiorpointindprop:c13}.
\Halmos\endproof

\subsection{Proof of \autoref{lem:generaldecompsimple}}

\subsubsection{Tools from convex analysis}

We first list several tools from convex analysis that we use to prove \autoref{lem:generaldecompsimple}.

We start with  the following affine separation result between a convex and a concave function  and two simple properties of convex functions.
\begin{lemma}\label{separationlemma}
Let $I\subseteq \Real^d$ be a convex set and $f,g:I \to \Real$ be such that $f$ and $-g$ are  convex functions and $g\bra{\bz}\leq f\bra{\bz}$ for all $\bz\in I$. Then there exist $\bm{v}\in \Real^d$ and $K\in \Real$ such that
\[ g\bra{\bz}\leq \bm{v}^T\bz +K \leq f\bra{\bz}\quad \forall \bz\in I.\]
\end{lemma}
\proof{\textbf{Proof}}
Follows from the first two pages of \cite[Section 31]{Rockafellar1997}.
\Halmos\endproof

\begin{lemma}\label{univariatelemma}
Let $f:[0,\infty) \to \Real$ be a convex function. Then
\begin{itemize}
    \item $f$ is continuous in $(0,\infty)$, and
    \item if $f$ is upper bounded, then $f$ is non-increasing.
    \end{itemize}
\end{lemma}
\proof{\textbf{Proof}}
Follows from Theorem 3.1.1 and Corollary 2.3.2 of~\cite{Hiriart-Urruty1996}.
\Halmos\endproof

A specific function we rely on is the support function of a set, which yields a precise description of it when the set is closed and convex (e.g. \cite[Section C.2]{hiriart-lemarechal-2001}).
\begin{definition} For $A\subseteq \Real^n$ we let its \textbf{support function} be $\sigma_A:\Real^n\to\Real\cup \set{+\infty}$ be defined by $\sigma_A\bra{\bm{c}}=\operatorname{sup}\{ {\bm{c}}^T\bm{x} : \bm{x} \in A \}$.
\end{definition}

\begin{lemma}\label{supportequalconvexset}
Let $A,B\subseteq \Real^n$ be convex sets. Then
\begin{itemize}
    \item $ \operatorname{cl}\bra{A}=\set{\bx\in \Real\,:\ \bm{c}^T\bx \leq \sigma_A\bra{\bm{c}}, \quad \forall \bm{c}\in \Real^n}$, and
\item  $\operatorname{cl}\bra{A}\subseteq \operatorname{cl}\bra{B}$ if and only if $ \sigma_A\bra{\bm{c}}\leq  \sigma_B\bra{\bm{c}}$ for all $\bm{c}\in \Real^n$.
\end{itemize}
\end{lemma}
\proof{\textbf{Proof}}
Follows from Theorems C.2.2.2 and C.3.3.1, of \cite{hiriart-lemarechal-2001}.
\Halmos\endproof

\subsubsection{Auxiliary Lemmas}

We also use the following lemma that describes properties of the  support function of the  $\bz$-projected sets viewed as a function of $\bz$ (for fixed $\bm{c}$).
\begin{lemma}\label{lem:concavesupport}
Let $M\subseteq \Real^{n+p+d}$ be a closed convex set inducing an MICP formulation with index set $I$ and  $\bz$-projected sets   $\set{A_{\bz}}_{\bz\in I}$. Furthermore, for each $\bm{c}\in \mathbb{R}^n$ let the three functions $g_{{\bm{c}}}, f_{{\bm{c}}}, w_{{\bm{c}}} : I \to \mathbb{R} \cup \{-\infty,+\infty\} $ be given by
\begin{itemize}
    \item $f_{{\bm{c}}}(\bm{z})=-\sigma_{A_{\bm{z}}}\bra{-\bm{c}} =  \operatorname{inf}\{ {\bm{c}}^T\bm{x} : \bm{x} \in A_{\bm{z}} \}$,
    \item $g_{{\bm{c}}}(\bm{z})=\sigma_{A_{\bm{z}}}\bra{\bm{c}} = \operatorname{sup}\{ {\bm{c}}^T\bm{x} : \bm{x} \in A_{\bm{z}} \}$, and
    \item $w_{{\bm{c}}}(\bm{z})=g_{{\bm{c}}}(\bm{z})-f_{{\bm{c}}}(\bm{z}) $.
\end{itemize}
Then $f_{{\bm{c}}}(\bm{z})$ is a convex function, and $g_{{\bm{c}}}(\bm{z})$ and $w_{{\bm{c}}}(\bm{z})$ are concave. Furthermore, for all $\bz\in I$, we have $f_{{\bm{c}}}(\bm{z})<\infty$,  $g_{{\bm{c}}}(\bm{z})>-\infty$, and $w_{{\bm{c}}}(\bm{z})\blue{\ge 0}$.

If, in addition, $I_{\mathbb{Z}}$ has an integer point in its relative interior, then for  any ${\bm{c}}\in \mathbb{R}^n$, we have
\begin{equation}\label{boundedzprojected}
   \sup_{\bm{z} \in I_{\mathbb{Z}}\cap \mathbb{Z}^d} w_{\bm{c}}(\bm{z}) <\infty \quad\quad \Leftrightarrow \quad\quad \sup_{\bm{z} \in I_{\mathbb{Z}}} w_{\bm{c}}(\bm{z})  <\infty.
\end{equation}
\end{lemma}
\proof{\textbf{Proof}}
Note $A_{\bm{z}}$ is nonempty for $\bm{z} \in I$, so $g_{{\bm{c}}}(\bm{z})$ is well defined. Choose any two $\bm{z}^1,\bm{z}^2 \in I$ and $\lambda \in [0,1]$. We will show $g_{{\bm{c}}}(\lambda \bm{z}^1 + (1-\lambda)\bm{z}^2) \ge \lambda g_{{\bm{c}}}(\bm{z}^1) + (1-\lambda) g_{{\bm{c}}}(\bm{z}^2)$. Let
$\set{\bm{x}^m}_{m\in \mathbb{N}}$ and $\set{\bm{y}^m}_{m\in \mathbb{N}}$ be sequences contained in $A_{\bm{z}^1}$ and $A_{\bm{z}^2}$ respectively such that $\lim_{m \to \infty} {\bm{c}}^T\bm{x}^m = g_{{\bm{c}}}(\bm{z}^1)$ and $\lim_{m \to \infty} {\bm{c}}^T\bm{y}^m = g_{{\bm{c}}}(\bm{z}^2)$. Since $M$ is convex, it follows that for each $i$, $\lambda \bm{x}^m + (1-\lambda) \blue{\bm{y}}^m \in A_{\lambda \bm{z}^1 + (1-\lambda) \bm{z}^2}$, so
\begin{equation*}
g_{{\bm{c}}}(\lambda \bm{z}^1 + (1-\lambda)\bm{z}^2) \ge \lim_{m \to \infty} {\bm{c}}^T(\lambda \bm{x}^m + (1-\lambda) \bm{y}^m) = \lambda g_{{\bm{c}}}(\bm{z}^1) + (1-\lambda) g_{{\bm{c}}}(\bm{z}^2).
\end{equation*}
The result for $f$ is analogous, and the result for $w$ follows because it is a sum of two concave functions. The fact that $A_{\bz}\neq\emptyset$ for all $\bz\in I$ implies $f_{{\bm{c}}}(\bm{z})<\infty$,  $g_{{\bm{c}}}(\bm{z})>-\infty$ and $w_{{\bm{c}}}(\bm{z})\blue{\ge 0}$.

For the final statement, first note that because $L=\aff_{\mathbb{Z}}\bra{I}=\aff\bra{I_{\mathbb{Z}}}$ is a rational affine subspace, we can use \autoref{fulldimred} to reduce to the case in which $I_{\mathbb{Z}}$ is full-dimensional and has an integer point in its interior. Then let $B=\set{\bra{\bz,t}\in \Real^{d+1}\,:\, \bz\in I_{\mathbb{Z}},\quad - w_{\bm{c}}(\bm{z}) \leq t}$ so that $B$ is a convex set such that $\Int\bra{B}\cap \bra{\mathbb{Z}^d\times \Real}\neq \emptyset$. Then, by Proposition 4.5 in \cite{Strong-Dual-for-Conic}, we have that $\inf\set{t\,:\, \bra{\bz,t}\in B}=-\infty$ if and only if $\inf\set{t\,:\, \bra{\bz,t}\in B\cap \bra{\mathbb{Z}^d\times \Real}}=-\infty$, which shows the result.
\Halmos\endproof

The fact that $I_{\mathbb{Z}}$ containing an integer point in its relative interior is needed for \eqref{boundedzprojected} is illustrated in the following example.
\begin{example}\label{interiorIZ:example}
  Let  $M=\set{(x,y,z)\in \Real^3 \,:\,x^2\leq y z,\quad (x-1)^2\leq y (1-z),\quad 0\leq z\leq 1}$ be a closed convex set\footnote{See \autoref{Jeroslowexample} for the closure and convexity of this set.}   inducing a MICP formulation of $S=\set{0,1}$. Then $I = \proj_z(M)=[0,1]$,
  \[f_{1}(z)=\begin{cases}0 & z=0\\ 1 & z=1\\ -\infty & z\in (0,1)\end{cases},\quad\text{ and }\quad  g_{1}(z)=\begin{cases}0 & z=0\\ 1 & z=1\\ +\infty & z\in (0,1)\end{cases}.\Halmos\]
\end{example}

To prove  \autoref{lem:generaldecompsimple}, we first need to understand the relaxation of the integrality constraints of a formulation along one of its rational unbounded directions. We then need to connect rational unbounded directions of the index set to periodic directions of $S$. We respectively achieve these tasks through the following two lemmas.

\begin{lemma}\label{relaxtransformlemma}
Let  $M\subseteq \Real^{n+p+d}$ be a closed convex set \blue{such that the pair $(M, d)$ induces} a rational MICP formulation with index set $I\subseteq \Real^d$ and $\bz$-projected sets $\set{A_{\bz}}_{\bz\in I}$.  If there exists $\bm{r} \in \bra{\mathbb{Z}^d\cap I_\infty}\setminus\set{\bm{0}}$, then there exists a closed convex set $\tilde{M}\subseteq \Real^{n+(p+1)+(d-1)}$ \blue{such that the pair $(\tilde{M},d-1)$ induces} a rational MICP formulation \blue{of some set with the property that} $\tilde{A}$ is a $\bz$-projected set of the formulation induced by $(\tilde{M}, d-1)$ if and only if there exist
$\bz^0\in \mathbb{Z}^d\cap I$ such that
\[\tilde{A}=\bigcup\nolimits_{\bz'\in \bra{\bz^0+ \bm{r} \Real}\cap I} A_{\bz'}.\]
In particular, $\bigcup\nolimits_{\bz'\in \bra{\bz^0+ \bm{r} \Real}\cap I} A_{\bz'}$ is a convex set for any $\bz^0\in \mathbb{Z}^d\cap I$.
\end{lemma}
\proof{\textbf{Proof}}
We may assume without loss of generality that $\operatorname{gcd}(r_1,\ldots,r_d) = 1$. Then, by \autoref{fulldimred}, there exists an invertible $\bm{U}\in \mathbb{Z}^{d\times d}$ such that $\bm{U} \mathbb{Z}^d=\bm{U}^{-1} \mathbb{Z}^d= \mathbb{Z}^d$ and $\bm{U}^{-1}\bm{e}(1)= \bm{r}$. Let
\[\tilde{M}=\set{\bra{\bx,\tilde{\by},\tilde{\bz}}\in \Real^{n+(p+1)+(d-1)}\,:\, \bra{\bx,\;\tilde{\by}_{1:p},\;\bm{U}^{-1}\bra{\tilde{y}_{p+1},\tilde{\bz}}}\in M}\]
where $\tilde{\by}_{1:p}\in \Real^n$ is the vector composed of the first $p$ components of $\tilde{\by}\in \Real^{p+1}$.
We claim \blue{$(\tilde{M}, d-1)$} induces the desired MICP formulation. Now let $\tilde{I}=\proj_{\tilde{\bz}}\bra{\tilde{M}}$ be the index set of $\tilde{M}$,
\[\bm{V}_M=\begin{pmatrix}\bm{I}_{(n+p)\times (n+p)} &\bm{0}_{(n+p)\times d}&\\ \bm{0}_{d\times (n+p)} & \bm{U}\end{pmatrix}\quad\text{and}\quad \bm{V}_I=\begin{pmatrix}\bm{0}_{(n+p)\times (n+p)} &\bm{0}_{(n+p)\times d}&\\ \bm{0}_{d\times (n+p)} & \bm{I}_{d\times d}\end{pmatrix}\]
so that $\tilde{M}=\bm{V}_MM$ and $I=\bm{V}_I M$.
Then
$\tilde{I}=\proj_{\tilde{\bz}}\bra{\bm{V}_M M}=\proj_{\tilde{\bz}}\bra{\bm{V}_M\bm{V_I}M}=\proj_{\tilde{\bz}}\bra{\bm{U}I}$, and hence $\tilde{I}$ is rationally unbounded as it is a rational affine transformation of $I$, which is rationally unbounded by assumption. Finally, let $\Lambda\bra{\tilde{\bz}}=\set{\lambda\in \Real\,:\, \bra{\lambda,\tilde{\bz}} \in \bm{U}I}$. Then, for any $\tilde{\bz}\in \tilde{I}$, we have
\begin{align*}
    \tilde{A}_{\tilde{\bz}}&=\proj_{\bx}\bra{\tilde{M} \cap \bra{\mathbb{R}^{n+(p+1)} \times \{ \tilde{\bm{z}} \}}}\\
    &=\bigcup_{\lambda\in \Lambda\bra{\tilde{\bz}}}\proj_{\bx}\bra{\tilde{M} \cap \bra{\mathbb{R}^{n+p} \times \set{ \bra{\lambda,\tilde{\bm{z}}} }}}\\
    &=\bigcup_{\lambda\in \Lambda\bra{\tilde{\bz}}}\proj_{\bx}\bra{\bm{V}_M\bra{M \cap \bra{\mathbb{R}^{n+p} \times \set{\bm{U}^{-1} \bra{\lambda,\tilde{\bm{z}}} }}  }}\\
        &=\bigcup_{\lambda\in \Lambda\bra{\tilde{\bz}}}\proj_{\bx}\bra{M \cap \bra{\mathbb{R}^{n+p} \times \set{\bm{U}^{-1} \bra{\lambda,\tilde{\bm{z}}} }}  }\\
             &=\bigcup_{\lambda\in \Lambda\bra{\tilde{\bz}}}A_{\bm{U}^{-1}\bra{\lambda,\tilde{\bm{z}}}}= \bigcup_{\lambda\in \Lambda\bra{\tilde{\bz}}}A_{\lambda \bm{r}+\bm{U}^{-1}\bra{0,\tilde{\bm{z}}}} =\bigcup\nolimits_{\bz'\in \bra{\bm{U}^{-1}\bra{0,\tilde{\bm{z}}}+ \bm{r} \Real}\cap I} A_{\bz'}.
    \end{align*}
To finish the result, note that $\bm{U}^{-1}\bm{e}(1)= \bm{r}$ and $\tilde{\bz}\in \tilde{I}$ if an only if there exists $\bz^0\in I$ such that $\bra{0,\tilde{\bz}}=\bra{\bm{I}-\bm{e}(1)\bm{e}(1)^T}\bm{U}\bz^0$. Then $\bm{U}^{-1}\bra{0,\tilde{\bm{z}}}+ \bm{r} \Real = \bz^0-\bm{r}\bra{\bm{e}(1)^T\bm{U}\bz}+ \bm{r} \Real=\bz^0+ \bm{r} \Real$ \Halmos\endproof

\begin{lemma}\label{newnewqraylemmasimple}
Let  $M\subseteq \Real^{n+p+d}$ be a closed convex set inducing a rational MICP formulation with index set $I\subseteq \Real^d$ and $\bz$-projected sets $\set{A_{\bz}}_{\bz\in I}$. If $D_{\text{MICP}}\bra{M}<\infty$ and $I_{\mathbb{Z}}$ has an integer point in its relative interior, then for any $\bm{r}\in \mathbb{Z}^d \cap (I_\infty \setminus \set{\bm{0}})$ there exists a unique $\bm{u}\in \Real^n $ such that
\begin{equation}\label{newnewqraylemmaperiodicitysimple}
   \bx+\lambda \bm{u} \in \operatorname{cl}\bra{A_{\bz+\lambda \bm{r}}} \quad \forall \; \lambda >0, \;\bm{z} \in\mathbb{Z}^d\cap I, \;\bx\in A_{{\bz}}.
\end{equation}
Furthermore, if  $\bm{u}=\bm{0}$, then for all $\bm{z} \in \mathbb{Z}^d\cap I$ we have that  $\tilde{A}_{\bz}=\bigcup\nolimits_{\bz'\in \bz+ \bm{r} \Real_+} A_{\bz'}$ is a convex set such that $D\bra{\tilde{A}_{\bz}}\leq D_{\text{MICP}}\bra{M}$ and
\[ \bigcup\nolimits_{\bz'\in \bra{\bz+ \bm{r} \Real_+}\cap \mathbb{Z}^d}A_{\bz'}\subseteq \tilde{A}_{\bz}\subseteq \bigcup\nolimits_{\bz'\in \bra{\bz+ \bm{r} \Real_+}\cap \mathbb{Z}^d} \cl\bra{A_{\bz'}}.\]
  \end{lemma}
\proof{\textbf{Proof}}
For $\bm{c}\in \mathbb{R}^n$ let $f_{{\bm{c}}}(\bm{z})=-\sigma_{A_{\bm{z}}}\bra{-\bm{c}}$, $g_{{\bm{c}}}(\bm{z})=\sigma_{A_{\bm{z}}}\bra{\bm{c}}$, and $w_{{\bm{c}}}(\bm{z})=g_{{\bm{c}}}(\bm{z})-f_{{\bm{c}}}(\bm{z})$ be the functions from \autoref{lem:concavesupport}. By \autoref{lem:concavesupport} and assumption $D_{\text{MICP}}\bra{M}<\infty$
we have that $f_{\bm{c}}$ and $-g_{\bm{c}}$ are finite-valued convex functions on $ I_{\mathbb{Z}}$, and for each ${\bm{c}}\in \mathbb{R}^n$, there exist $K_{\bm{c}}\geq 0$ such that $w_{\bm{c}}(\bm{z}) \leq K_{\bm{c}}$ for all $\bm{z} \in I_{\mathbb{Z}}$. We then have that
\begin{equation}\label{newnewqraylemmaperiodicitysimpleineq}
f_{\bm{c}}(\bm{z})\leq g_{\bm{c}}(\bm{z}) \leq f_{\bm{c}}(\bm{z}) + K_{\bm{c}}\quad \forall \bz\in  I_{\mathbb{Z}},
\end{equation}
and then by  \autoref{separationlemma} (for $g(\bm{z})=g_{\bm{c}}(\bm{z})$ and $f(\bm{z})=f_{\bm{c}}(\bm{z}) + K_{\bm{c}}$), there exist $\bm{v}({\bm{c}})\in \Real^n$ and $\underline{K}_{\bm{c}},\,\overline{K}_{\bm{c}}\in \Real$ such that
\begin{equation}\label{decompboundnewsimple}\bm{v}({\bm{c}})^T\bm{z}+\underline{K}_{\bm{c}}\leq f_{\bm{c}}(\bm{z})\leq g_{\bm{c}}(\bm{z}) \leq \bm{v}({\bm{c}})^T\bm{z}+\overline{K}_{\bm{c}} \quad \forall \bz\in  I_{\mathbb{Z}}.\end{equation}

By \autoref{supportequalconvexset},  for any $\lambda\geq0$, $\tilde{\bm{z}} \in I \cap \mathbb{Z}^d$, $\bx\in A_{{\tilde{\bm{z}}}}$ and $\bm{u}\in \Real^n$ we have that  $\bx+\lambda \bm{u} \in \operatorname{cl}\bra{A_{\tilde{\bz}+\lambda \bm{r}}}$ is equivalent to  $ \bra{\bx+\lambda \bm{u},\lambda }\in \tilde{M}_{\tilde{\bz}}$ for
\begin{alignat*}{4}
\tilde{M}_{\tilde{\bz}}:&=\bigl\{ \bra{\bx,\lambda}\in  \Real^{n+1}\,:\, \lambda\geq 0,&\quad f_{\bm{c}}\bra{\tilde{\bz}+\lambda \bm{r}}&\leq \bm{c}^T \bx &&\leq  g_{\bm{c}}\bra{\tilde{\bz}+\lambda \bm{r}}\quad &\forall \bm{c}\in \Real^n\bigr\}\\
&=\bigl\{ \bra{\bx,\lambda}\in \Real^{n+1}\,:\, \lambda\geq 0,&\quad  f_{\bm{c},\tilde{\bz}}^1\bra{\lambda}&\leq \bm{c}^T \bx &&\leq  g_{\bm{c},\tilde{\bz}}^1\bra{\lambda}\quad &\forall \bm{c}\in \Real^n\bigr\}
\end{alignat*}
 with $f_{\bm{c},\tilde{\bz}}^1\bra{\lambda}=f_{\bm{c}}\bra{\tilde{\bz}+\lambda \bm{r}}$ and $g_{\bm{c},\tilde{\bz}}^1\bra{\lambda}=g_{\bm{c}}\bra{\tilde{\bz}+\lambda \bm{r}}$. By \autoref{univariatelemma} (and the fact that $f_{\bm{c},\tilde{\bz}}^1\bra{\lambda}$ and $-g_{\bm{c},\tilde{\bz}}^1\bra{\lambda}$ are convex), these two functions are continuous for $\lambda>0$. Then, for all  $\lambda>0$, we have  $ \bra{\hat{\bx},\lambda }\in \tilde{M}_{\tilde{\bz}}$ is equivalent to $ \bra{\hat{\bx},\lambda }\in \cl\bra{\tilde{M}_{\tilde{\bz}}}$, which yields
 \begin{equation}\label{newnewqraylemmasimpleEE}
     \forall \lambda >0,\; \tilde{\bm{z}} \in I \cap \mathbb{Z}^d,\;
     \quad \quad
     \hat{\bx}\in \operatorname{cl}\bra{A_{\tilde{\bz}+\lambda \bm{r}}}\Leftrightarrow \bra{
     \hat{\bx},\lambda }\in \cl\bra{\tilde{M}_{\tilde{\bz}}}.
 \end{equation}
 Using \eqref{decompboundnewsimple}, we have
$\tilde{M}_{\tilde{\bz}}\subseteq V$
for
\[V=\set{ \bra{\bx,\lambda}\in  \Real^{n}\times \Real_+\,:\,\lambda\geq 0,\quad \bm{v}({\bm{c}})^T\bra{\tilde{\bz}+\lambda \bm{r}}+\underline{K}_{\bm{c}}\leq \bm{c}^T \bx \leq\bm{v}({\bm{c}})^T\bra{\tilde{\bz}+\lambda \bm{r}}+\overline{K}_{\bm{c}}\quad \forall \bm{c}\in \Real^n},\]
and because $V$ is closed, we also have $\cl\bra{\tilde{M}_{\tilde{\bz}}}\subseteq V$. Then, by \autoref{rationaldefcprop}, we have
\begin{align*}\bra{\cl\bra{\tilde{M}_{\tilde{\bz}}}}_\infty\subseteq V_\infty =&\bigl\{ \bra{\bx,\lambda}\in  \Real^{n}\times \Real_+\,:\,  \lambda\geq 0,&\quad \bm{v}({\bm{c}})^T\bm{r} \lambda&=  \bm{c}^T \bx   \quad\quad\quad &&\forall   \bm{c}\in \Real^n\bigr\}\\
\subseteq &\bigl\{ \bra{\bx,\lambda}\in  \Real^{n}\times \Real_+\,:\,  \lambda\geq 0,&\quad \bm{v}({\bm{e}(i)})^T\bm{r} \lambda&= \bm{e}(i)^T \bx=x_i   \quad &&\forall   i\in \sidx{n}\bigr\}.\\
=&\set{\bra{\lambda\bm{u},\lambda}\,:\, \lambda \geq 0}
\end{align*}
for $\bm{u}=\sum_{i=1}^n \bm{e}(i) \bm{r}^T\bm{v}({\bm{e}(i)})$. \blue{In particular, $\bra{\cl\bra{\tilde{M}_{\tilde{\bz}}}}_\infty$ has dimension at most 1, and hence} we either have $\bra{\cl\bra{\tilde{M}_{\tilde{\bz}}}}_\infty=\set{\bra{\lambda\bm{u},\lambda}\,:\, \lambda \geq 0}$ or  $\bra{\cl\bra{\tilde{M}_{\tilde{\bz}}}}_\infty= \set{\bm{0}}$. However, because $\bm{r}\in \mathbb{Z}^d \cap (I_\infty \setminus \set{\bm{0}})$, then for any $\lambda > 0$ there exists $\bx\bra{\lambda} \in  A_{\tilde{\bz}+\lambda \bm{r}}\subseteq \operatorname{cl}\bra{A_{\tilde{\bz}+\lambda \bm{r}}}$. Together with \eqref{newnewqraylemmasimpleEE}, $\bx\bra{\lambda} \in \operatorname{cl}\bra{A_{\tilde{\bz}+\lambda \bm{r}}}$ for all $\lambda > 0$ implies   $\bra{\bx\bra{\lambda},\lambda} \in \cl\bra{\tilde{M}_{\tilde{\bz}}}$ for all $\lambda > 0$. Then, by taking $\lambda\to \infty$, we conclude  that $\cl\bra{\tilde{M}_{\tilde{\bz}}}$ is unbounded and hence  $\bra{\cl\bra{\tilde{M}_{\tilde{\bz}}}}_\infty\neq \set{\bm{0}}$ (e.g. \cite[Proposition A.2.2.3]{hiriart-lemarechal-2001}). Hence $\bra{\cl\bra{\tilde{M}_{\tilde{\bz}}}}_\infty=\set{\bra{\lambda\bm{u},\lambda}\,:\, \lambda \geq 0}$. For any $\tilde{\bx}\in A_{\tilde{\bz}}$ we have $ \bra{\tilde{\bx},0 }\in \tilde{M}_{\tilde{\bz}}\subseteq \cl\bra{\tilde{M}_{\tilde{\bz}}}$, which implies that $ \bra{\tilde{\bx}+\lambda\bm{u},\lambda }\in\cl\bra{\tilde{M}_{\tilde{\bz}}}$ for all $\lambda> 0$. Then \eqref{newnewqraylemmasimpleEE} for $\hat{\bx}=\tilde{\bx}+\lambda\bm{u}$ implies that $\bm{u}$ satisfies  \eqref{newnewqraylemmaperiodicitysimple}.

Furthermore, because $\cl\bra{\tilde{M}_{\tilde{\bz}}}\subseteq V$, we also have
\begin{equation}\label{newnewqraylemmasimpleEEtwo}
\bm{v}({\bm{c}})^T\bra{\tilde{\bz}+\lambda \bm{r}}+\underline{K}_{\bm{c}}\leq \bm{c}^T \bra{\tilde{\bx}+\lambda\bm{u}} \leq\bm{v}({\bm{c}})^T\bra{\tilde{\bz}+\lambda \bm{r}}+\overline{K}_{\bm{c}}\quad \forall \bm{c}\in \Real^n,\; \lambda >0
\end{equation}
so if $\bm{u}=\bm{0}$, then  $\bm{v}({\bm{c}})^T\bm{r}=0$ for all $\bm{c}\in \Real^n$ (otherwise, we obtain a contradiction with \eqref{newnewqraylemmasimpleEEtwo} by letting $\lambda\to \infty$). Together with \eqref{decompboundnewsimple} this implies that both $f_{\bm{c},\tilde{\bz}}^1\bra{\lambda}=f_{\bm{c}}\bra{\tilde{\bz}+\lambda \bm{r}}$ and $g_{\bm{c},\tilde{\bz}}^1\bra{\lambda}=g_{\bm{c}}\bra{\tilde{\bz}+\lambda \bm{r}}$ are  upper and lower bounded. Then, by convexity of $f_{\bm{c},\tilde{\bz}}^1\bra{\lambda}$ and $-g_{\bm{c},\tilde{\bz}}^1\bra{\lambda}$ and by \autoref{univariatelemma}, both $f_{\bm{c},\tilde{\bz}}^1\bra{\lambda}$ and $-g_{\bm{c},\tilde{\bz}}^1\bra{\lambda}$ are non-increasing. In particular, by the definition of $g_{\bm{c},\tilde{\bz}}^1$ we have that $ \sigma_{A_{\tilde{\bz}+\lambda \bm{r}}}\bra{\bm{c}}$ is non-decreasing as a function of $\lambda\geq0$  for all $\bm{c}\in\Real^n$. Then, by \autoref{supportequalconvexset}, we have
\begin{equation}\label{newnewqraylemmasimple:blue}
\cl\bra{A_{\tilde{\bz}+\lambda \bm{r}}}\subseteq \cl\bra{A_{\tilde{\bz}+\lambda' \bm{r}}} \quad \forall 0\leq \lambda \leq \lambda',
\end{equation}
which implies $\bigcup\nolimits_{\bz'\in \bz+ \bm{r} \Real_+} \cl\bra{A_{\bz'}}=\bigcup\nolimits_{\bz'\in \bra{\bz+ \bm{r} \Real_+}\cap \mathbb{Z}^d} \cl\bra{A_{\bz'}}$ for all $\bm{z} \in \mathbb{Z}^d\cap I$. Then,
\begin{equation}\label{newnewqraylemmasimple:bluetoo}
    \tilde{A}_{\bz}\subseteq \bigcup\nolimits_{\bz'\in \bra{\bz+ \bm{r} \Real_+}\cap \mathbb{Z}^d} \cl\bra{A_{\bz'}} \quad \forall \bm{z} \in \mathbb{Z}^d\cap I
\end{equation}
for  $\tilde{A}_{\bz}=\bigcup\nolimits_{\bz'\in \bz+ \bm{r} \Real_+} A_{\bz'}$ as required.  By \autoref{relaxtransformlemma}, we have that  $\tilde{A}_{\bz}$ is convex. \blue{Using the equivalent characterization of diameter using support functions as} $D(C)=\sup\set{\sigma_C\bra{\bm{c}} + \sigma_C\bra{-\bm{c}} \,:\, \bm{c}\in \Real^n,\quad \norm{\bm{c}}_2=1}$, we have \blue{also} that \eqref{newnewqraylemmasimple:blue} and \autoref{supportequalconvexset} imply
\[
D\bra{\tilde{A}_{\bz}}\leq \sup_{\bz'\in \bra{\bz+ \bm{r} \Real_+}\cap \mathbb{Z}^d} D\bra{A_{\bz'}}\leq D_{\text{MICP}}\bra{M}.
\]
\Halmos\endproof

\subsubsection{Proof of \autoref{lem:generaldecompsimple}}
\generaldecompsimpleprop*
\proof{\textbf{Proof}}
\label{generaldecompsimpleprop:proof}
Let $I\subseteq \Real^d$ be the index set of the MICP formulation induced by $M$ and  $\set{A_{\bz}}_{\bz\in I}$ be the $\bz$-projected sets of this formulation.  If  $I\cap \mathbb{Z}^d$ is finite, then $S$ is a finite union of projections of closed convex sets, so \autoref{lem:balasform} implies option \ref{lem:generaldecompc1simple}. If $I\cap \mathbb{Z}^d$ is infinite, by the rational unboundedness of $I$  there exists $\bm{r} \in \bra{\mathbb{Z}^d\cap I_\infty}\setminus\set{\bm{0}}$, for which we may assume without loss of generality that $\operatorname{gcd}\bra{r_1,\ldots,r_d}=1$.

Because $I_{\mathbb{Z}}$ has an integer point in its relative interior, then by \autoref{newnewqraylemmasimple} there exists $\bm{u}\in \Real^n$ such that  $\bx+\lambda \bm{u} \in \operatorname{cl}\bra{A_{\bz+\lambda \bm{r}}}$ for all $\lambda >0$, $\bm{z} \in I\cap \mathbb{Z}^d$ and $\bx\in A_{{\bz}}$.  If $\bm{u}\neq \bm{0}$, we get option \ref{lem:generaldecompc2simple} by noting that $\bz+\lambda \bm{r}\in \mathbb{Z}^d$ if and only if $\lambda\in \mathbb{Z}$ (because $\operatorname{gcd}\bra{r_1,\ldots,r_d}=1$). If $\bm{u}= \bm{0}$, by \autoref{relaxtransformlemma} and \autoref{newnewqraylemmasimple} there exists a closed convex set $\tilde{M}\subseteq \Real^{n+(p+1)+(d-1)}$ \blue{such that the pair $(\tilde{M}, d-1)$ induces a} rational MICP formulation of a set $\tilde{S}\subseteq \Real^n$ such that
\begin{itemize}
    \item $D_{\text{MICP}}\bra{\tilde{M}}\leq D_{\text{MICP}}\bra{{M}}$,
    \item if $\tilde{I}$ and $\set{\tilde{A}_{\bz}}_{\bz\in \tilde{I}}$ are the index set and $\bz$-projeted sets of the formulation induced by $\tilde{M}$, then
    \begin{equation}\label{generaldecompsimpleprop:eq1}
    S=\bigcup_{\bz\in I\cap \mathbb{Z}^d}\bigcup_{\bz'\in \bra{\bz+ \bm{r} \Real_+}\cap  \mathbb{Z}^d}A_{\bz'}\subseteq \tilde{S}=\bigcup_{\tilde{\bz}\in \tilde{I}\cap \mathbb{Z}^{d-1}}\tilde{A}_{\tilde{\bz}}\subseteq \bigcup_{\bz\in I\cap \mathbb{Z}^d} \bigcup_{\bz'\in \bra{\bz+ \bm{r} \Real_+}\cap \mathbb{Z}^d} \cl\bra{A_{\bz'}}\subseteq \cl\bra{S}.
    \end{equation}
\end{itemize}

We then get option \ref{lem:generaldecompc3simple} by noting that $\tilde{S}\subseteq \cl\bra{S}$ implies $\cl\bra{\tilde{S}}\subseteq \cl\bra{S}$.
\Halmos\endproof

\section{Conclusion}\label{conclusion:sec}
We would like to highlight two important modeling insights evident from the characterizations of MICP representability obtained in this work. First, unlike for the simpler case of MILP-R sets, the class of MICP-R sets is closed under unions because of features of convex sets that do not appear when considering only polyhedra. In particular, any finite union of polyhedra is MICP-R but not necessarily MILP-R. Second, under reasonable assumptions on the MICP formulation, unbounded integer variables induce periodicity and are therefore necessary only for modeling unbounded sets, a statement which was known to be true only for MILP formulations.

We conclude by posing three open questions whose answers we believe could yield more elegant characterizations of MICP-R sets:
\begin{enumerate}
    \item The ``midpoint lemma'' (\autoref{midpointlemmaref}) developed a simple necessary condition for MICP representability. Is it also a sufficient condition?
    \item Is the assumption of rationality in the statement of \autoref{thm:compactfinite} necessary? That is, are all compact MICP-R sets finite unions of compact convex sets?
    \item Can the geometric assumption on the diameter of convex subsets in \autoref{periodictheoref} be relaxed?
\end{enumerate}

\section*{Acknowledgements}
The second author was supported by the National Science Foundation  under Grant No. CMMI-1351619 and by the Office of Naval Research under Grant No. N00014-18-1-2079

\bibliographystyle{ormsv080}
\bibliography{main}

\begin{APPENDICES}

\section{Rational conic representation of a non-rational polyhedron}\label{app:conicratpolyirrat}
\blue{
 To prove \autoref{nonperiodicirrationalex:lem} we use the following auxiliary lemma.
 \begin{lemma}\label{nonperiodicirrationalex:lem0}
 For $s\in \set{1,2}$, let $E\bra{s}=\set{\bm{w}\in \Real^2\,:\, (-1)^s \bra{w_2- \sqrt{2}w_1} \geq 0}$. Then
 \[E(s)=\proj_{\bm{w}}\bra{\set{\bra{\bm{w},\by}\in \Real^2\times \Real^2\,:\, \norm{\bra{y_1,y_1}}_2\leq (-1)^s (w_2-y_2),\quad \norm{\bra{y_2,y_2}}_2\leq (-1)^{s+1} 2 (w_1-y_1)}}.\]
 In particular, if $E=\set{\bm{w}\in \Real^2\,:\, w_2- \sqrt{2}w_1 = 0}$, then
  \[E=\proj_{\bm{w}}\bra{\set{\bra{\bm{w},\by}\in \Real^2\times \Real^4\,:\, \begin{alignedat}{3}
  \norm{\bra{y_1,y_1}}_2\leq y_2-w_2,\quad \norm{\bra{y_2,y_2}}_2\leq  2 (w_1-y_1)\\
   \norm{\bra{y_3,y_3}}_2\leq w_2-y_4 ,\quad \norm{\bra{y_4,y_4}}_2\leq2 (y_3-w_1)
  \end{alignedat}}}.\]
 \end{lemma}
 \proof{\textbf{Proof}}

First note that $E(s)=\set{\bm{w}\in \Real^2\,:\, \norm{\bra{w_1,w_1}}_2\leq (-1)^s w_2} +\set{\bm{w}\in \Real^2\,:\, \norm{\bra{w_2,w_2}}_2\leq (-1)^{s+1} 2 w_1} $. Then, by Section 3.3 in \cite{Ben-Tal2001}, we have
\[E(s)=\proj_{\bm{w}}\bra{\set{\bra{\bm{w},\underline{\bm{w}},\overline{\bm{w}}}\in \Real^2\times \Real^2\times \Real^2\,:\, \begin{alignedat}{3}\norm{\bra{\underline{w}_1,\underline{w}_1}}_2\leq (-1)^s \underline{w}_2,\quad \norm{\bra{\overline{w}_2,\overline{w}_2}}_2\leq (-1)^{s+1} 2 \overline{w}_1,\\ w_1=\underline{w}_1+\overline{w}_1,\quad w_2=\underline{w}_2+\overline{w}_2\end{alignedat}}}.\]
The first result follows by letting $y_1=\underline{w}_1$, $y_2=\overline{w}_2$ and eliminating $(\overline{w}_1,\underline{w}_2)$ through the equations $\overline{w}_1=w_1-y_1$ and $\underline{w}_2=w_2-y_2$. The second result follows by noting that $E=E(1)\cap E(2)$.
 \Halmos\endproof
  \begin{lemma}\label{nonperiodicirrationalex:lem}
Let $M=\set{\bra{x,\bz}\in \mathbb{R}\times \Real^2\,:\, z_1=x,\; -0.4\leq \sqrt{2} z_1 - z_2 \leq     0.4}$. Then
\[M=\proj_{x,\bz}\bra{\set{\bra{x,\by,\bz}\in \mathbb{R}\times \Real^5\times \Real^2\,:\,\begin{alignedat}{3} z_1=x,\; -0.4\leq y_5 - z_2 \leq     0.4\\  \norm{\bra{y_1,y_1}}_2\leq y_2-y_5,\quad \norm{\bra{y_2,y_2}}_2\leq  2 (z_1-y_1)\\
   \norm{\bra{y_3,y_3}}_2\leq y_5-y_4 ,\quad \norm{\bra{y_4,y_4}}_2\leq 2 (y_3-z_1)
\end{alignedat}}}.\]
 \end{lemma}
 \proof{\textbf{Proof}}
 Follows by using \autoref{nonperiodicirrationalex:lem0} to model $y_5=\sqrt{2}z_1$.
 \Halmos\endproof
}
\section{Jeroslow restricted binary MICP-R}\label{appjeroslowbinaryequiv}

The following proposition shows that pure binary MICP-R coincides with the restricted version of binary MICP-R that Jeroslow defines in Section 2 of  \cite{jeroslow1987representability} through level sets of homogeneous convex functions.
\begin{restatable}{proposition}{jeroslowbinaryequiv}\label{prop:jeroslowbinaryequiv}
A non-empty set $S\subseteq \Real^n$ is pure binary MICP-R if and only if there exist $m\geq 1$, $\bm{b}\in \Real^m$ and $F:\Real^{n+p'+d'}\to \bra{\Real\cup\set{\infty}}^m$ such that   \begin{equation}\label{jeroslow_condition}
    F\bra{\bm{0},\by,\bz}\leq \bm{0} \quad\Rightarrow\quad  \bra{\by,\bz} = \bm{0},
\end{equation}
$F_i$ is closed, positively homogeneous and convex for each $i\in\sidx{m}$, and
\begin{equation}\label{jeroslow_restriction}
   S=\proj_{\bx}\bra{\set{\bra{\bx,\by, \bz}\in \Real^{n+p'}\times \set{0,1}^{d'}\,:\, F\bra{\bx,\by,\bz}\leq \bm{b}}}.
\end{equation}
\end{restatable}
\proof{\textbf{Proof}}
If $S$ is a non-empty pure binary MICP-R set, then there exists a non-empty closed convex set  $M\subseteq \mathbb{R}^{n+d}$ such that $S=\proj_{\bx}
    \bra{M\cap \bra{\mathbb{R}^{n} \times \mathbb{Z}^d}}$ and $\proj_{\bz}\bra{M\cap \bra{\mathbb{R}^{n} \times \mathbb{Z}^d}}\subseteq\set{0,1}^d$. Without loss of generality, we may assume there exists $\bx^0\in \Real^n$ such that $\bra{\bx^0,\bm{0}}\in M$. Let $\gamma:\Real^{n+p}\to \Real\cup\set{\infty}$ be the \textbf{gauge} of $M-\bra{\bx^0,\bm{0}}$ defined by $\gamma\bra{\bx,\bz}=\inf\set{\lambda\,:\, \bra{\bx,\bz}\in \lambda\bra{M-\bra{\bx^0,\bm{0}}}}$. Then, by Theorem C.1.2.5 in \cite{hiriart-lemarechal-2001},  $\gamma\bra{\bx,\bz}$ is non-negative, closed, positively homogeneous, and convex,  $M=\set{\bra{\bx,\bz}\in \Real^{n+d}\,:\, \gamma\bra{\bx-\bx^0,\bz} \leq 1}$ and $M_\infty=\set{\bra{\bx,\bz}\in \Real^{n+d}\,:\, \gamma\bra{\bx,\bz} \leq 0}$. Let $p'=n$, $d'=d$, $m=1+2n+2d$,  $F:\Real^{n+p'+d'}\to \bra{\Real\cup\set{\infty}}^m$ be such that
    $F\bra{\bx,\by,\bz}=\bra{\gamma\bra{\bx-\by,\bz},\by,-\by,\bz,-\bz}$, and $\bm{b}=\bra{1,\bx^0,-\bx^0,\bm{1},\bm{0}}$. Then  $F_i$ is closed, positively homogeneous and convex for each $i\in\sidx{m}$, and \eqref{jeroslow_condition} holds trivially. Finally, $M\cap \bra{\mathbb{R}^{n} \times [0,1]^d}=\proj_{\bx,\bz}\set{\bra{\bx,\by,\bz}\in \Real^{n+p'+d'}\,:\, F\bra{\bx,\by,\bz} \leq \bm{b}}$, which implies \eqref{jeroslow_restriction} because, under assumption $\proj_{\bz}\bra{M\cap \bra{\mathbb{R}^{n+p} \times \mathbb{Z}^d}}\subseteq\set{0,1}^d$, we have that $M\cap \bra{\mathbb{R}^{n} \times [0,1]^d}$ also induces an MICP formulation of $S$.

  For the converse, it suffices to let $p=p'$, $d=d'$ and show that \[M=\proj_{\bx,\bz}\bra{\set{\bra{\bx,\by, \bz}\in \Real^{n+p}\times [0,1]^d\,:\, F\bra{\bx,\by,\bz}\leq \bm{b}}}\] is a closed convex set. Convexity of $M$ follows from convexity of $F$ and standard properties of the projection operation. The closure results from the following arguments in the proof of Theorem 2.1 in \cite{jeroslow1987representability}. Let $\set{\bra{\bx^m,\by^m,\bz^m}}_{m\in\mathbb{N}}$ be such that $F\bra{\bx^m,\by^m,\bz^m}\leq b$ for all $m\in\mathbb{N}$ and $\bra{\bx^0,\bz^0}=\lim_{m\to\infty}\set{\bra{\bx^m,\bz^m}}_{m\in\mathbb{N}}$. If $\set{\bra{\by^m}}_{m\in\mathbb{N}}$ is unbounded, we may assume without loss of generality that $\by^m\neq \bm{0}$ for all $m\in\mathbb{N}$, $\lim_{m\to\infty}\norm{\by^m}=\infty$ and $\lim_{m\to\infty}\by^m/\norm{\by^m}=\by'$ for some $\by'\in \Real^p$ with $\norm{\by'}=1$. By positive homogeneity, of $F$ we have $F\bra{\bx^m/\norm{\by^m},\by^m/\norm{\by^m},\bz^m/\norm{\by^m}}\leq b/\norm{\by^m}$, and by the closure of $F$ we can take the limit when $m$ goes to $\infty$ to conclude $F\bra{\bm{0},\by',\bm{0}}\leq \bm{0}$, which contradicts \eqref{jeroslow_condition}. Hence, $\set{\bra{\by^m}}_{m\in\mathbb{N}}$ is bounded, and without loss of generality $\bra{\bx^0,\by^0,\bz^0}=\lim_{m\to\infty}\set{\bra{\bx^m,\by^m,\bz^m}}_{m\in\mathbb{N}}$ for some $\by^0\in \Real^n$. The closure of $F$ implies that $F\bra{\bx^0,\by^0,\bz^0}\leq \bm{b}$ so $\bra{\bx^0,\bz^0}\in M$ and $M$ is closed.
\Halmos\endproof

\section{Illustrating the technicalities in the proof of \autoref{periodictheoref}.}\label{defappendix}

The following example illustrates the need for \autoref{nearlyperiodicdef} and the need to consider non-closed sets in \autoref{lem:generaldecompsimple}  and \autoref{peroidictheomoregeneral}.

\begin{example}\label{nearlyperiodicclosureex}
Consider the closed convex set\footnote{See \autoref{Jeroslowexample} for the closure and convexity of this set.}
 \begin{align*}
  M&=\set{\bra{x,y,\bz}\in \Real_+\times \Real_+\times \Real^2_+\,:\,
 (z_1+2z_2)+(z_1+2z_2)^2/y\leq x\leq (z_1+2z_2)+1,\quad z_1\leq 1}\\
  &=\set{\bra{x,y,\bz}\in \Real_+\times \Real_+\times \Real^2_+\,:\,(z_1+2z_2)^2 \leq
 \bra{x-(z_1+2z_2)}y,\quad x\leq (z_1+2z_2)+1,\quad z_1\leq 1}
      \end{align*}
 The index set of the MICP formulation induced by $M$ is $I=[0,1]\times[0,\infty)$, and the $\bz$-projected sets of this formulation are $A_{\bz}=[0,1]$ if $z_1=z_2=0$ and $A_{\bz}=(z_1+2z_2,z_1+2z_2+1]$ otherwise. Then $M$ induces a rational MICP formulation of $S=[0,\infty)=[0,1]\cup\bigcup_{j=1}^\infty (j,j+1]$. The following observations are in order.
 \begin{itemize}
 \item The formulation induced by $M$ is nearly periodic, and $S$ is a closed periodic set.  However, the $\bz$-projected sets of the formulation are not closed, and the formulation does not satisfy \eqref{nearlyperiodiceq} with $\operatorname{cl}\bra{A_{\bz+\lambda \bm{r}}}$ replaced by $A_{\bz+\lambda \bm{r}}$.
 \item We have that $\bra{{I}}_{\mathbb{Z}}$ does not have an integer point in its relative interior, and \autoref{interiorpointind} yields two closed convex sets $\set{M_i}_{i=1}^2$ that satisfy
 \[M_i\cap \bra{\Real^2\times \mathbb{Z}^2}= M\cap \bra{\mathbb{R}^2\times \set{i-1}\times\mathbb{Z}} \quad \forall i\in \sidx{2}. \]
 In particular, $M_2$ induces a formulation of $S_2=\bigcup_{j=0}^\infty (2j+1,2j+2]$, which is not a closed set.
 \end{itemize}
 \Halmos
 \end{example}

\section{Tools from the geometry of numbers}\label{Khintchine:sec}
\fulldimredlemma*
\proof{\textbf{Proof}}
\label{fulldimredlemma:proof}
Let $\bm{A}\in\mathbb{Z}^{d\times k}$ be a rank $k$ matrix and $\bm{c}\in \mathbb{Z}^d\cap F$ be such that $L=\set{\bm{A}\bz:\bz\in \mathbb{R}^k}$ and $F=L+\bm{c}$. By Corollary 4.3b in \cite{schrijver1998theory}, there exist a unimodular matrix $\bm{U}\in \mathbb{Z}^{d\times d}$ and an invertible $\bm{B}\in \mathbb{Z}^{k\times k}$ such that
\begin{equation}\label{hermite}
    \begin{pmatrix}\bm{B}\\\bm{0}_{\bra{d-k}}\end{pmatrix} = \bm{U}\bm{A}.
    \end{equation}Since $\bm{U}L=\set{\bm{U}\bm{A}\bz:\bz\in \mathbb{R}^k}$, \eqref{hermite} and the fact that $B$ is invertible concludes the first bullet point.

For the second bullet point, first notice that from standard linear algebra we have $L^{\perp}=\set{\bz\in \mathbb{R}^d\,:\, \bm{A}^{T}\bz=\bm{0}}$. Now fix any $\bz \in L^{\perp}$, and consider an arbitrary $\bw=(U^T)^{-1}\bz$. The property that $\bw \in (U^T)^{-1}L^{\perp}$ is equivalent to the property that  $\bm{A}^T\bm{U}^Tw=\bm{A}^T\bz=0$. Since $\bm{A}^T \bm{U}^T=\begin{pmatrix}\bm{B}^T\quad \bm{0}_{k\times \bra{d-k}}\end{pmatrix} $ and $\bm{B}^T$ are invertible, we conclude that this is equivalent with the first $k$ coordinates of $\bw$ being zero. This concludes the second bullet point.

For the final statement, note that for $k=1$ we have that \eqref{hermite} is equivalent with $B=B_{1,1}$ and $U$ satisfying $\bm{e}(1)=\bm{U}\bra{\bm{r}/B_{1,1}}$. Using that $U$ is unimodular, this implies $\bm{r}/B_{1,1}\in \mathbb{Z}^d$. Together with $\operatorname{gcd}(r_1,\ldots,r_d) = 1$, we get $B_{1,1}=1$, and therefore $U^{-1}\bm{e}(1)=r$ as we wanted.
\Halmos\endproof

To prove \autoref{latticefreecoro}, we will use a standard combination of Lov\'asz's characterization of maximal lattice-free convex set\footnote{i.e. maximal convex sets without an integer point in their interior} \cite{lovasz1989geometry} and Khintchine's Flatness Theorem (e.g. \cite[Theorem 8.3]{barvinok2002course}). In particular, we combine the  version of Lov\'asz's characterization from \cite{basu2010maximal}  given by
\begin{theorem}\label{latticefreetheoremoriginal}
Let $F\subseteq \Real^d$ be a rational affine space such that $F\cap\mathbb{Z}^d\neq \emptyset$, and let $C$ be a \emph{maximal lattice-free convex set in $F$}. That is, let $C\subseteq F$ be such that
\begin{itemize}
    \item $C$ is convex,
    \item $C$ does not contain an integer point in its interior relative to $F$, and
    \item $F$ is inclusion-wise maximal with respect to the previous three properties and $C\subseteq F$.
\end{itemize}

If $\dim\bra{C}=\dim\bra{F}$, then there exist a  polytope $P\subseteq F$ and a rational linear space $L\subseteq F$ such that $C=P+L$, $\dim\bra{C}=\dim\bra{P}+\dim\bra{L}$, and $C$ contains at least one integer point in its boundary.
\end{theorem}
with the simplification of Theorem 8.3 in Chapter VII of \cite{barvinok2002course} given by
\begin{theorem}\label{flatnesstheo}
Let $C\subseteq \Real^d$ be a full-dimensional compact convex set such that $C\cap\mathbb{Z}^d=\emptyset$. Then there exist $\bm{\pi}\in \mathbb{Z}^d\setminus \bm{0}$ and  $\pi_L,\pi_U\in \mathbb{R}$ such that $C\subseteq \set{\bz\in \Real^d\,:\, \pi_L\leq  \bm{\pi}^T \bz\leq \pi_U}$.
\end{theorem}

\latticefreecorostatement*
\proof{\textbf{Proof}}
\label{latticefreecorostatement:proof}
Let $k=\dim\bra{I_{\mathbb{Z}}}$ and $\bm{c}\in I\cap\mathbb{Z}^d$. By \autoref{fulldimred}, there exists a unimodular  $\bm{U}_1\in \mathbb{Z}^{d\times d}$    such that
\begin{equation}\label{latticefreecoro1}
    \tilde{I}=\bm{U}_1\bra{I_{\mathbb{Z}}-\bm{c}}\subseteq \bm{U}_1\bra{\aff_{\mathbb{Z}}\bra{I}-\bm{c}} = \mathbb{R}^{k}\times \set{\bm{0}}^{d-k}.
\end{equation}
If $I_{\mathbb{Z}}$ does not have an integer point in its relative interior, then $\tilde{I}$ does not have an integer point in its interior relative to $\mathbb{R}^{k}\times \set{\bm{0}}^{d-k}$. Then, there exists a maximal lattice free convex set $C\subseteq \mathbb{R}^{k}\times \set{\bm{0}}^{d-k}$ such that $\tilde{I}\subseteq C$, and by \autoref{latticefreetheoremoriginal} there exists a polytope $\tilde{P}\subseteq \mathbb{R}^{k}\times \set{\bm{0}}^{d-k}$ and a rational linear space $\tilde{L}\subseteq \mathbb{R}^{k}\times \set{\bm{0}}^{d-k}$ such that $C=\tilde{P}+\tilde{L}$ and $k=\dim\bra{C}=\dim\bra{\tilde{P}}+\dim\bra{\tilde{L}}$. Let $m=k-\dim\bra{\tilde{L}}$, and let $\tilde{L}_k$ be the projection of $\tilde{L}$ onto the first $k$ coordinates. Then,  by \autoref{fulldimred}, there exists an unimodular   $\overline{\bm{U}}_2\in \mathbb{Z}^{k\times k}$  such that $\overline{\bm{U}}_2\tilde{L}_k= \set{\bm{0}}^{m} \times \mathbb{R}^{k-m}$. Then $\overline{\bm{U}}_2 C= \hat{P} + \bra{\set{\bm{0}}^{m} \times \mathbb{R}^{k-m}}$ where  $\hat{P}=\overline{\bm{U}}_2 C \cap \bra{\mathbb{R}^{m}\times \set{\bm{0}}^{k-m}}\subseteq \Real^k$ is a polytope such that $\dim\bra{\hat{P}}=m$. Let $P=\hat{P}\times \set{\bm{0}}^{d-k}$ and $\bm{U}_2\in \mathbb{Z}^{d\times d}$ be the block-diagonal matrix whose first diagonal block is equal to $\overline{\bm{U}}_2$ and whose remaining diagonal block is the $(d-k)\times(d-k)$ identity matrix. Then  $\bm{U}_2$ is an unimodular matrix  such that
\begin{equation}\label{latticefreecoro2}
  \bm{U}_2 \tilde{I}\subseteq  P+ \bra{\set{\bm{0}}^{m} \times \mathbb{R}^{k-m}\times \set{\bm{0}}^{d-k}} \text{ and } \bm{U}_2\bra{\mathbb{R}^{k}\times \set{\bm{0}}^{d-k}} = \mathbb{R}^{k}\times \set{\bm{0}}^{d-k}.
\end{equation}

Let $P_m$ be the projection of $P$ (or equivalently $\hat{P}$) onto the first $m$ coordinates. We have that $\dim\bra{P_m}=m$, and because $C$ is lattice free and  $\overline{\bm{U}}_2$ is unimodular, we have that $P_m$ does not contain an integer point in its  interior, but it may contain a point in its boundary ($C$ does contain such points because of its maximality and they may be inherited by $P$). Let $\bar{\bz}\in \mathbb{R}^{m}$ be a point in the interior of $P_m$, $\varepsilon\in (0,1)$ and  $P_\varepsilon= (1-\varepsilon) \bra{P_k-\bar{\bz}}+\bar{\bz}$. Because $P_m$ does not contain an integer point in its  interior, we have that $P_\varepsilon\cap\mathbb{Z}^m=\emptyset$. Then, because of \autoref{flatnesstheo}, there exist  $\bm{\pi}\in \mathbb{Z}^m\setminus \bm{0}$ and  $\bar{\pi}_L,\bar{\pi}_U\in \mathbb{R}$ such that $P_\varepsilon\subseteq \set{\bz\in \mathbb{R}^{m}\,:\, \bar{\pi}_L\leq  {\bm{\pi}}^T \bz\leq\bar{\pi}_U}$, and hence $P_m\subseteq \set{\bz\in \mathbb{R}^{m}\,:\, a\leq {\bm{\pi}}^T \bz\leq b}$ for $a= \left\lfloor \bra{{\bar{\pi}}_L - \varepsilon{\bm{\pi}}^T  \bar{\bz}}/\bra{1-\varepsilon} \right\rfloor$ and $b=\left\lceil \bra{{\bar{\pi}}_U - \varepsilon {\bm{\pi}}^T  \bar{\bz}}/\bra{1-\varepsilon} \right\rceil$.

Finally, by \autoref{fulldimred}, there exists a unimodular matrix  $\overline{\bm{U}}_3\in \mathbb{Z}^{m\times m}$  such that $\overline{\bm{U}}_3{\bm{\pi}}/\operatorname{gcd}({\pi}_1,\ldots,{\pi}_d) = \bm{e}(1)$. Then $\overline{\bm{U}}_3^{-T}P_m\subseteq \set{\bz\in \mathbb{R}^{m}\,:\, a\leq z_1 \leq b}$. Let $\bm{U}_3\in \mathbb{Z}^{d\times d}$ be the block diagonal matrix whose first diagonal block is $\overline{\bm{U}}_3^{-T}$ and whose remaining diagonal block is the $(d-m)\times(d-m)$ identity matrix.  Then  $\bm{U}_3$ is an unimodular matrix  such that
\begin{equation}\label{latticefreecoro3}
\bm{U}_3\bra{P+ \bra{\set{\bm{0}}^{m} \times \mathbb{R}^{k-m}\times \set{\bm{0}}^{d-k}}}\subseteq \set{\bz\in \mathbb{R}^{d}\,:\, a\leq z_1 \leq b} \text{ and } \bm{U}_3\bra{\mathbb{R}^{k}\times \set{\bm{0}}^{d-k}} = \mathbb{R}^{k}\times \set{\bm{0}}^{d-k}.
\end{equation}
The result follows from \eqref{latticefreecoro1}--\eqref{latticefreecoro3} and \autoref{intafflemma} by letting $\mathcal{R}(\bz)=\bm{U}_3\bm{U}_2\bm{U}_1\bra{\bz-\bm{c}}$, $l=a+c_1$ and $u=b+c_1$.
\Halmos\endproof

\section{Recession cones of closed and non-closed convex sets}\label{sec:recessionconeex}

\recessionconeprop*
\proof{\textbf{Proof}}
\label{recessionconeprop:proof}
The fact that  $C_\infty$ is a convex cone containing the origin is shown in Theorem 8.1 of \cite{Rockafellar1997}. The fact that $\bra{\cl\bra{C}}_\infty=\bra{\relint\bra{C}}_\infty$ is shown in Corolarry 8.3.1 in \cite{Rockafellar1997}, and the fact that $C_\infty$ is closed when $C$ is closed is shown in Theorem 8.2 of \cite{Rockafellar1997}. The fact that $C_\infty=C_\infty\bra{\bx}$ for all $\bx \in C$ if $C$ is closed is shown in  Theorem 8.3 of \cite{Rockafellar1997}. Finally, for the containment statement, let $\bm{u}\in C_\infty$ and $\bx\in C$ so that $\bx+\lambda\bm{u}\in C$ for all $\lambda\geq 0$. Then $\bx+\lambda\bm{u}\in C'$ for all $\lambda\geq0$, and hence by Proposition A.2.2.1 of \cite{hiriart-lemarechal-2001} we have $\bm{u}\in C'_\infty$.
\Halmos\endproof

The following example shows  how some common properties for recession cones of closed convex sets may fail to hold for  non-closed convex sets, even for the index set of an MICP formulation, which is the  projection of a closed convex set.
\begin{example}\label{recessionconeexex}
Consider the closed convex set $M=\set{\bra{\bz,y}\in \Real^2\times\Real_+\,:\, z_2^2 \leq z_1 y}$\footnote{See \autoref{Jeroslowexample} for the closure and convexity of this set.}
 and its projection $I=\proj_{\bz}\bra{M}$ onto the $\bz$ variables. We have  $I=\bra{(0,\infty)\times (-\infty,\infty)} \cup \set{\bra{0,0}}$, $I_\infty = I$ and $\cl\bra{I_\infty}=I_{\overline{\infty}}=[0,\infty)\times (-\infty,\infty)$. Then, both $I$ and $I_\infty$ are non-closed, and $I_\infty\bra{(0,0)}=I_\infty\subsetneq I_{\overline{\infty}}=I_\infty\bra{(1,1)}$. Finally, while $(0,\infty)\times (-\infty,\infty)=\relint\bra{I}\subsetneq I$, we have $I_\infty \subsetneq \bra{\relint\bra{I}}_\infty=[0,\infty)\times(-\infty,\infty)$. \Halmos
\end{example}

\section{Rational unboundedness is not preserved by rational affine sections}\label{rationalnotaffine} The following example illustrates that the intersection of a rationally unbounded set with an affine subspace need not remain rationally unbounded.

\begin{lemma}\label{rationalaffinesections}
The closed convex set $I=\set{\bz\in \Real^4\,:\, \bra{z_1+\sqrt{2}z_2}^2\leq z_3,\quad \bra{z_2-\sqrt{2}z_1}^2\leq z_4 }$
is rationally unbounded, but  $\set{\bz\in I\,:\, z_4 =0.4}$ is not rationally unbounded.
\end{lemma}
\proof{\textbf{Proof}}
Let  $\mathcal{R}:\Real^4\to \Real^d$ be an arbitrary rational affine transformation. Without loss of generality (by possibly translating by $-\mathcal{R}\bra{\bm{0}}$) we may assume that $\mathcal{R}\bra{\bm{0}}=\bm{0}$ and that $\mathcal{R}$ is a linear transformation. Under this linearity assumption, we have that  $\mathcal{R}\bra{\bm{e}\bra{3}},\mathcal{R}\bra{\bm{e}\bra{4}}\in \mathcal{R}\bra{I}_\infty\cap \mathbb{Q}^d$, so if $\mathcal{R}\bra{\bm{e}\bra{3}},\mathcal{R}\bra{\bm{e}\bra{4}}\neq \bm{0}$, then there exists $\bm{r}\in \mathbb{Z}^d\cap \bra{\mathcal{R}\bra{I}_\infty\setminus \set{\bm{0}}}$. If $\mathcal{R}\bra{\bm{e}\bra{3}}=\mathcal{R}\bra{\bm{e}\bra{4}}= \bm{0}$,
then $\mathcal{R}\bra{I}=\mathcal{R}\bra{\proj_{z_1,z_2}\bra{I}\times \set{0}^2}=\mathcal{R}\bra{\Real^2\times \set{0}^2}$ and hence $I$ is rationally unbounded because $\Real^2\times \set{0}^2$ is rationally unbounded.

Finally, $\proj_{z_1,z_2}\bra{\set{\bz\in I\,:\, z_4 =0.4}}=\set{\bz\in \Real^3\,:\, \bra{z_1+\sqrt{2}z_2}^2\leq z_3,\quad \bra{z_2-\sqrt{2}z_1}^2\leq 0.4 }$ is the set from \autoref{almostfinalexampleagain}, which is not rationally unbounded. Then $\set{\bz\in I\,:\, z_4 =0.4}$ is not rationally unbounded, because the orthogonal projection is a rational affine transformation.
\Halmos\endproof

\end{APPENDICES}

\end{document}